\documentclass[12pt]{amsart}
\usepackage{amsmath, amsthm, amssymb}
\usepackage{enumerate}

\input{diagrams}
\diagramstyle[scriptlabels,height=8mm,width=8mm]

\def\pdfsyncstart{}
\def\pdfsyncstop{}

\def\bdi{\pdfsyncstop\begin{diagram}}
\def\edi{\end{diagram}\pdfsyncstart}

\theoremstyle{plain}

\newtheorem{thm}{Theorem}[subsection]
\newtheorem{cor}[thm]{Corollary}
\newtheorem{lem}[thm]{Lemma}
\newtheorem{prop}[thm]{Proposition}

\theoremstyle{definition}
\newtheorem{defi}[thm]{Definition}
\newtheorem{defis}[thm]{Definitions}
\newtheorem{conj}[thm]{Problem}
\newtheorem{conv}[thm]{Convention}
\newtheorem{nota}[thm]{Notation}
\newtheorem{rem}[thm]{Remark}
\newtheorem{rems}[thm]{Remarks}
\newtheorem{exa}[thm]{Example}
\newtheorem{exas}[thm]{Examples}
\newtheorem{sit}[thm]{}

\newcommand{\brem}{\begin{rem}}
\newcommand{\brems}{\begin{rems}}
\newcommand{\erem}{\end{rem}}
\newcommand{\erems}{\end{rems}}
\newcommand{\bexa}{\begin{exa}}
\newcommand{\bexas}{\begin{exas}}
\newcommand{\eexa}{\end{exa}}
\newcommand{\eexas}{\end{exas}}
\newcommand{\bdefi}{\begin{defi}}
\newcommand{\edefi}{\end{defi}}
\newcommand{\bdefis}{\begin{defis}}
\newcommand{\edefis}{\end{defis}}
\newcommand{\bcor}{\begin{cor}}
\newcommand{\ecor}{\end{cor}}
\newcommand{\blem}{\begin{lem}}
\newcommand{\elem}{\end{lem}}
\newcommand{\bconv}{\begin{conv}}
\newcommand{\econv}{\end{conv}}
\newcommand{\bconj}{\begin{conj}}
\newcommand{\econj}{\end{conj}}
\newcommand{\bprop}{\begin{prop}}
\newcommand{\eprop}{\end{prop}}
\newcommand{\bthm}{\begin{thm}}
\newcommand{\ethm}{\end{thm}}
\newcommand{\bnota}{\begin{nota}}
\newcommand{\enota}{\end{nota}}
\newcommand{\bsit}{\begin{sit}}
\newcommand{\esit}{\end{sit}}
\newcommand{\be}{\begin{equation}}
\newcommand{\ee}{\end{equation}}
\newcommand{\bproof}{\begin{proof}}
\newcommand{\eproof}{\end{proof}}
\def\ba{\begin{array}}
\def\ea{\end{array}}
\def\bea{\begin{eqnarray}}
\def\eea{\end{eqnarray}}

\def\bnum{\begin{enumerate}}
\def\enum{\end{enumerate}}
\newcommand{\no}{\noindent}
\def\ds{\displaystyle}

\def\lto{\longrightarrow}

\def\ext{{\rm ext}}
\def\sto{\rightsquigarrow}

\newcommand{\Spec}{\operatorname{Spec}}

\newcommand{\Pic}{\operatorname{Pic}}
\newcommand{\hot}{\operatorname{h.o.t.}}

\newcommand{\Frac}{\operatorname{Frac}}

\newcommand{\id}{\operatorname{id}}
\newcommand{\discr}{\operatorname{discr}}

\newcommand{\Ker}{\operatorname{Ker}}

\newcommand{\Aut}{{\operatorname{Aut}}}

\newcommand{\reg}{{\operatorname{reg}}}

\def\supp{{\rm supp\,}}

\renewcommand{\div}{{\operatorname{div}}}

\def\fF{{\mathfrak F}}

\def\fM{{\mathfrak M}}

\def\cM{{\mathcal M}}

\def\cO{{\mathcal O}}

\def\cX{{\mathcal X}}

\def\cY{{\mathcal Y}}

\def\bV{{\bar V}}

\def\tQ{{\tilde Q}}
\def\tO{{\tilde O}}

\def\tV{{\tilde V}}

\newcommand{\la}{\label}

\newcommand{\A}{{\mathbb A}}

\newcommand{\GG}{{\mathbb G}}
\newcommand{\PP}{{\mathbb P}}

\newcommand{\C}{{\mathbb C}}
\newcommand{\Q}{{\mathbb Q}}
\newcommand{\Z}{{\mathbb Z}}
\newcommand{\N}{{\mathbb N}}
\newcommand{\T}{{\mathbb T}}

\newcommand{\G}{{\Gamma}}

\newcommand{\p}{{\partial}}

\def\To{\Rightarrow}

\def\reg{{\rm reg}}

\newcommand{\nlin}{\unitlength1mm\begin{picture}(0,9.25)
                     \put(0,0.75){\line(0,1){8.5}}
                    \end{picture}}

\newcommand{\vlin}[1]{\hspace{0.75mm}\unitlength1mm\begin{picture}
(#1,0)
                     \put(0,0){\line(1,0){#1}}
                    \end{picture}\hspace{0.75mm}\rule[-3mm]{0mm}{4mm}}

\def\llin{\vlin{11.5}}
\newcommand{\lin}{\vlin{8.5}}

\newcommand{\co}[1]{\unitlength1mm\begin{picture}(0,8)
  \put(0,0){\circle{1.5}}
  \put(0,3){\makebox(0,5)[b]{$#1$}}
                    \end{picture}}
\newcommand{\mybox}{\unitlength1mm\begin{picture}(0,1.5)
  \put(-0.75,-0.75){\line(0,1){1.5}}
  \put(-0.75,-0.75){\line(1,0){1.5}}
  \put(0.75,0.75){\line(0,-1){1.5}}
  \put(0.75,0.75){\line(-1,0){1.5}}
  \end{picture}}
\newcommand{\boxo}[1]{\unitlength1mm\begin{picture}(0,8)
  \put(0,0){\mybox}
  \put(0,3){\makebox(0,5)[b]{$#1$}}
                    \end{picture}}
\newcommand{\xbox}{\unitlength1mm\begin{picture}(0,1.5)
  \put(0,0){$\mybox$}
  \put(-0.75,0){\line(1,0){1.5}}
  \put(0,-0.75){\line(0,1){1.5}}
  \end{picture}}

\newcommand{\cou}[2]{\unitlength1mm\begin{picture}(0,8)
  \put(0,0){\circle{1.5}}
  \put(0,3){\makebox(0,5)[b]{$#1$}}
  \put(0,-7){\makebox(0,4)[t]{$#2$}}
    \end{picture}
    \rule[-7mm]{0mm}{7mm}}

\newcommand{\boxou}[2]{\unitlength1mm\begin{picture}(0,8)
  \put(0,0){\mybox}
  \put(0,3){\makebox(0,5)[b]{$#1$}}
  \put(0,-7){\makebox(0,4)[t]{$#2$}}
    \end{picture}
    \rule[-7mm]{0mm}{7mm}}

\newcommand{\xbrl}[2]{\unitlength1mm\begin{picture}(0,8)
  \put(0,0){\xbox}
  \put(-5,0){\makebox(0,5)[b]{$#1$}}
 \put(5,0){\makebox(0,5)[b]{$#2$}}
    \end{picture}
    \rule[-7mm]{0mm}{7mm}}

\newcommand{\boxrl}[2]{\unitlength1mm\begin{picture}(0,8)
  \put(0,0){\mybox}
  \put(-5,0){\makebox(0,5)[b]{$#1$}}
 \put(5,0){\makebox(0,5)[b]{$#2$}}
    \end{picture}
    \rule[-7mm]{0mm}{7mm}}

\newcommand{\xbshiftup}[2]{\unitlength1mm\begin{picture}(0,9.25)
                     \put(0,10){\xbrl{#1}{#2}}
                    \end{picture}}

\newcommand{\boxshiftup}[2]{\unitlength1mm\begin{picture}(0,9.25)
                     \put(0,10){\boxrl{#1}{#2}}
                    \end{picture}}

\title{Smooth affine
surfaces with non-unique $\C^*$-actions}

\author{Hubert Flenner}
\address{Fakult\"at f\"ur Mathematik,
Ruhr Universit\"at Bochum, Geb.\ NA 2/72, Universit\"ats\-str.\
150, 44780 Bochum, Germany}
\email{Hubert.Flenner@rub.de}

\author{Shulim Kaliman}
\address{Department of Mathematics,
University of Miami, Coral Gables, FL  33124, U.S.A.}
\email{kaliman@math.miami.edu}

\author{Mikhail Zaidenberg}
\address{Universit\'e
Grenoble I, Institut Fourier, UMR 5582 CNRS-UJF, BP 74, 38402
St.\ Martin d'H\`eres c\'edex, France}
\email{zaidenbe@ujf-grenoble.fr}

\thanks{
{\bf Acknowledgements:} This research was done during a visit of
the first and the second authors at the Institut Fourier,
Grenoble and of all three authors at
the Max-Planck-Institute of Mathematics, Bonn.
They thank these institutions for the
generous support and excellent working conditions.}

\thanks{
\mbox{\hspace{11pt}}{\it 1991 Mathematics Subject
Classification}:
14R05, 14R20, 14J50.\\
\mbox{\hspace{11pt}}{\it Key words}: $\C^*$-action, $\C_+$-action,
affine surface}

\begin{document}

\begin{abstract}
In this paper we complete
the classification of effective $\C^*$-actions on smooth affine
surfaces up to conjugation in the full automorphism group and up
to inversion $\lambda\mapsto \lambda^{-1}$ of $\C^*$. If a smooth
affine surface $V$ admits more than one $\C^*$-action then it is
known to be Gizatullin i.e.,  it can be completed by a
linear chain of smooth rational curves. In \cite{FKZ3} we gave a
sufficient condition, in terms of the Dolgachev-Pinkham-Demazure
(or DPD) presentation, for the uniqueness of a $\C^*$-action on a
Gizatullin surface. In the present paper we show that this
condition is also necessary, at least in the smooth case. In fact,
if the uniqueness fails for a smooth Gizatullin surface $V$ which
is neither toric nor Danilov-Gizatullin, then $V$ admits a
continuous family of pairwise non-conjugated $\C^*$-actions
depending on one or two parameters. We give an explicit description of
all such surfaces and their $\C^*$-actions in terms of DPD
presentations. We also show that for every $k>0$ one can find a
Danilov-Gizatullin surface $V(n)$ of index $n=n(k)$ with
a family of pairwise non-conjugate
$\C_+$-actions depending on $k$ parameters.
\end{abstract}

\maketitle

\tableofcontents

\section{Introduction}
The classification of $\C^*$-actions on normal affine surfaces up
to equivariant isomorphism is a widely studied subject and by now
well understood, see e.g., \cite{FlZa1}.
However from this classification it is not clear which of these
surfaces are {\em abstractly} isomorphic. This leads to the
question of classifying all equivalence classes of
$\C^*$-actions on a given surface $V$ under the equivalence
relation generated by conjugation in $\Aut(V)$ and inversion
$t\mapsto t^{-1}$ of $\C^*$. In this paper we give a complete
solution to this problem for smooth affine surfaces. In particular
we obtain the following result.

\bthm \label{0.1} Let  $V$ be a smooth affine $\C^*$-surface. Then
its $\C^*$-action is unique up to equivalence if and only if $V$
does not belong to one of the following classes. \bnum[(1)] \item
$V$ is a toric surface;

\item $V=V(n)$ is a Danilov-Gizatullin surface of index $n\ge 4$
(see below);

\item $V$ is a special smooth  Gizatullin surface of type I or II
(see Definition \ref{0.2} below). \enum Furthermore, $V$ admits at
most two conjugacy classes of $\A^1$-fibrations $V\to\A^1$
if and only if  $V$ is not one of the surfaces in (2) or (3).
\ethm

\bsit\la{1000.000}
Here two $\A^1$-fibrations $\varphi_1,\varphi_2:V\to\A^1$ are
called conjugated if $\varphi_2=\beta\circ\varphi_1\circ\alpha$
for some $\alpha\in\Aut (V)$ and $\beta\in\Aut(\A^1)$. Let us
describe in more detail
the exceptions (1)-(3) in Theorem \ref{0.1}.

Obviously uniqueness of $\C^*$-actions fails for affine toric
surfaces. Restricting the torus action to one-dimensional subtori
yields an infinite number of equivalence classes of $\C^*$-action
on such a surface.

A surface $V=V(n)$ as in (2) is by definition the complement to an
ample section $C$ of a Hirzebruch surface $\Sigma_k\to\PP^1$ with
\footnote{Our enumeration of the Danilov-Gizatullin surfaces
differs from that in \cite{FKZ2}. In this enumeration e.g.,
$V_2\simeq \PP^1\times\PP^1\backslash\Delta$, where $\Delta$ is
the diagonal.} $C^2=n$. By a remarkable theorem
of Danilov and
Gizatullin \cite[Theorem II.5.8.1]{DaGi} \footnote{See
Corollary 4.8 in \cite{CNR}, \cite{FKZ4},
or Corollary \ref{DG-uniproof} below
for alternative proofs.} this surface only depends on $n$ and
neither on $k$ nor on the choice of the section. As was observed
by Peter Russell, for $n\ge 4$ there are
non-equivalent
$\C^*$-actions on $V(n)$; see \cite[5.3]{FKZ2}
for a full classification of them.
It is also shown in {\it loc.cit.} that there are at least
$n-1$ non-conjugated $\A^1$-fibrations $V(n)\to\A^1$.
In Section 6.3 we give a complete classification of all
$\A^1$-fibrations on $V(n)$. It turns out that there are even
families of pairwise non-conjugated $\A^1$-fibrations
depending on an arbitrary number of parameters,
if $n$ is sufficiently large;
see  Corollary \ref{7000.010}.
\esit

\bsit\la{1000.001} To describe the special surfaces as in (3)
we recall\footnote{See \cite{Gi} or \cite{FKZ3}.} that
a normal affine surface is said to be {\em
Gizatullin} if it can be completed by a zigzag that is, by
a linear chain of smooth rational curves.
Any normal non-Gizatullin surface admits at most one $\C^*$-action
up to equivalence (see \cite{Be} for the smooth and
\cite{FlZa2} for the general case). Thus, if for a normal
affine $\C^*$-surface uniqueness of $\C^*$-action fails, it must
be Gizatullin.

By a result of Gizatullin \cite{Gi} (see also \cite{Du, FKZ2}) any
non-toric Gizatullin surface has a completion $(\bV, D)$ with a
boundary $D=\bV\backslash V$ which is a {\em standard} zigzag.
This means that $D$ is a zigzag with dual graph \be \G_D:\quad\quad
\cou{0}{C_0}\lin\cou{0}{C_1}\lin\cou{w_2}{C_2}
\lin\ldots\lin\cou{w_n}{C_n}\quad, \ee
where $C_0^2=C_1^2=0$  and
$w_i=C_i^2\le-2,\,\,i=2,\ldots,n$.
%
Such a completion $\bV$ can be constructed starting from the
quadric $Q=\PP^1\times\PP^1$ and the curves
$$
C_0=\{\infty\}\times\PP^1, \,\quad
C_1=\PP^1\times \{\infty\},\,\quad
C_2=\{0\}\times\PP^1\,
$$
by a sequence of blowups on $C_2\backslash C_1$ and
infinitesimally  near points, see \ref{1.2}. An exceptional curve
not belonging to the zigzag is called a {\em feather} \cite{FKZ2}.
If in this process a feather $F$ is created by a blowup on
component $C_\mu$ of the zigzag then we call $C_\mu$ the {\em
mother component} of $F$.\esit

\bdefi\label{0.2} A smooth Gizatullin surface $V$ is said to be
{\em special} if it admits a standard completion $(\bV,D)$ such that

\bigskip

(a) every component  $C_{i+1}$, $i\ge 2$, is created by a blowup
on $C_i\backslash C_{i-1}$, and

(b) $n\ge 3$ and the divisor formed by the feathers  can be written as
$$
F_2+F_{t1}+\cdots+ F_{tr}+F_n, \quad r\ge 0,
$$
where $F_2$, $F_{t\rho}$, $1\le \rho\le r$,  and $F_n$
have mother component $C_2$, $C_t$ and $C_n$, respectively.
Such a surface  $V$ is called a special surface of
\begin{itemize}
\item type I if either $r=1$ or $r\ge 2$ and  $t\in\{2,n\}$;
\item type II if $r\ge 2$ and $2<t<n$.
\end{itemize}
\edefi

The remaining special surfaces with $r=0$ are just the
Danilov-Gizatullin surfaces $V_n$ with $n\ge 3$. They are neither
of type I nor II.

In \cite{FlZa3, FKZ3} we have shown that a $\C^*$-action on a
smooth affine surface $V$ is unique up to equivalence unless $V$
belongs to one of the classes (1)-(3) in Theorem \ref{0.1} (see
also \ref{stcond} and \ref{crit}). A similar uniqueness theorem
holds for $\A^1$-fibrations $V\to\A^1$, see \cite[5.13]{FKZ3} and
Proposition \ref{stren} in Sect.\ 6.4. Thus Theorem \ref{0.1} is a
consequence of the following result.

\bthm\label{0.3} Let $V$ be a special smooth Gizatullin surface.
If $V$ is of type I then the equivalence classes of $\C^*$-actions
on $V$ form in a natural way a $1$-parameter family, while in case
of type II they form a $2$-parameter family. Similarly,
conjugacy classes of $\A^1$-fibrations $V\to\A^1$ contain, in the
case of type I special surfaces, a one-parameter family,
while in case of type II they contain two-parameter families.
\ethm

In particular any special Gizatullin surface admits a
$\C^*$-action. To construct a one- or two-parameter family of such
$\C^*$-actions let us recall \cite{FlZa2} that on a non-toric
Gizatullin surface there can exist only hyperbolic $\C^*$-actions. These
actions can be described via the following DPD-presentation
\cite{FlZa1}.

\bsit\label{eph} Any  hyperbolic $\C^*$-surface can be presented as
$$
V=\Spec A,\qquad\mbox{where}\quad
A=A_0[D_+,D_-]=A_0[D_+]\oplus_{A_0} A_0[D_-]\,
$$
for a pair  of $\Q$-divisors $(D_+,D_-)$
on a smooth affine curve $C=\Spec A_0$ satisfying the condition
$D_++D_-\le 0$. Here
$$
A_0[D_\pm]=\bigoplus_{k\ge 0} H^0(C,
\cO_C(\lfloor kD_\pm\rfloor))u^{\pm k}
\subseteq \Frac(A_0)[u,u^{-1}]\,,
$$
where $\lfloor D\rfloor$ stands for the integral part of a divisor
$D$ and $u$ is an independent variable. Two pairs $(D_+,D_-)$ and
$(D'_+,D'_-)$ are said to be {\em equivalent} if $D'_\pm=D_\pm\pm
\div\,\varphi$ for a rational function $\varphi$ on $C$.
\esit

In terms of these DPD-presentations we can reformulate Theorem \ref{0.3}
in the following more precise form.

\bthm\label{0.5} (a) A smooth Gizatullin surface $V$ equipped
with a hyperbolic $\C^*$-action is special if and only if $V$
admits a DPD-presentation $V=\Spec \C[z][D_+, D_-]$ with
\be\label{010} \left(D_+,D_-\right)=
\left(-\frac{1}{t-1}[p_+]\,\,, -\frac{1}{n-t+1}[p_-]-D_0\right)
\,,\ee where $p_+\ne p_-$, $2\le t\le n$, $n\ge 3$ and
$D_0=\sum_{i=1}^r [p_i]$ is a reduced divisor on $C\cong \A^1$
such that $r\ge 0$ \footnote{For $r=0$ we obtain the
Danilov-Gizatullin surfaces.} and all points $p_i$ are
different from $p_\pm$.

(b) \cite[Theorem
4.3(b)]{FlZa1} Two such $\C^*$-surfaces $V$, $V'$
given by pairs of divisors
$(D_+,D_-)$ and $(D'_+, D'_-)$ as in (a) are equivariantly
isomorphic if and only if $(D_+,D_-)$ and $(D'_+, D'_-)$ are
equivalent up to interchanging $D_+$ and $D_-$,
if necessary, and  up to
an automorphism of the underlying curve $C=\A^1$.

(c) Two surfaces $V$, $V'$  as in (b) are (abstractly) isomorphic
if and only if the unordered pairs
\footnote{ $\{D\}$ denotes the fractional part of the $\Q$-divisor $D$.}
$(\deg\{D_+\},\deg\{D_-\})$ and
$(\deg\{D'_+\},\deg\{D'_-\})$ coincide
and the integral part $\lfloor -D_+- D_-\rfloor$
is equivalent to $\lfloor -D'_+ - D'_-\rfloor$ up to an
automorphism of the underlying curve $C=\A^1$. \ethm

In case $r=0$, (c) is just the theorem of Danilov and Gizatullin
cited above. The generalized Isomorphism Theorem of (c)
is our principal new result, which
occupies the major part of the paper.

The standard boundary zigzag of a special smooth
Gizatullin surface $V$ as in (a) is
\be \label{00w} [[0,0,
(-2)_{t-2}, -2-r, (-2)_{n-t}]]\,,
\ee
where $[[(-2)_{a}]]$
represents a chain of $(-2)$-curves of length $a$. Here the
numbers  $n,r,t$ have the same meaning as in Definition
\ref{0.2}. In particular,
\begin{itemize}
\item $V$ is of type I if and only if either $r=1$ or $r\ge 2$ and one
of the coefficients of $p_\pm$ in (\ref{010}) equals $-1$,
\item $V$ is of type II if and only if $r\ge 2$ and both coefficients
of $p_\pm$ are in the interval $]-1,0[$.
\end{itemize}
Comparing part (b) and  (c) of the theorem, the position of $p_+$
is essential for the equivariant isomorphism type of $V$ while it
does not affect the abstract isomorphism type of $V$ unless the
coefficient of $p_+$ in $D_+$ is integral.  Dually the same holds
for $p_-$. Thus Theorem \ref{0.5} implies Theorem \ref{0.3}, since
fixing $D_0$ we can vary $p_+$ and $p_-$ for surfaces of type II
and one of these points for surfaces of type I.  This preserves
the isomorphism type of $V$, but changes the equivalence class of
the $\C^*$-action and of the $\A^1$-fibration.

\bsit\label{0099} To make the result above more concrete, consider
for instance  the normalization $V$ of the singular surface in
$\A^3$ given by the equation
$$
x^{n-1}y=(z-p_-)(z-1)^{n-1}q^{n-1}(z),\quad\mbox{where}\quad
q(z)=\prod_{i=1}^{r} (z-p_i),\quad r\ge 0,\,\, n\ge 3,
$$
and $1,p_-,p_1\ldots, p_r\in\C$ are pairwise different. Such a surface
carries a hyperbolic $\C^*$-action
$$
\lambda . (x,y,z)=(\lambda x, \lambda^{1-n} y, z)\,,
$$
which amounts to a DPD-presentation (\ref{010}) with $t=2$ and
$p_+:=1$; see Example 4.10 in \cite{FlZa1}. For $r\ge 1$ these
are special smooth Gizatullin surfaces of type I whereas for
$r=0$, $V=V(n)$ is the Danilov-Gizatullin surface of index $n$.

According to Theorem \ref{0.5}(b),(c) fixing a polynomial
$q\in\C[z]$ of degree $r\ge 2$ and varying  $p_-\in\A^1\setminus\{
p_1,\ldots,p_{r+1}\}$ the resulting surfaces $V$ are all
abstractly but not equivariantly
isomorphic, in general. \esit

Let us give a brief overview of the contents
of the various sections.
After recalling in Section 2 some necessary preliminaries, we
prove in Section 3 that the configuration of points in
$\lfloor-D_+-D_-\rfloor$ represents an invariant of the (abstract)
isomorphism type of a Gizatullin $\C^*$-surface $V$, see Corollary
\ref{specialconfig} and Remark \ref{nonspecial}. This remains valid more
generally without assuming the existence of a $\C^*$-action, see
Theorem \ref{matchmain}.

In Section 6 we establish that for special $\C^*$-surfaces this
 configuration of points
together with the numbers $\deg\{D_\pm\}$ is the only invariant of
the isomorphism type of $V$, see Theorem \ref{main2}. This yields
Theorem \ref{0.5}. The proof proceeds in two steps. In
Section 4 we deal with special Gizatullin surfaces of $(-1)$-type,
which are characterized by the property, that all feathers are
$(-1)$-curves. The hard part in the proof of the general case is
to transform any standard completion of a special surface into one
of $(-1)$-type. This problem is solved in Sections 5 and 6.1
using an explicit coordinate description of special surfaces.

Sections 6.2-6.4 contain applications of the main results.
For instance, we show in Section 6.3 how to classify
$\C_+$-actions or, equivalently, $\A^1$-fibrations on Gizatullin
surfaces, see Theorem \ref{fibration.140}.
In Section 6.4 we strengthen our
previous uniqueness result \cite{FlZa3} for $\A^1$-fibrations on
(singular, in general) Gizatullin surfaces.

The authors are grateful to Peter Russell for
inspiring discussions concerning the results of Section 6.3.

\section{Preliminaries}
In this section we recall some necessary notions and facts from
\cite{DaGi} and  \cite{FKZ3}.

\subsection{Standard zigzags and reversions}\label{ss1.1}

\bsit\label{zigz} Let $X$ be a complete normal algebraic surface,
and let $D$ be an SNC (i.e.\ a simple normal crossing) divisor $D$
with rational components contained in the smooth part $X_\reg$ of
$X$. We say that $D$ is a {\it zigzag} if the dual graph $\G_D$ of
$D$ is linear i.e.,
\be
 \G_D:\quad\quad
\cou{w_0}{C_0}\lin\cou{w_1}{C_1}\lin\ldots\lin\cou{w_n}{C_n}\quad,
\ee where $w_i=C_i^2,\,\,i=0,\ldots,n,$ are the weights of $\G_D$.
We abbreviate this chain by $[[w_0,\ldots, w_n]]$. We also write
$[[\ldots, (w)_k,\ldots ]]$ if a weight $w$ occurs at $k$
consecutive places. \esit

\bsit\label{stzi} A zigzag $D$ is called {\em standard}  if its
dual graph $\G_D$ is either $[[w_0,\ldots, w_n]]$ with all $w_i\le
-2$, or $\G_D$ is one of the chains\footnote{ The case of the
zigzag $[[w_0,\ldots, w_n]]$ with all $w_i\le -2$ was
unfortunately forgotten in \cite[Lemma 2.17]{FKZ1}, in
\cite[2.8]{FKZ2} and in \cite[1.2]{FKZ3}. However, such a chain
and also $[[0]]$ cannot appear as boundaries of affine surfaces.
Hence this omission does not affect any of the results of these
papers. By abuse of notation, we often denote an SNC divisor and
its dual graph by the same letter.} \be \label{standardzigzag}
[[0]]\,,\quad [[0,0]]\,, \quad [[0,0,0]] \,\,\mbox{or}\,\,
  [[0,0,w_2,\ldots,w_n]],\,\, \mbox{where}
\,\, n\ge 2 \,\, \mbox{and}\,\, w_j\le-2\,\,\,\forall j.
\ee
A linear chain $\Gamma$ is said to be {\em semi-standard} if it is
either standard or one of
\be \label{sstandard}
 [[0,w_1,w_2,\ldots,w_n]],\quad [[0,w_1,0]]\,, \,\, \mbox{where}
\,\, n\ge 1,\,\,w_1\in \Z,\,\,\mbox{and}\,\,
w_j\le-2\,\,\,\forall j\ge 2\,. \ee
\esit

\bsit\label{reversion}
Every Gizatullin surface $V$ admits a  {\em standard
completion} $(\bV, D)$ i.e., a completion by a standard zigzag
$D$, see \cite{DaGi, Du} or
Theorem 2.15 in \cite{FKZ1}. Similarly we call $(\bV, D)$ a {\em semi-standard
completion} if its boundary zigzag is semi-standard.

The standard boundary  zigzag is unique up to
reversion
\be\la{reverse} D=[[0,0,w_2,\ldots,w_n]]\rightsquigarrow
[[0,0,w_n,\ldots,w_2]]=:D^\vee\,. \ee We say that $D$ is {\em
symmetric} if $D=D^\vee$.

The reversion of a zigzag, regarded as a birational transformation
of the weighted dual graph, admits a factorization into a sequence
of inner elementary transformations (see 1.4 in \cite{FKZ3} and
\cite{FKZ1}). By an {\em inner elementary transformation} of a
weighted graph we mean blowing up at an edge incident to a
$0$-vertex of degree $2$ and blowing down the image of this
vertex. Given $[[0,0,w_2,\ldots,w_n]]$ we can successively move
the pair of zeros to the right
$$
[[0,0,w_2,\ldots,w_n]]\sto [[w_2, 0,0,w_3,\ldots,w_n]]
\sto \ldots \sto [[w_2,\ldots,w_n, 0,0]]
$$
by inner elementary transformations, which gives the reversion. An
{\em outer elementary transformation} consists in blowing up at a
$0$-vertex of degree $\le 1$ and blowing down the image of this
vertex. A birational inner elementary transformation on a surface
is rigid, whereas an outer one depends on the choice of the center
of blowup.

If $(\bV,D)$ is a standard completion of a Gizatullin surface $V$,
then reversing the zigzag $D$ by a sequence of inner elementary
transformations we obtain from $(\bV, D)$ a new completion
$(\bV^\vee, D^\vee)$, which we call the {\em reverse standard
completion}. It is uniquely determined by $(\bV, D)$. \esit

\subsection{Symmetric reconstructions}\label{ss1.2}
Given a Gizatullin surface, any two SNC completions are related
via a birational transformation which is called a  {\em
reconstruction}, see \cite[Definition 4.1]{FKZ3}. For further use
we recall the necessary notions and facts.

\bsit\label{reco} Given weighted graphs
$\Gamma$ and $\Gamma'$, a
{\em reconstruction} $\gamma$  of $\Gamma$ into $\Gamma'$ consists in a
sequence
\bdi \,\qquad\gamma:\quad
\Gamma=\Gamma_0&\rDotsto^{\gamma_1}& \Gamma_1 & \rDots^{\gamma_2}
&\cdots & \rDotsto^{\gamma_n}& \Gamma_n=\Gamma'\,, \edi where each
arrow $\gamma_i$ is either a blowup or a blowdown. The graph
$\Gamma'$ is called the {\em end graph} of $\gamma$. The inverse
sequence $\gamma^{-1}= (\gamma_n^{-1}, \ldots, \gamma_1^{-1})$
yields a reconstruction of $\Gamma'$ with end graph $\Gamma$.

Such a reconstruction is said to be \bnum[-] \item {\em
admissible} if it involves only blowdowns of at most linear
vertices, inner blowups or outer blowups at end vertices;

\item {\em symmetric} if it can be written in the form $(\gamma,
\gamma^{-1})$. Clearly in this case the end graph is again
$\Gamma$. \enum\esit

For symmetric reconstructions the following hold (see Proposition
4.6 in \cite{FKZ3}).

\blem\label{equivariant.6}

Given two standard
completions $(X,D)$ and $(Y,E)$ of a normal Gizatullin surface
$V\not\cong \A^1\times \A^1_*$, \footnote{Here
$\A^1_*:=\A^1\setminus\{0\}$. } after replacing, if necessary,
$(X,D)$ by its reversion $(X^\vee, D^\vee)$, there exists a
symmetric reconstruction of $(X,D)$ into $(Y,E)$.

\elem

\subsection{Generalized reversions}\label{ss1.3}

\bdefi\label{1.4} We say that two standard completions $(\bV,D)$,
$(\bV',D')$ of a Gizatullin surface $V$ are {\em evenly linked} if
there is a symmetric reconstruction of  $(\bV,D)$ into
$(\bV',D')$. In particular, the dual graphs of $D$ and $D'$ are
then the same. Otherwise $(\bV,D)$ and $(\bV',D')$ are called {\em
oddly linked}.

By Lemma \ref{equivariant.6}, $(\bV',D')$ is always evenly
linked to one of the completions $(\bV,D)$ or $(\bV^\vee,
D^\vee)$. \edefi

%

\bdefi\label{1.5} We let $(\bV,D)$ be a semi-standard completion a
Gizatullin surface $V$ with boundary zigzag $[[0,-m,w_2,\ldots,
w_n]]$, where $m\ge 0$. Moving the zero vertex to the right by
elementary transformations we can transform this to a
semi-standard zigzag $[[w_2, \ldots, w_m, -k,0]]$ for every $k\ge
0$. We call the resulting semi-standard completion $(\bV',D')$ a
{\em generalized reversion} of $(\bV,D)$.

Transforming $[[0,-m,w_2,\ldots, w_n]]$ into the standard zigzag
$[[0,0,w_2,\ldots, w_n]]$ requires outer elementary
transformations; see \ref{reversion}. To transform
further
the latter zigzag into $[[w_2, \ldots, w_m, -k,0]]$ only inner
elementary transformations are needed. Thus the resulting
semi-standard completion $(\bV',D')$ depends on parameters, namely
on the choice of the centers of outer blowups when passing from
$[[0,-m,w_2,\ldots, w_n]]$ to $[[0,0,w_2,\ldots, w_n]]$. \edefi

The following proposition follows from  the connectedness part of
Theorem I.1.2 in \cite{DaGi}. We provide an independent proof
relying on \cite{FKZ1}.

\bprop\label{1.6} For any two semi-standard completions $(\bV,D)$,
$(\bV',D')$ of a Gizatullin surface $V$, $(\bV', D')$ can be
obtained from $(\bV, D)$ by a sequence of generalized reversions
$$
(\bV,D)=(\bV_0, D_0)\quad\leadsto\quad (\bV_1,D_1)\quad\leadsto
\quad\ldots\quad\leadsto
\quad (\bV_l,D_l)=(\bV',D')\,.
$$
\eprop

\bproof By Lemma 3.29 in \cite{FKZ1} we can find a sequence of
semi-standard completions
$$
(\bV,D)=(\bV_0, D_0),\quad (\bV_1,D_1), \quad\ldots,
\quad (\bV_l,D_l)=(\bV',D')
$$
such that every step $(\bV_i, D_i)\leadsto  (\bV_{i+1}, D_{i+1})$
is dominated by a completion  $(W_i, E_i)$ of $V$ with a linear
zigzag $E_i$. Thus it is sufficient to show the assertion in case
where $l=1$ i.e., $(\bV,D)$ and $(\bV',D')$ are dominated by a
completion $(W,E)$ with a linear zigzag $E$. We can perform
elementary transformations of $(\bV,D)$ and $(\bV',D')$ to obtain
standard completions $(\bV_0,D_0)$ and $(\bV_0', D_0')$,
respectively, such that all these surfaces are dominated by a
suitable admissible blowup of $(W,E)$. Replacing $(\bV,D)$ and
$(\bV',D')$ by $(\bV_0,D_0)$ and $(\bV_0', D_0')$, respectively,
we are reduced to the case where both $(\bV,D)$ and $(\bV',D')$
are standard completions of $V$. The result follows now from
Proposition 3.4 in \cite{FKZ1}, which says that a birational
transformation between standard graphs $\bdi \G_D&\rDotsto&
\G_{D'}\edi$ dominated by a linear graph is either the identity or
the reversion. \eproof

\section{The principle of matching feathers}
Consider a Gizatullin surface  $V$. By Gizatullin's Theorem
\cite{Gi} (see also \cite{FKZ1}), the sequence of weights
$[[w_2,\ldots,w_n]]$ (up to reversion) of a standard boundary
zigzag $D$ of $V$ is a discrete invariant of the abstract isomorphism type
of $V$. In this section we introduce a more subtle continuous
invariant of $V$ called the {\it configuration invariant}. This is
a point in the product of certain configuration spaces, up to
reversing the order in the product. Although it is defined using a
standard completion of $V$, in Corollary \ref{190} below we
establish that this point is an invariant of the open surface $V$.

\subsection{Configuration spaces}
The configuration invariant takes values in configurations
of points on $\A^1$ and $\A^1_*=\A^1\setminus\{0\}$. We
recall shortly the necessary notions.

\bsit\label{conf+} We
let $\cM^+_s$ denote the configuration space of all $s$-points
subsets $\{\lambda_{1},\ldots,\lambda_{s}\}$ of the affine line
$\A^1$. This is a Zariski open subset of the Hilbert scheme of
$\A^1$. By the main theorem on symmetric functions $\cM^+_s$ can
be identified with the set of all monic polynomials
$P=X^s+\sum_{j=1}^{s} a_jX^{s-j}$ of degree $s$, whose
discriminant is nonzero. This identification
$$\{\lambda_{1},\ldots,\lambda_{s}\}\mapsto P=\prod_{j=1}^s (X-\lambda_j)\,$$
sends
$\cM^+_s$ onto the principal Zariski open subset
$D(\discr):=\A^{s}\backslash\{\discr (P)=0\}$ of $\A^{s}$.

The affine group $\Aut(\A^1)$ acts on $\cM^+_s$ in a natural way.
By restriction we obtain an action on $ \cM^+_s$ of the normal
subgroup $\GG_a$ of translations. The quotient $ \cM^+_s/\GG_a\, $
can be identified with the space, say, $U_0$ of all monic
polynomials $P=X^s+\sum_{j=2}^{s} a_jX^{s-j}$ with $a_1=0$ and
with nonzero discriminant. The residual action of the
multiplicative group
$\C^*\simeq\Aut(\A^1)/\GG_a$ on $U_0$ is given by
$$
P=X^s+\sum_{j=2}^s a_j X^{s-j}\mapsto t.P:=X^s+\sum_{j=2}^s t^ja_j
X^{s-j}\,, \qquad t\in \C^*\,.
$$
Consequently, the quotient
$$
\fM^+_s=\cM^+_s/\Aut(\A^1)
$$
exists and is an affine variety of dimension $s-2$. More
precisely, we can identify $\fM^+_s$ with the principal Zariski
open subset of the weighted projective space
$\PP(2,3,\ldots,s)$ given by the discriminant i.e.,
$\fM^+_s=D_+(\mbox{discr})$, see e.g., \cite{Li} or \cite[Ch. 1,
Example 2]{ZaLi}. \esit

\bsit\label{conf*} Let $\cM_s^*$ be the part of $\cM_s^+$
consisting of all subsets of $\A^1_*$. Similarly as before this
can be identified with the space of all monic polynomials $P$ of
degree $s$ with $P(0)\ne 0$. The group $\C^*$ acts on $\cM_s^*$,
and the quotient $\fM_s^*$ embeds as a principal Zariski open
subset into the weighted projective space $\PP(1,\ldots, s)$.
\esit

\subsection{The configuration invariant}
\label{ss2.2} We consider a smooth Gizatullin surface
$V=\bV\setminus D$ completed by a semi-standard zigzag
\be\label{zigzag} D: \qquad\quad\cou{0}{C_0}\llin
\cou{-m}{C_1}\lin\cou{w_2}{C_2} \llin \ldots\lin\cou{w_i}{C_i}
\llin\ldots\llin \cou{w_n}{C_n} \quad,\quad w_i\le -2
\,\,\,\forall i\ge 2. \ee We associate
to $(\bV, D)$ a point in the product of configuration spaces
$$
\fM:=\fM^{\tau_2}_{s_2}\times\ldots\times \fM^{\tau_n}_{s_n}
$$
for suitable numbers $s_i$, $2\le i\le n$, where $\tau_i\in \{+,*\}$
depends on the component $C_i$ as described below. This point occurs to be
an invariant of $V$, i.e.\ it depends only on the isomorphism
class of $V$ and not on the choice of a semi-standard
completion $(\bV, D)$.

\bsit\label{1.2} To define this invariant we need to recall the
notion of extended divisor. Let $(\bV,D)$ be a semi-standard
completion of a Gizatullin surface $V$. Then the linear system
$|C_0|$ on $\bV$ defines a morphism $\Phi_0:\bV\to \PP^1$ with at
most one degenerate fiber, say, $\Phi_0^{-1}(0)$ while the fiber
$C_0=\Phi^{-1}(\infty)$ is non-degenerate. The reduced SNC
divisor $D_{\ext}=D\cup \Phi_0^{-1}(0)$ is called the {\it
extended divisor} of the completion $(\bV, D)$. By
Proposition 1.11 of \cite{FKZ3}, this divisor has dual
graph \vskip 3mm

\be\label{extended} D_\ext: \qquad\quad\cou{0}{C_0}\llin
\cou{-m}{C_1}\lin\cou{}{C_2} \nlin\xbshiftup{}{\,\fF_2} \llin
\ldots\lin\cou{}{C_i}\nlin\xbshiftup{}{\,\fF_i} \llin\ldots\llin
\cou{}{C_n} \nlin\xbshiftup{}{\,\fF_n} \quad,\qquad\qquad \,\ee
where $\fF_i=\{F_{i\rho}\}_{1\le \rho\le r_i}$ ($2\le i\le n$) is
a collection of pairwise disjoint {\em feathers} attached to the
component $C_i$, $i\ge 2$. A feather is a linear chain of smooth
rational curves on $\bV$. In our particular case, where $V$ is
assumed to be smooth and affine, each of these feathers
consists
of just one smooth rational curve $F_{ij}$ with
self-intersection $w_{ij}:=F_{ij}^2\le -1$. Given $i$ we have $w_{ij}\le
-2$ for at most one of the feathers $F_{ij}$, see Proposition 2.6
in \cite{FKZ3}.

As in \cite{FKZ3} we let $D_\ext^{\ge i}$ denote the branch of
$D_\ext$ at the vertex $C_{i-1}$ containing $C_i$,
while $D_\ext^{> i}$ stands for $D_\ext^{\ge i}-C_i$.

$$
D_\ext^{\ge i}:
\qquad\quad\cou{}{C_i}\nlin\xbshiftup{}{\,\fF_i} \llin\ldots\llin
\cou{}{C_n} \nlin\xbshiftup{}{\,\fF_n} \quad,\quad\qquad
D_\ext^{> i}: \quad\quad\xbshiftup{}{\,\fF_i}\quad
\cou{}{C_{i+1}}\llin\ldots\llin \cou{}{C_n}
\nlin\xbshiftup{}{\,\fF_n} \quad.\qquad\qquad $$
Likewise, we let
$D^{\ge i}=D\cap D_\ext^{\ge i}$ and $D^{>i}=D\cap D_\ext^{> i}$.

Assume now that $(\bV, D)$ is a standard completion so that $m=0$
in (\ref{extended}). The linear systems $|C_0|,|C_1|$ on $\bV$
define a morphism
$$
\Phi=\Phi_0\times\Phi_1:\bV\to \PP^1\times\PP^1,
$$
 called  the {\it standard morphism},
 which is birational according to \cite{FKZ2}, Lemma 2.19.

Decomposing $\Phi$ into a sequence of blowups we can grow $\bV$
starting with the quadric $\PP^1\times\PP^1$, see \cite{FKZ3}. For
a feather $F_{i\rho}$ we let $C_{\mu_{i\rho}}$ denote its {\it
mother component}. The latter means that $F_{i\rho}$ was born by a
blowup with center on $C_{\mu_{i\rho}}$ under this decomposition
of $\Phi$. Since the zigzag $D$ is connected this defines
$C_{\mu_{i\rho}}$ in a unique way.
By Proposition 2.6 in \cite{FKZ3}, $\mu_{i\rho}=i$ (i.e., $C_i$
is the mother component for $F_{i\rho}$) if and only if
$F_{i\rho}^2=-1$, otherwise $\mu_{i\rho} < i$. In the latter case
we say that the feather $F_{i\rho}$ jumped.

\bdefi\label{coin} (1) We let
$$
s_\mu:=\#\{(i,\rho): 2\le i\le n,\,\, 1\le\rho\le r_i,\,\,
\mbox{and }\, \mu_{i\rho}=\mu\}
$$
be the number of feathers $F_{i\rho}$ whose mother component is $C_\mu$.

(2) We say that $C_\mu$ ($2\le\mu\le n$)
is a component {\em of type} $*$, or a
{\em $*$-component} for short, if
\bnum [(i)] \item $D_\ext^{\ge \mu+1}$
is not contractible, and \item $D_\ext^{\ge \mu+1}-F_{ij}$ is not
contractible for every feather $F_{ij}$ of $D_\ext^{\ge \mu+1}$
with mother component $C_\tau$, where $\tau< \mu$.
\enum
Otherwise $C_\mu$ is called a component of type $+$, or simply a
{\em $+$-component}. For instance, $C_2$ and $C_n$
are always components of type $+$.
We let $\tau_\mu=*$ in the first case and
$\tau_\mu=+$ in the second one.
\edefi

\brem\label{inou} It is easily seen that,
in the process of blowing up starting from the
quadric,
every $*$-component $C_\mu$ with $\mu\ge 3$
appears as a result of an inner blowup of the previous
zigzag, while an outer blowup of a
zigzag creates a $+$-component (cf.\ \ref{reversion} and
\ref{1000.001} in the Introduction).
\erem

Given a component $C_\mu$ we denote by $p_{\mu\rho}$, $1\le
\rho\le s_\mu$, the following collection of points on
$C_\mu\backslash C_{\mu-1}\cong \A^1$. For every feather
$F_{\mu\rho}$ of self-intersection $-1$ we let $p_{\mu\rho}$ be
its intersection point with $C_\mu$. This gives $r_\mu$ or
$r_\mu-1$ points on $C_\mu\backslash C_{\mu-1}$ depending on
whether the feathers $F_{\mu r_\mu}$ attached to $C_\mu$ are
all $(-1)$-curves or not. If there is a feather $F_{ij}$  with
mother component $C_\mu$ and with $i>\mu$ then we also add the
intersection point $c_{\mu+1}$ of $C_\mu$ and $ C_{\mu+1}$
to our collection. Thus
$s_\mu$ is one of the numbers $r_\mu-1$, $r_\mu$ or
$r_\mu+1$,
and the points
$$
p_{\mu\sigma}\in C_\mu,\,\, 1\le \sigma\le s_\mu
$$
are just the locations on $C_\mu$ in which the feathers with
mother component $C_\mu$ are born by a blowup.
We call them the {\em base
points} of the associated feathers.

The collection $(p_{\mu\sigma})_{1\le \sigma\le s_\mu}$
defines a point $Q_{\mu}$ in the configuration space $\fM_{s_\mu}^+$.

Suppose further that $C_\mu$ is a component of type $*$. We
consider then $Q_{\mu}$ as a collection of points in
$C_\mu\backslash(C_{\mu-1}\cup C_{\mu+1})$. Note that the
intersection point $c_{\mu+1}$ of $C_\mu$ and $C_{\mu+1}$
cannot be one of the points $p_{\mu\sigma}$ because of (ii) in
Definition \ref{coin}(2). Identifying
$C_\mu\backslash(C_{\mu-1}\cup C_{\mu+1})$ with $\C^*$ in such a
way that $c_{\mu+1}$ corresponds to $0$ and
$c_\mu$ to $\infty$, we obtain a point in
the configuration space $\fM_{s_\mu}^*$. Thus in total we obtain a
point
$$
Q(\bV, D):=(Q_2, \ldots Q_n)\in \fM=\fM^{\tau_2}_{s_2}\times \ldots
\times\fM^{\tau_n}_{s_n}\,
$$
called the {\em configuration invariant} of $(\bV, D)$.

Performing in $(V,D)$ elementary transformations with centers at
the component $C_0$ in (\ref{zigzag}) does neither change $\Phi_0$
nor the extended divisor (except for the self-intersection
$C_1^2$) and thus leaves $s_i$ and $Q(\bV,D)$ invariant.
Hence we can define these invariants for any semi-standard
completion $(\bV, D)$ of $V$ by sending it via elementary
transformations with centers on $C_0$ into a standard completion.
\esit

\subsection{Matching feathers}\label{ss2.3}
In the following proposition we show that reversion of the
boundary zigzag of length $n$ leads to the same configuration
invariant. To formulate this result, it is convenient to use
systematically the notation \be \label{invol} t^\vee= n-t+2 \ee
for an integer $t\in \Z$.

\bprop\label{matchthm} {\rm (Matching Principle)} Let
$V=\bV\setminus D$ be a smooth Gizatullin surface completed by a
standard zigzag $D$. Consider the reversed completion $(\bV^\vee,
D^\vee)$ with boundary zigzag $D^\vee=C_0^\vee\cup \ldots\cup
C_n^\vee$, associated numbers $s_2', \ldots, s_n'$ and
types $\tau_2', \ldots, \tau'_n$. Then
$s_i=s'_{i^\vee}$ and $\tau_i=\tau'_{i^\vee}$ for
$i=2,\ldots, n$. Furthermore, the associated points $Q(\bV,
D)$ and $Q(\bV^\vee, D^\vee)$ in $\fM$ coincide under the natural
identification
$$
\fM=\fM_{s_2}^{\tau_2}\times \ldots \times \fM^{\tau_n}_{s_n}
\cong \fM^{\tau_2'}_{s'_2}\times
\ldots \times \fM^{\tau_n'}_{s_n'}\,.
$$
\eprop

The proof is given in \ref{1.9}-\ref{1.91} below. It uses the
following construction.

\bsit\label{1.9} {\bf Correspondence fibration.} Let us consider a
standard completion $(\bV,D)$ of $V$ and the reversed completion
$(\bV^\vee,D^\vee)$. Thus $D=C_0\cup \ldots\cup C_n$ is a standard
zigzag $[[0,0,w_1,\ldots,w_n]]$ as in (\ref{zigzag}) and
$D^\vee=C^\vee_0\cup \ldots\cup C^\vee_n$ is the standard zigzag
$[[0,0,w_n,\ldots,w_2]]$. Let $D_\ext$ and $D^\vee_\ext$ denote
the corresponding extended divisors and let $F_{i\rho}$ and
$F^\vee_{j\rho}$ be the feathers attached to $C_i$, $C^\vee_j$,
respectively.

Using inner elementary transformations we can move the pair of
zeros in the zigzag $[[0,0,w_2,\ldots, w_n]]$ several places to
the right. In this way we obtain a new completion, say, $(W, E)$
of $V$ with boundary zigzag $E=[[w_2,\ldots, w_{t-1}, 0, 0, w_t,
\ldots,w_n]]$ for some $t\in\{2,\ldots,n+1\}$. For $t=2$,
$E=D=[[0, 0,w_2,\ldots, w_n]]$ is the original zigzag, while for
$t=n+1$, $E=D^\vee=[[ w_2,\ldots, w_n,0, 0]]$ is the reversed
one. The transformed components of $E$ are
$$
E=C_n^\vee \cup \ldots\cup C_{t^\vee}^\vee \cup C_{t-1} \cup
C_t\cup \ldots \cup C_n\,,
$$
where we identify $C_i\subseteq \bV$ and $C_j^\vee\subseteq
\bV^\vee$ with their proper transforms in $W$ ($t-1\le i\le n$,
$t^\vee\le j\le n$). In particular $E=D^{\ge t-1}\cup D^{\vee \ge
t^\vee}$ with new weights $C_{t-1}^2=C^{\vee 2}_{t^\vee}=0$.
There are natural isomorphisms
\be\label{isom}\ba{l}
W\backslash D^{\vee \ge t^\vee}=W\backslash
(C^\vee_n\cup\ldots\cup C^\vee_{t^\vee}) \cong
\bV\backslash (C_0\cup\ldots \cup C_{t-2})\qquad\mbox{and}\\[0pt]
W\backslash D^{\ge t-1} =W\backslash (C_{t-1}\cup\ldots \cup C_n)
\cong \bV^\vee\backslash (C^\vee_0\cup\ldots \cup
C^\vee_{t^\vee -1})\,. \ea
\ee

In the proof of the Matching Principle \ref{matchthm} we
use the following fibration.

\bdefi\label{corr} The map
$$
\psi: W\to \PP^1
$$
defined by the linear system $|C_{t-1}|$ on $W$  will be called
the  {\em correspondence fibration} for the pair of curves $(C_t,
C_{t^\vee}^\vee)$. \edefi

The components $C_t$ and $C_{t^\vee}^\vee$ represent sections of
$\psi$. Since the feathers of $D_\ext$ and $D_\ext^\vee$ are not
contained in the boundary zigzags they are not contracted in $W$.
We denote their proper transforms in $W$ by the same letters. It
will be clear from the context where they are considered. \esit

We use below the following technical facts.

\blem\label{mem} \bnum[(a)]\item The divisor $D_\ext^{\ge t+1}$ is
contained in some fiber $\psi^{-1}(q)$,  $q\in\PP^1$. Similarly,
$D_\ext^{\vee \ge t^\vee +1}$ is contained in some fiber
$\psi^{-1}(q^\vee)$. The points $q$ and $q^\vee$ are uniquely
determined unless $D_\ext^{\ge t+1}$ and $D_\ext^{\vee\ge t^\vee
+1}$ are empty, respectively. \item A fiber $\psi^{-1}(p)$ can
have at most one component $C$ not belonging to $D_\ext^{> t}\cup
D_\ext^{\vee> t^\vee}$. Such a component $C$ meets both $D^{\ge
t}$ and $D^{\vee\ge t^\vee}$.\enum \elem

\bproof (a) follows immediately for the divisor $D_\ext^{\ge t+1}$
since it is connected and disjoint from the full fiber $C_{t-1}$
of $\psi$. By symmetry the assertion holds also for the divisor
$D_\ext^{\vee\ge t^\vee +1}$.

(b) Let $C$ be a component of the fiber $\psi^{-1}(p)$ belonging
neither to $D_\ext^{> t}$ nor to $D_\ext^{\vee> t^\vee}$. We claim
that it meets both $D^{\ge t}$ and $D^{\vee\ge t^\vee}$. Assume on
the contrary that it does not meet e.g., $D^{\ge t}$. Since the
affine surface $V$ does not contain complete curves and
$V=W\backslash E$, we have $C\cdot E=C\cdot D^{\vee\ge t^\vee}\neq
0$. Thus the proper transform $C'$ of $C$ on $\bV^\vee$ must be a
feather of $D_\ext^{\vee}$. Indeed, $C'\cdot C_0^\vee=0$, see
(\ref{isom}). Hence $C'$ is a component of the only degenerate
fiber $(\Phi_0^\vee)^{-1}(0)$ of $\Phi_0^\vee: \bV^\vee\to\PP^1$
and does not belong to $D^\vee$. Since $C'\cdot D^{\vee\ge
t^\vee}\neq 0$, we must have $C'\subseteq D_\ext^{\vee> t^\vee}$
on $\bV^\vee$. This contradicts our assumption that $C$ does not
belong to $D_\ext^{\vee> t^\vee}$ on $W$, and so the claim
follows.

Finally, there can be at most one such fiber component $C$ since
the fiber $\psi^{-1}(p)$ does not contain cycles and meets only
once each of the sections $C_t$ and $C^\vee_{t^\vee}$ of
$\psi$. \eproof

\bcor\label{claim}
If, in the notation as in Lemma
\ref{mem}(a), $q\ne q^\vee$ then each of the divisors $D_\ext^{\ge
t+1}$ and $D_\ext^{\vee\ge t^\vee+1}$ is either empty or
contractible.
\ecor

\bproof We suppose that $q\ne q^\vee$ and $D_\ext^{\ge t+1}\ne
\emptyset$. By Lemma \ref{mem}(a) the latter divisor is contained
in the fiber $\psi^{-1}(q)$. This fiber contains also a component
$C$ meeting the section $C^\vee_{t^\vee}$.  Clearly such a curve
$C$ is neither a component of the zigzag nor a feather of
$D_\ext^{\ge t+1}$ and so not a component of $D_\ext^{\ge t+1}$.
Since $q\ne q^\vee$ the fiber over $q$ cannot contain any
component of $D^{\vee\ge t^\vee +1}_\ext$. Thus by Lemma
\ref{mem}(b) $\psi^{-1}(q)=C\cup D_\ext^{\ge t+1}$. Since the
multiplicity of $C$ in the fiber is $1$, the remaining part
$D_\ext^{\ge t+1}$ of the fiber can be blown down. Symmetrically,
the same holds for $D_\ext^{\vee\ge t^\vee+1}$. \eproof

The following lemma is crucial in the proof of Proposition
\ref{matchthm}, see \ref{1.91} below.

\blem\label{matchlem} Let $F_{i\rho}$ be a feather of the extended
divisor $D_\ext$ attached to $C_i$ and with mother component
$C_\tau$, where $\tau \le t\le i$. Then the following hold.
\bnum[(a)] \item  $F_{i\rho}$ is contained in a fiber
$\psi^{-1}(q_{i\rho})$ on $W$ for some point $q_{i\rho}\in\PP^1$.

\item The fiber $\psi^{-1}(q_{i\rho})$ contains as well a
feather $F_{j\sigma}^\vee$ of $D^{\vee> t^\vee}_\ext$ meeting
$F_{i\rho}$. This feather $F_{j\sigma}^\vee$ has mother component
$C^\vee_{\tau^\vee}$.

\item The feather $F_{j\sigma}^\vee$ in
(b) is uniquely determined by $F_{i\rho}$, and the points
$q_{i\rho}$ in (a) are all different.
\enum \elem

\bproof Let $F_{i\rho}$ be a feather with mother component
$C_\tau$, where $\tau \le t\le i$. Since $F_{i\rho}$ does not meet
$C_{t-1}$ it is vertical with respect to $\psi$ and so contained
in a fiber over some point $q_{i\rho}\in\PP^1$, proving (a).  We
note that by the same reasoning, any feather $F_{j\sigma}^\vee$ of
$D_\ext^{\vee}$ with $j\ge t^\vee$ is a fiber component of $\psi$
on $W$.

To deduce (b), let us start with the case $i=t=\tau$ so that
$F_{t\rho}$ is a $(-1)$-feather. In this case the fiber
$\psi^{-1}(q_{i\rho})$ cannot be irreducible and so
$F_{t\rho}$ meets some other component, say, $C$ of
$\psi^{-1}(q_{i\rho})$. Clearly, $q_{t\rho}\neq q$ (see Lemma
\ref{mem}(a)) and so $C\cdot D_\ext^{\ge t}=0$. By Lemma \ref{mem}
(b) $C$ belongs either to $D_\ext^{\vee>t^\vee}$ or to
$D_\ext^{>t}$. Since $F_{t\rho}$ cannot meet any other feather of
$D_\ext$ and cannot meet the boundary zigzag twice, $C$ must be
one of the feathers, say, $F_{j\sigma}^\vee$ of
$D_\ext^{\vee>t^\vee}$.

Let us show that $F_{j\sigma}^\vee$ has mother component
$C^\vee_{t^\vee}$. Since $F_{t\rho}\cdot C_t=1$, the feather
$F_{t\rho}$ has multiplicity 1 in the fiber
$\psi^{-1}(q_{t\rho})$. Thus the remaining part
$\psi^{-1}(q_{t\rho})- F_{t\rho}$ can be blown down to a
$(-1)$-curve. After this contraction we must still have
$F_{t\rho}^2=-1$. Hence this remaining $(-1)$-curve must be the
image of $C=F_{j\sigma}^\vee$. Moreover, after this contraction
$F_{j\sigma}^\vee\cdot C^\vee_{t^\vee}=1$, hence
$F_{j\sigma}^\vee$ appears under a blowup with center on
$C^\vee_{t^\vee}$.  Thus the mother component of
$F_{j\sigma}^\vee$ is indeed $C^\vee_{t^\vee}$, as stated.

Consider further the case where $i>t=\tau$ so that
$F_{i\rho}$ is contained in the fiber $\psi^{-1}(q)$, see Lemma
\ref{mem}(a). According to Proposition 2.6(b) in \cite{FKZ3} the
divisor $A:=D_\ext^{\ge t+1}- F_{i\rho}$ is contractible to
a point on $C_t$. Since $C_t$ is the mother component of
$F_{i\rho}$ in $D_\ext$, after this contraction $F_{i\rho}$
becomes a $(-1)$-curve with $F_{i\rho}\cdot E=F_{i\rho}\cdot
C_t=1$. Replacing $W$ by the contracted surface $W/A$ and arguing
as before the result follows as well in this case.

If $i,t,\tau$ with $\tau\le t\le i$ are arbitrary, then we pass to
the correspondence fibration $\psi':W'\to\PP^1$ for the pair
$(C_\tau, C^\vee_{\tau^\vee})$, see Definition \ref{corr}. By what
was shown already there is a feather $F_{j\sigma}^\vee$ of
$D_\ext^{\vee
>\tau^\vee}$ with mother component $C^\vee_{\tau^\vee}$ meeting
$F_{i\rho}$ on $W'$. Since $D_\ext^{\vee
>\tau^\vee}\subseteq D_\ext^{\vee
>t^\vee}$ these feathers $F_{i\rho}$ and
$F_{j\sigma}^\vee$ also meet on the surface $W$. Being both fiber
components of $\psi$ (see the proof of (a) above), they meet
within the same fiber. This completes the proof of (b).

Finally,
(c) is a simple consequence of the fact
that the fibers of $\psi$ cannot contain cycles and intersect with
index $1$ each of the sections $C_t,\,C_{t^\vee}^\vee$ of $\psi$.
\eproof

Lemma \ref{matchlem} motivates the following definition.

\bdefi\label{mafe} Consider a pair $(F_{i\rho},\,
F^\vee_{j\sigma})$, where $F_{i\rho}$ is a feather of $D_\ext$
attached to component $C_i$  and $F^\vee_{j\sigma}$ is a feather
of $D_\ext^{\vee}$ attached to component $C_j^\vee$. This pair is
called a {\em pair of  matching feathers}, or simply a {\em
matching pair}, if $i+j\ge n+2$ and $F_{i\rho}$ and
$F^\vee_{j\sigma}$ meet on $V$.\edefi

\brem\label{2800}
Thus if two feathers $F_{i\rho}$ and
$F^\vee_{j\sigma}$ with $i+j\ge n+2$ meet on $V$, then
$(F_{i\rho},\,F^\vee_{j\sigma})$ is a matching pair. By
Lemma \ref{matchlem} every feather $F_{i\rho}$ of $D_\ext$ has a
unique matching feather $F^\vee_{j\sigma}$ of $D_\ext^{\vee}$, and
vice versa. Moreover, if $F_{i\rho}$ has mother component $C_\tau$
then its matching feather $F^\vee_{j\sigma}$ has mother component
$C^\vee_{\tau^\vee}$.

The condition $i+j\ge 2$ here is essential. Indeed, every
feather $F_{t-1,\rho}$ represents a section of $\psi$, hence it
meets every fiber of $\psi$. Since it cannot meet $D_\ext^{\ge
t}$, it meets every feather $F^\vee_{t^\vee,\sigma}$ with $F^{\vee
2}_{t^\vee,\sigma}=-1$ on the affine surface $V$.
\erem



\brem\label{2801}
One can treat in a similar way feathers of
arbitrary length instead
of length one.  Such feathers appear under the minimal
resolution of a singular Gizatullin surface. These are linear
chains
$$
F_{i\rho}:\qquad \co{B_{i\rho}}\lin \co{F_{i\rho 1}}\lin\ldots\lin
\co{F_{i\rho k_{i\rho}}}\quad ,
$$
where the chain $F_{i\rho}-B_{i\rho}$ (if non-empty)
contracts to a singular
point on $V$ and $B_{i\rho}$ is attached to component $C_i$. The
curve $B_{i\rho}$ is also called the {\em bridge curve} of the
feather $F_{i\rho}$.  For instance, an $A_k$-singularity on $V$
leads to an $A_k$-feather, where $k_{i\rho}=k$ and $F_{i\rho}-
B_{i\rho}$ is a chain of $(-2)$-curves of length $k$. Again,
a matching principle provides a one-to-one correspondence between
feathers $F_{i\rho}$ and $F^\vee_{j\sigma}$ such that the mother
component of the bridge curve $B_{i\rho}$ of $F_{i\rho}$ is equal
to the mother component of the tip of $F^\vee_{j\sigma}$. Moreover
$F^\vee_{j\sigma}$ has dual graph
$$
\co{B^\vee_{j\sigma}}\llin \co{F_{i\rho k_{i\rho}}}\llin
\co{F_{i\rho k_{i\rho}-1}}\llin \ldots \llin \co{F_{i\rho
1}}\quad,
$$
so that the tip of $F_{j\sigma}^\vee$ is just $F_{i\rho 1}$.\erem

In the next lemma we show that the reversion respects the type
of components of the zigzag.

\blem \label{2.144} $C_t$ is a $*$-component   if and only if
$C^\vee_{t^\vee}$ is. Furthermore in this case the points $q$
and $q^\vee$ in Lemma \ref{mem}(a) are equal. \elem

\bproof Let $C_t$ be a $*$-component of the zigzag $D$; see
Definition \ref{coin}.(2). Let us first deduce the second
assertion. If on the contrary $q\neq q^\vee$ then by Corollary
\ref{claim} above, $D_\ext^{\ge t+1}$ is contractible
contradicting (i) in Definition \ref{coin}. Thus  $q=q^\vee$.

It remains to check that also $C^\vee_{t^\vee}$ is a
$*$-component i.e., conditions (i) and (ii) in Definition
\ref{coin} are satisfied.

To check (i), we assume on the contrary that $D_\ext^{\vee\ge
t^\vee+1}$ is contractible. After contracting this divisor
in the fiber of $\psi$ over $q=q^\vee$ there is a component $F$
that meets the section $C^\vee_{t^\vee}$ with multiplicity $1$.
The rest, say, $R$ of the remaining fiber is as well
contractible.  Clearly $F$ cannot be a component of the
zigzag $D^{\ge t+1}$. If $F$ is a feather of $D_\ext^{\ge t+1}$
then $R=D_\ext^{\ge t+1}-F$ is contractible, which is only
possible if $F$ has mother component $C_\tau$, $\tau \le t$. The
latter contradicts condition (ii) in Definition \ref{coin}.
Otherwise $F=C$ is an extra component of the fiber $\psi^{-1}(q)$
(see Lemma \ref{mem}(b)), and the argument as in the proof of
Corollary \ref{claim} shows that $R=D_\ext^{\ge t+1}$. Since $R$
is contractible, again we get a contradiction, this time to (i) of
Definition \ref{coin}.

Let us finally  check that (ii) in Definition
\ref{coin} holds. We need to show that for every feather $F^\vee$ of
$D_\ext^{\vee\ge t^\vee+1}$ with mother component $C^\vee_\tau$,
where $\tau \le t^\vee$, the divisor $D_\ext^{\vee\ge
t^\vee+1}-F^\vee$ cannot be contracted. Indeed, otherwise
after contracting this divisor, $F^\vee $ meets the section
$C^\vee_{t^\vee}$. The remaining fiber is $D_\ext^{\ge
t+1}+F^\vee$, since $F^\vee$ meets a matching feather $F$ in
$D_\ext^{\ge t+1}$. Hence $D_\ext^{\ge t+1}$ is contractible, and
again we   arrive at a contradiction.
\eproof

Now we are ready to deduce Proposition \ref{matchthm}.

\bsit\label{1.91}
\bproof[Proof of Proposition \ref{matchthm}] By
Lemma \ref{matchlem}(b) the map $\psi$ provides a one-to-one
correspondence between the feathers of $D_\ext$ with mother
component $C_t$ and the feathers of $D_\ext^{\vee}$ with mother
component $C_{t^\vee}^\vee$. Moreover using Lemma
\ref{matchlem}(c) it provides a  one-to-one correspondence between the
base points $p_{t\rho}$ and $p_{t^\vee,\rho}^\vee$ of the feathers
as considered in \ref{1.2}. By virtue of Lemma \ref{2.144} $C_t$
is a $*$-component  if and only if  $C^\vee_{t^\vee}$ is, proving
the assertion. \eproof \esit

\subsection{Invariance of the configuration invariant}\label{ss2.4}

Theorem \ref{matchmain} and its Corollary \ref{190} below are the
main results of Section 3.

\bthm\label{matchmain} Given two semi-standard completions $(\bV,D)$, $(\bV', D')$   of a Gizatullin surface $V$, for the corresponding
configuration invariants $s_i$, $s_i'$ and $Q(\bV, D)\in\fM$,
$Q(\bV', D')\in\fM'$ as introduced in \ref{1.2} the following
hold.
\begin{enumerate}
\item[(1)] If $(\bV,D)$ and $(\bV',D')$ are evenly linked then
$s_i=s_i'$ for $i=2,\ldots, n$ and the points $Q(\bV, D)$ and
$Q(\bV',D')$ of $\fM=\fM'$ coincide.

\item[(2)] If $(\bV,D)$ and $(\bV',D')$ are oddly linked then
$s_i=s_{i^\vee}'$ for $i=2,\ldots, n$ and the points $Q(\bV, D)$
and $Q(\bV',D')$ of $\fM$ and $\fM'$ coincide under the
identification
$$
 \fM=\fM^{\tau_2}_{s_2}\times \ldots \times\fM^{\tau_n}_{s_n}
\cong \fM^{\tau'_n}_{s'_n}\times \ldots \times\fM^{\tau'_2}_{s_2'}=\fM'\,.
$$\end{enumerate}
\ethm

\bproof (1) By Proposition \ref{1.6} $(\bV',D')$ can be obtained
by a sequence of generalized reversions
$$
(\bV,D)=(\bV_0, D_0)\quad\leadsto\quad
(\bV_1,D_1)\quad\leadsto\quad \ldots \quad\leadsto\quad
(\bV_l,D_l)=(\bV',D')\,,
$$
where $(\bV_i,D_i)$, $0\le i\le l$, are semi-standard completions
of $V$. Moreover $l$ is even if $(\bV,D)$ and $(\bV,D)$ are evenly
linked, and odd otherwise. Hence it suffices to show the
theorem for a generalized reversion of semi-standard
completions. Since
 elementary transformations on $(\bV,D)$ with centers
on the component $C_0$ in (\ref{zigzag}) do not change the
extended divisor (except for the weight $C_1^2$) and leave
$s_i$ and $Q(V,D)$ invariant, we can reduce the statement to the
case where $(V,D)$ and, symmetrically, $(V',D')$ are standard. The
assertion now follows from Proposition \ref{matchthm}. \eproof

\bdefi\label{sycoin}
Given a configuration space
$\fM=\fM^{\tau_2}_{s_2}\times \ldots \times\fM^{\tau_n}_{s_n}$
we consider the reversed product
$$\fM^\vee=\fM^{\tau_n}_{s_n}\times \ldots \times\fM^{\tau_2}_{s_2}\,.
$$ By the {\it symmetric configuration
invariant} of a completion $(\bV,D)$ of a Gizatullin surface $V$
we mean the unordered pair
$$
\tQ(\bV, D)=\left\{ Q(\bV, D), Q(\bV^\vee,
D^\vee)\right\},\quad\mbox{where}\,\, Q(\bV,
D)\in\fM\,\,\mbox{and}\,\, Q(\bV^\vee, D^\vee)\in\fM^\vee.
$$
\edefi

Theorem \ref{matchmain} leads immediately to the following result.

\bcor\label{190} The sequence $(s_i)_{2\le i\le n}$ (up to
reversion) and the pair $\tQ(V):=\tQ(\bV,D)$ are invariants of the
isomorphism  type of $V$. \ecor

\subsection{The configuration invariant for $\C^*$-surfaces}

According to \cite{FlZa1} a normal non-toric $\C^*$-surface
admits a hyperbolic DPD-presentation
$V=\Spec \C[u][D_+,D_-]$. If, moreover, $V$
is Gizatullin, then  there are (not necessarily different)
points $p_\pm$ with $\supp\{ D_\pm\}\subseteq \{p_\pm\}$,
see \cite[Section 4]{FlZa2}. By \cite{FKZ2} $V$
admits an equivariant standard completion $(\bV, D)$,
which is unique up to reversion.
Concerning the structure of this completion
we can summarize the main results from \cite{FKZ2}
and Section 3 in \cite{FKZ3} for smooth $\C^*$-surfaces as follows.
We recall that the {\em parabolic} component is the unique component $C_t$
of the zigzag with $t\ge 2$ consisting
of fixed points of the $\C^*$-action.

\bprop\label{eboundary} If $V$ is non-toric and smooth
then it admits a unique equivariant standard completion $(\tV, D)$
with extended divisor

\vskip0.3truecm

\be\label{ezigzag}
D_{\rm ext}: \qquad\quad\cou{0}{C_0}\lin
\cou{0}{C_1}\lin \cou{w_2}{C_2}\lin\ldots \lin \cou{w_{t-1}}{C_{t-1}}
\llin
\cou{w_t\quad\;\;}{C_t}\nlin\xbshiftup{}{\qquad
\{F_{t\rho}\}_{\rho=0}^r}\lin \cou{w_{t+1}}{C_{t+1}}
\lin\ldots\lin\cou{w_n\quad\;\;}{C_n}
\nlin\boxshiftup{}{ F_n} \quad\quad
\ee
and with boundary zigzag
$D$ represented by the bottom line in (\ref{ezigzag}) such that $C_t$
is an attractive parabolic component. Here
$w_t=\deg (\lfloor D_+ \rfloor+\lfloor D_-\rfloor)\le -2$, $F_n$ is a
single feather (possibly empty) and $\{F_{t\rho}\}_{\rho=0}^r$ is a non-empty
collection of feathers with all $F_{t\rho}$, $\rho\ge 1$,
being $(-1)$-curves. Furthermore the following hold:

\begin{enumerate} [(a)]

\item Suppose that $p_+\ne p_-$ or one of the
fractional parts $\{D_\pm\}$ of the divisors $D_\pm$ is zero.
Then $F_n$ is a $(-1)$-curve\footnote{Hence it is non-empty.},
$F_{t0}$ with $F_{t0}^2=1-t$ has mother component
$C_2$ and $w_i=-2$
for $i\ne 0,1,t$.
Up to equivalence the pair $(D_+,D_-)$ is
\be\label{smoothDPD}
D_+=-\frac{1}{t-1}[p_+], \qquad D_-=-\frac{1}{n-t+1}[p_-]-\sum_{i=1}^r[p_i]
\ee
with pairwise different points  $p_+,p_-,p_1, \ldots ,p_r$, where
$r=-2-w_t\ge 0$.
The feathers $F_{t0},\ldots,F_{tr}$ are attached to the points
$p_+,p_1, \ldots ,p_r$ of $C_t\backslash C_{t-1}\cong \A^1$ whereas $p_-$
corresponds to the intersection point $C_t\cap C_{t+1}$ if $t<n$ and to
$C_n\cap F_n$ if $t=n$.

\item $\{D_+(p_+)\}=0$ iff $t=2$ and, similarly, $\{D_-(p_-)\}=0$ iff $t=n$.

\item Assume that $p_+=p_-=:p$. Then the $F_{t\rho}$ are $(-1)$-feathers
$\forall \rho\ge 0$ while $F_n=\emptyset$ if and only if $D_+(p)+D_-(p)=0$.
Moreover the feathers $F_{t0},\ldots,F_{tr}$ are attached to
the points of the reduced divisor $\lfloor -D_+-D_-\rfloor=
\sum_{i=0}^r[p_i]$ considered as points
of $C_t\backslash C_{t-1}\cong \A^1$ while $p_-$
corresponds to the intersection point $C_t\cap C_{t+1}$.
\end{enumerate}
\eprop

For the proof we refer the reader to \cite[\S 3]{FKZ3},
in particular to Proposition 3.10 and Remark 3.11(2).

\bsit\label{stcond} We recall the following conditions ($\alpha_*$) and
($\beta$) of Theorem 0.2 in \cite{FKZ3}.

\bnum \item[($\alpha_*$)] $\supp \{D_+\}\cup \supp \{D_-\}$ is
either empty or consists of one point  $p$, where $D_+(p) + D_-(p)
\le -1$ or both fractional parts $\{D_+(p)\}$, $\{ D_-(p)\}$ are
nonzero.

\item[($\beta$)] $\supp \{D_+\}=\{p_+\}$ and $\supp
\{D_-\}=\{p_-\}$ for two distinct points $p_+,p_-$, where
$D_+(p_+) + D_-(p_+)\le-1$ and $D_+(p_-) + D_-(p_-)\le-1.$
\enum
\esit

By Theorem 0.2
of \cite{FKZ3} for a normal $\C^*$-surface
satisfying one of these condition the $\C^*$-action is unique up
to conjugation and inversion. In the next proposition we clarify
for which smooth surfaces these conditions do not hold. This
yields in particular part (a) of Theorem \ref{0.5} from the
Introduction.

\bprop\label{crit} Let $V$ be a smooth Gizatullin $\C^*$-surface
with DPD-presentation $V=\Spec \C[u][D_+,D_-]$.
Then the following conditions are equivalent.
\bnum

\item[(i)] Neither ($\alpha_*$) nor ($\beta$) is fulfilled.

\item[(ii)] $(D_+,D_-)$ is equivalent to a pair as in
(\ref{smoothDPD}) with $n\ge 3$.


\item[(iii)] $V$ is special with $n\ge 3$.
\enum
\eprop

\bproof
If (i) holds then $\supp \{D_+\}\cup \supp\{D_-\}$ is non-empty and so
by Proposition \ref{eboundary}(b) $n\ge 3$ in (\ref{ezigzag}).
Moreover $p_+\ne p_-$ or one of the fractional parts
$\{D_\pm\}$ is zero.
Hence (ii) follows from \ref{eboundary}(a).
The implication (ii)$\To$(i) is easy and left to the reader.

If (ii) holds then inspecting (\ref{ezigzag})
$V$ is special with $n\ge 3$ and so (iii) holds.
To prove the converse, assume that $V$ is special  with $n\ge 3$.
If $p_+\ne p_-$ or one of
the divisors $\{D_\pm\}$ is zero, then we can conclude by Proposition
\ref{eboundary}(a). So assume that both divisors $\{D_\pm\}$
are non-zero and $p_+=p_-$. In particular by \ref{eboundary}(b)
$3\le t\le n-1$.

According to Definition \ref{0.2},
the extended divisor has a feather $F_2$ with
mother component $C_2$ and another one $F_n$ with mother component
$C_n$. Comparing with (\ref{ezigzag})  $F_2$ must be attached to
$C_t$ with $F_2^2\le-2$. This contradicts Proposition
\ref{eboundary}(c). \eproof

\brem\label{nrm} If $V$ is as in
(\ref{smoothDPD}) with $r=0$, then $V$ is a Danilov-Gizatullin
surface, see \cite[5.2]{FKZ2} or \ref{dgs} below. Furthermore, if
$n\ge 3$ then $V$ is special of type I if either $r=1$ or $r\ge 2$
and one of the fractional parts $\{D_\pm\}$ vanishes (i.e., $t=2$ or $t=n$).
Otherwise it is of type II. \erem

We note the following important consequence.

\bcor\la{specialconfig}
Assume that $V=\Spec \C[u][D_+,D_-]$ is a special smooth $\C^*$-surface with
$$
\lfloor -D_+-D_-\rfloor=q_1+\ldots+ q_s\,.
$$
Then the following hold.
\bnum
\item The configuration invariant of $V$ is given by the point
$(q_1,\ldots, q_s)\in \fM_s^+$.
\item The numbers
$\deg\,\{D_+\}$, $\deg\,\{D_-\}$ and $s$ uniquely determine
the zigzag $D$ of a standard completion of $V$ up to reversion.
\enum
\ecor

\bproof
According to Proposition \ref{crit} $(D_+,D_-)$ is up to
equivalence a pair
as in (\ref{smoothDPD}).
Using the description of the zigzag in Proposition \ref{eboundary},
(2) follows.

To deduce (1), we assume first that $2<t<n$.
With the notations as in Proposition \ref{eboundary}(a),
the feather $F_{t0}$ has mother component $C_2$ while
the feathers $F_{t1}, \ldots, F_{tr}$ are $(-1)$-feathers
attached to $C_t\backslash C_{t-1}\cong \A^1$
in the points $p_1,\ldots, p_r$. Since in this case
$\lfloor D_++D_-\rfloor=p_1+\ldots+ p_r$, the assertion follows.

In the case $t=2<n$ we have
$\lfloor D_++D_-\rfloor=p_++p_1+\ldots+ p_r$ while $F_0$
is an additional $(-1)$-feather attached to $C_2=C_t$ in $p_+$.
Hence we can conclude as before.
Replacing $F_0$ by $F_n$ and $p_+$ by $p_-$
the same argument works also in the case $2<t=n$.
\eproof

\brem\la{nonspecial}
Let $V=\Spec \C[u][D_+,D_-]$ be a non-special
smooth $\C^*$-surface
such that the divisors $\{D_+\}$ and $\{D_-\}$ are
both nonzero and supported on the same point $p=p_+=p_-$.
In this case the parabolic component $C_t$
in (\ref{ezigzag})
is of $*$-type. Indeed, condition (ii) in Definition
\ref{coin} is empty. Condition (i) follows
in the case $D_+(p)+D_-(p)\ne 0$
by Lemma 3.21 in \cite{FKZ3} and in the case
$D_+(p)+D_-(p)=0$ from the fact
that $F_n=\emptyset$, see \ref{eboundary}(c).

Now again the conclusion of Corollary \ref{specialconfig} holds.
Indeed, part (1) with $\fM^+_s$ replaced by $\fM^*_s$ is a
consequence of Proposition \ref{eboundary}(c), while part (2)
follows from the fact that the weights $w_2,\ldots, w_{t-1}$ and
$w_{t+1},\ldots, w_{n}$ of the boundary zigzag in (\ref{ezigzag})
are uniquely determined by $\deg \{D_+\}$  and  $\deg \{D_-\}$,
respectively (see \cite[Proposition 3.10]{FKZ3}). \erem

The preceding remark can be used to deduce uniqueness of
$\C^*$-actions for all non-special
smooth Gizatullin surfaces;
cf.\ the more general Theorem 0.2 in \cite{FKZ3},
which covers as well
the singular case.

\bcor\la{unique*}
If $V$ is a non-special smooth $\C^*$-surface
then its $\C^*$-action is unique up to equivalence.
\ecor

\bproof
Assume that $V=\Spec \C[t][D_+,D_-]$ is non-special.
Then either both $D_\pm$ are integral, or both fractional parts $\{D_\pm\}$
are non-zero and supported by the same point
$p=p_+=p_-$.

Suppose first that both $D_\pm $ are integral.
Comparing with Proposition \ref{eboundary}(a) $D$
is then a zigzag with $n=t=2$, and the configuration invariant
of $V$ is given by $(p_-,p_+,p_1,\ldots, p_r)$.
Hence the pair $(D_+,D_-)$ is uniquely determined up to equivalence.

Suppose now that both fractional parts $\{D_\pm\}$
are non-zero and supported by the same point $p$. Let
$\lfloor -D_+-D_-\rfloor=\sum_{i=0}^rp_i.$
By Remark \ref{nonspecial} and Theorem \ref{eboundary}(c)
the values $\{D_+(p)\}$ and $\{D_-(p)\}$ are, up to interchanging,
uniquely determined by the boundary zigzag and so
by the abstract isomorphism type of $V$.
Since $V$ is smooth, using Theorem 4.15 in \cite{FlZa1}
we have $-1<D_+(p)+D_-(p)\le 0$. Hence if we require that
$-1< D_+(p)<0$ then $D_\pm (p)$ are uniquely determined.

Applying again Remark \ref{nonspecial}, the parabolic component
$C_t$ is of $*$-type. By  Theorem \ref{eboundary}(c) under
a suitable isomorphism $C_t\backslash C_{t-1}\cong \A^1$
the configuration invariant of $V$ is the point in $\fM^*_{r+1}$
given by the subset
$\{p_0,\ldots, p_r\}$ of $C_{t}\backslash ( C_{t-1}\cup C_{t+1})$.
Since
by {\em loc.cit.\ } $p$ corresponds to the intersection point
$C_t\cap C_{t+1}$ the abstract isomorphism type of $V$
determines the pair $(D_+,D_-)$ up to equivalence.
\eproof

\section{Special Gizatullin surfaces of ($-1$)-type}

The main result of this section says that the isomorphism
type of a special surface (as introduced in Definition \ref{0.2})
is
uniquely determined by its configuration invariant provided
that all feathers are $(-1)$-curves, see Proposition  \ref{main}. To
this purpose we introduce presentations of Gizatullin surfaces.
Given such a surface $V$, we define certain natural group actions
on the set of all presentations of $V$. They change the completion
while leaving the affine surface unchanged. The most important
ones are a 2-torus action, elementary shifts, and backward
elementary shifts. Given two special smooth Gizatullin surfaces
$V,\,V'$of $(-1)$-type with the same zigzag and configuration invariant, we
show that they admit a common presentation, hence are isomorphic
(see the proof of Proposition  \ref{main}). To achieve this we
gradually change two given presentations of $V$ and $V'$ by means
of the above actions until they become equal.

\subsection{Presentations}\label{ss3.5}
Every smooth Gizatullin surface $V$ can be constructed along with
a standard completion $(\bV,D)$ via a sequence of blowups starting
from the quadric $\PP^1\times\PP^1$. If all components $C_i$ of $D$
with $i\ge 2$ are of
$+$-type then the necessary sequence of blowups can be described
in the following way (cf.\  Corollary \ref{DG+special}).

\bsit\label{basic0} We let $Q=\PP^1\times \PP^1$ denote the
quadric, where $\PP^1=\A^1\cup \{\infty\}$. In $Q$ we consider the
curves
$$
C_0=\{\infty\}\times \PP^1\,,\quad C_1= \PP^1\times
\{\infty\}\,,\quad \mbox{and}\quad C_2=\{0\}\times \PP^1\,.
$$
We choose finite sets of points
$$
M_2, c_3, M_3, \ldots, c_n, M_n
$$
among the points of $C_2$ and infinitesimally near points as
follows. \bnum \item[(a$_2$)] $M_2\subseteq C_2\backslash C_1$ is
a finite subset.  Blowing up $X_2=Q$ with centers at the points of
$M_2$ we obtain a surface $X_3$.

\item[(b$_2$)] $c_3\in X_3$ is a point on $C_2\backslash C_1$,
where by abuse of notation we identify the curves $C_1$, $C_2$
with their proper transforms in $X_3$. Blowing up $X_3$ with
center at $c_3$ leads to a surface $\tilde X_3$ with an
exceptional curve $C_3$ over $c_3$.

\item[(a$_3$)] $M_3\subseteq C_3\backslash C_2$ is a finite
subset, where we identify the curves $C_2$ and $C_3$ with their
proper transforms in $\tilde X_3$. We assume that none of the
points of $M_3$ is contained in an exceptional curve over $M_2$.
Blowing up $\tilde X_3$ with centers at the points of $M_3$ we
obtain a surface $X_4$.

\item[(b$_3$)] $c_4\in X_4$ is a point on $C_3\backslash C_2$,
where we identify again the curves $C_2$ and $C_3$ with their
proper transforms in $X_4$. Blowing up $c_4$ leads to a surface
$\tilde X_4$ with an exceptional curve $C_4$. \enum Iterating this
procedure we finally arrive at a smooth  rational projective
surface \be\label{manylet} X_n=X(M_2, c_3, M_3, \ldots,
c_n,M_n)\,. \ee We emphasize that in each step (a$_i$) we require
that \be\la{*} \mbox{none of the points of $M_i$ is in an
exceptional curve over $M_j$ with $j<i$. } \ee As before we
identify the curves $C_i$ with their proper transforms in $X_n$.
The smooth open surface
$$
V_n= V(M_2, c_3, M_3, \ldots, c_n,M_n)=X_n\backslash D\,,
$$
where $D=C_0\cup\ldots\cup C_n\subseteq X_n$ stands for the
boundary zigzag, is an affine Gizatullin surface, see Lemma
\ref{2.6}.

\esit \bdefi\label{basic} We call $X_n$ as in (\ref{manylet}) a
{\em presentation} of $V=V_n$. It is called a presentation of {\em
$(-1)$-type}, or simply a {\it $(-1)$-presentation}, if
$c_{i+1}\not\in M_i$  $\forall i$. Equivalently, this means that
at each step (b$_i$) the point $c_{i+1}$ is not contained in any
of the exceptional curves over the points of $M_i$.

We say that a standard completion $(\bV,D)$ of a Gizatullin
surface $V$ is of {\em $(-1)$-type}, or simply a {\it
$(-1)$-completion}, if all feathers of the extended divisor
$D_\ext$ are $(-1)$-feathers or, equivalently, are attached to
their mother components. It is easily seen that $(X_n,D)$ as above
gives a $(-1)$-completion of $V_n$ if and only if $X_n$ is a
$(-1)$-presentation.
\edefi

\bsit\label{prsp}
We let $\Pi (s_2,\ldots,s_n)$ denote
the space of all presentations $X(M_2, c_3,\ldots ,
c_n,M_n)$ with $|M_i|=s_i$. The tower of smooth fibrations
$$\Pi (s_2,\ldots,s_n)\to \Pi (s_2,\ldots,s_{n-1})\to
\ldots\to\{\cdot\}
$$ shows that $\Pi (s_2,\ldots,s_n)$ is a smooth
quasiprojective  variety (non-affine, in general).
\esit

The next lemma is immediate from the construction in
\ref{basic}; we leave the details to the reader.

\blem\label{2.6} For a presentation as in \ref{basic} and with
$s_i=|M_i|$ the following hold.
\bnum[(a)] \item If $n\ge 3$
then $D=C_0\cup \ldots \cup C_n\subseteq X_n$ represents a
zigzag
\be\label{2.60} \G_D=[[0,0,-s_2-1, -s_3-2,\ldots,
-s_{n-1}-2, -s_n-1]]\,,
\ee
while for $n=2$ we have
$\G_D=[[0,0,-s_2]]$. Consequently, if $s_2,s_n\ge 1$ for $n\ge
3$ and $s_2\ge 2$ for $n=2$, then $V_n$ is a Gizatullin surface
with standard completion $(X_n, D)$.

\item Letting $M_i=\{p_{ij}\}$, for a point $p_{ij}\in M_i$ we let
$F_{ij}$ denote\footnote{Attention: now $i$ stands for the
index of the mother component, whereas in Section 3 it means the
index of the component of attachment.} the proper transform in
$X_n$ of the exceptional curve over $p_{ij}$ that was generated in
step $(a_i)$. Then $F_{ij}$ is a feather of $D_\ext$ with mother
component $C_i$, and every feather of $D_\ext$ appears in this
way.

\item
All components $C_2,\ldots, C_n$ of the zigzag $D$
are $+$-components.\footnote{See Definition
\ref{coin}.}

\item The configuration invariant $Q(\bV,
D)\in\fM=\fM^+_{s_2}\times\ldots\times\fM^+_{s_n}$ of $V$ as in
\ref{1.2} is given by the sequence $(M_2,\ldots, M_n)$,
where $M_i$ is viewed as a point in the configuration space
$\fM^+_{s_i}$ of $s_i$-tuples of distinct points in
$C_i\backslash C_{i-1}\cong\A^1$.
\enum \elem

To show that $V_n$ in (a) is affine it suffices to observe
that the zigzag $D$ supports an ample divisor $\sum_{i=0}^n m_i
C_i$ with $0<m_0\ll m_1\ll \ldots\ll m_n$, due to the Nakai-Moishezon
criterion  (cf.\ e.g., \cite{Gi,Du}).

We have the following criteria for a surface to admit a
presentation.

\blem\label{2.7} Let $(\bV,D)$ be a standard completion of a
smooth Gizatullin surface $V$, and let $\Phi:\bV\to
Q=\PP^1\times\PP^1$ be the standard morphism so that $\bV$ is
obtained from $Q$ by a sequence of blowups. Then the following
conditions are equivalent.

\bnum[(a)]
\item Every standard completion $(\bV,D)$ of $V$ arises
from a presentation $\bV=X_n$.

\item $V$ arises from a presentation as in \ref{basic}.

\item Every component $C_{s+1}$, $s=2,\ldots,n-1$,
of the zigzag is created by a blowup on $C_s\backslash
C_{s-1}$.

\item The dual graph of $D$ is as in (\ref{2.60}), where
$s_i$ is the number of feathers with mother component $C_i$.
\enum
\elem
\bproof
(a)$\To$(b) is trivial while
(b)$\To$(c) follows from the definitions. To deduce
(c)$\To$(d) and (d)$\To $(a) we
proceed by induction on the length $n$ of the zigzag. In case
$n=2$ both implications are evident. Clearly every feather
with mother component $C_n$ is
a $(-1)$-feather. If $n\ge 3$, blowing down all $(-1)$-feathers
of $C_n$ this component becomes a $(-1)$-curve and can be
blown down too. This results in a zigzag of  shorter length
and so the induction argument works.

Since by Proposition  \ref{main} condition (d) does not depend
on the choice of the completion,
(d)$\To$(a) follows as well.
\eproof

\bcor\la{DG+special} Every standard completion of a special
surface (see Definition \ref{0.2} (a))  admits a presentation.
\ecor

\subsection{Reversed presentation}

In this subsection we study how a presentation changes when
we reverse the boundary zigzag.
The matching principle \ref{matchlem} yields the following result.

\bcor\label{2.8}
Let $X_n$ as in \ref{basic} be a presentation
of a special smooth Gizatullin surface $V=X_n\backslash D$, and
let $(X_n^\vee, D^\vee)$ be the reversion of the completion $(X_n,
D)$ with reversed zigzag $D^\vee=C_0^\vee\cup\ldots\cup C_n^\vee$.
Then there is a presentation
$$
X_n^\vee\cong X(M_n, c_3^\vee,M_{n-1}, \ldots, c_n^\vee, M_2)
$$ for suitable points $c_3^\vee, \ldots, c_n^\vee$,
where we identify $M_\mu$ with a subset of
$C_{\mu^\vee}^\vee\backslash C_{\mu^\vee-1}^\vee$ via the
correspondence fibration $\psi$ for the pair $(C_\mu,
C^\vee_{\mu^\vee})$ as in Definition \ref{corr}. \ecor



Next
we describe the positions of matching
feathers in $X_n$ and $X_n^\vee$.

\bprop\label{4-6}
Let $X=X_n$ be as in Definition
\ref{basic}, and let $F$ and $F^\vee$ be a
pair of matching feathers with mother components $C_\mu$ and
$C^\vee_{\mu^\vee}$, respectively, and with the same base
point $p\in M_\mu$. Then the following hold.

\bnum[(a)]
\item
$F$ is attached to component $C_{k+1}$ with $\mu< k+1< n$ if
and only if
$$
(i)\quad p=c_{\mu+1},\quad c_{i+1}=c^\vee_{i^\vee+1}
\quad\mbox{for}\quad \mu< i\le k, \quad\mbox{and}\quad (ii)\quad
c_{k+2}\ne c^\vee_{k^\vee}\,.
$$
Similarly, $F$ is attached to component $C_{n}$ with $\mu< n$
if and only (i) is satisfied with $k=n-1$.

\item $F^\vee$ is attached to
component $C^\vee_{l^\vee+1}$ with $\mu^\vee< l^\vee+1< n$ if and
only if
$$
(i)\quad p=c^\vee_{\mu^\vee+1},  \quad c_{i+1}=c^\vee_{i^\vee+1}
\quad\mbox{for}\quad l\le i<\mu, \quad \mbox{and}\quad
(ii)\quad c_{l}\ne c^\vee_{l^\vee+2}\,.
$$
Similarly, $F^\vee$ is attached to component $C^\vee_{n}$ with
$\mu^\vee< n$
if and only (i) is satisfied with $l^\vee=n-1$.
\enum \eprop

\bproof To deduce (a), assume that $F$ is attached to $C_{k+1}$
with $\mu\le k$. It is clear that then $p=c_{\mu+1}$ since
otherwise in the construction the feather $F$ would be attached to
$C_\mu$ as a $(-1)$-feather. Let us consider the correspondence
fibration $\psi:(W,E)\to \PP^1$ for the pair $(C_i, C_{i^\vee})$
as in Definition \ref{corr}. In the notation of Lemma \ref{mem}(a)
we have $q=c_{i+1}$ and $q^\vee=c^\vee_{i^\vee+1}$. If $F$ is
attached to $C_{k+1}$ and $\mu+1\le i\le k$ then $F$ appears in the
fiber $\psi^{-1}(c_{i+1})$ (cf.\ Lemma \ref{mem}(a)). Since $F$
meets its matching feather $F^\vee$ and the latter sits  in the
fiber $\psi^{-1}(c_{i^\vee +1}^\vee)$, this forces $c_{i+1}=
c_{i^\vee +1}^\vee$ in this range. However, if $i=k+1$ and $k+1\le n-1$
then $F$ is not any longer contained in $\psi^{-1}(c_{k+2})$ while
$F^\vee$ and then also $F$ is still in
$\psi^{-1}(c_{k^\vee}^\vee)$. This shows that indeed $c_{k+2}\ne
c_{k^\vee}^\vee$.

Dualizing (a) also (b) follows.
\eproof

\bcor\label{4-7} Suppose that $F$ is a feather of $X$ with
mother component $C_\mu$ and $G^\vee$ a feather of $X^\vee$ with
mother component $C^\vee_{\nu^\vee}$. If $F$ and $G^\vee$ are attached
to the components $C_{k+1}$ of $X$ and $C^\vee_{l^\vee+1}$ of
$X^\vee$, respectively, then  the intervals of integers
$[\mu,k+1]$ and $[l-1,
\nu]$ have at most one point in common.  \ecor

\bproof
If $\mu=k+1$ or $\nu^\vee=l^\vee +1$, i.e.\ $\nu=l-1$,
then the assertion is trivial.
So assume for the rest of the proof that $\mu<k+1$ and $\nu>l-1$.
Let $p\in M^\vee_{\nu^\vee}\cong M_\nu$ be the base point of $G^\vee$.
We claim that:

(1) $\nu\not\in ]\mu,k]$ and, dually, $\mu\not\in [l,\nu[$.
Clearly it suffices to show the first part. If on the contrary
$\nu\in ]\mu,k]$ then by Proposition \ref{4-6}(a)
$c_{\nu+1}=c^\vee_{\nu^\vee+1}$ while by \ref{4-6}(b)
$p=c^\vee_{\nu^\vee+1}$ hence $c_{\nu+1}\in M_\nu$. This
contradicts (\ref{*}) in \ref{basic0}.

(2) $l-1\not\in ]\mu,k]$ and, dually, $k+1\not\in [l,\nu[$.
Again it suffices to show the first part.
If on the contrary $l-1\in ]\mu,k]$
then by (i) in \ref{4-6}(a) $c_{l}=c^\vee_{l^\vee+2}$,
contradicting (ii) in \ref{4-6}(b).

(3) $k+1\ne\nu$ and, dually, $l-1\ne\mu$.
As before it suffices to show the first part.
If on the contrary $\nu=k+1$
then by Proposition \ref{4-6}(b) $p=c^\vee_{\nu^\vee+1}$.
The feather $F$ is attached to $C_\nu$ and contained in the
fiber over $c^\vee_{\nu^\vee+1}$,
since it has to meet its dual feather.
Moreover $G^\vee$ and then also
its dual feather $G$ are in this fiber.
Thus the two feathers $G$ and $F$
are attached to the same point of $C_\nu$,
which gives a contradiction.

Obviously (1)-(3) imply our assertion.
\eproof

\subsection{Actions of elementary shifts on
presentations}\label{ss3.6} Here we develop our principal tool in
the proof of the main theorem.

\bsit\label{2.9}
We fix a coordinate system $(x,y)$ on the
affine plane $\A^2$, and we let $\Aut_y(\A^2)$ denote the group
of all automorphisms $h:\A^2\to \A^2$ stabilizing the
$y$-axis $\{x=0\}$. Such an automorphism can be written as
\be\label{auto} h:(x,y)\lto(ax,by+P(x)), \quad\mbox{where }a,b\in
\C^* \mbox{ and } P\in \C[x]. \ee
Clearly $h$ extends to a
birational transformation of the quadric $Q=\PP^1\times\PP^1$,
which is biregular on $Q\setminus C_0$ and stabilizes the curves
$C_1$, $C_2$
 \footnote{With notation as
in Definition \ref{basic}, we identify the $y$-axis with
$C_2\setminus C_1$.}.
\esit

The group $\Aut_y(\A^2)$ acts on presentations
$X_n:=X(M_2, c_3,\ldots, c_n,M_n)$. Indeed,
given $h\in\Aut_y(\A^2)$,
in the inductive
construction of \ref{basic0} the set $M_2\subseteq C_2$
is moved by $h$ into a new set of point,
say $M_2'\subseteq C_2$. Thus $h$ induces a morphism
$X(M_2)\backslash C_0\to X(M_2')\backslash C_0$.
Under this map $c_3$ is mapped to a point $c_3'$,
yielding again a morphism
$X(M_2, c_3)\backslash C_0\to X(M_2', c_3')\backslash C_0$.
Continuing in this way we obtain finally a transformed presentation
$$
h_*(X_n):=X_n':=X(M_2', c_3', \ldots,c_n', M_n')
$$
together with an isomorphism (also denoted by $h$)
\be\la{map}
h: X(M_2, c_3,\ldots, c_n,M_n)\backslash C_0
\stackrel{\cong}{\longrightarrow} X(M_2', c_3',
\ldots,c_n', M_n')\backslash C_0\,.
\ee
In particular the affine surfaces
$V=X_n\backslash D$ and $V':=X_n'\backslash D'$
are isomorphic under $h$. Furthermore,
$h$ maps $D_\ext- C_0$ isomorphically onto $D'_\ext- C_0$, where
$D_\ext$ and $D_\ext'$ are the extended divisors of $X_n$ and $X_n'$,
respectively.

\brems\label{2.13}
1. The automorphism $h$ as in (\ref{auto})
extends to an automorphism of the Hirzebruch surface
$\Sigma_{t}$, where $t=\deg (P)$. We can replace $X_n$ and $X_n'$
by the corresponding semi-standard completions of $V$ and $V'$,
respectively, with boundary zigzags $[[0,-t,\ldots]]$, by
performing on both surfaces a sequence of inner elementary
transformations with centers at $C_0\cap C_1$. Then $h$ extends to
a biregular map between these new completions, sending $D_\ext$
isomorphically onto $D'_\ext$.

2. It is easy to see that the assignment
$$
X(M_2, c_3,\ldots, c_n,M_n)\mapsto X(M_2', c_3', \ldots,c_n', M_n')
$$
defines a regular action of the group
$\Aut_y (\A^2)$ on the presentation space $\Pi (s_2,\ldots,s_n)$
as in \ref{prsp}.
\erems

Let us study in detail the action of the {\em
elementary shifts} $h=h_{a,t}\in\Aut_y(\A^2)$, where
\be\label{elmap} h_{a,t}: (x,y)\mapsto (x,y+ax^{t-2})\qquad\mbox{
with}\quad a\in \C\quad\mbox{and}\quad t\ge 2\,.  \ee

\blem\label{2.12} Let $X_n$ as in (\ref{basic}) be a
$(-1)$-presentation. Then for every $a\in \C$ and $t\ge 2$ the
elementary shift $h=h_{a,t}$ induces the identity on
$C_2\cup \ldots \cup C_{t-1}$ and a translation $x\mapsto x+a$ in
a suitable coordinate on $C_t\setminus C_{t-1}$. In particular,
$$ X_n'=
h_*(X_n)= X(M_2, c_3,\ldots, M_{t-1},c_{t}, a+M_t,
a+c_{t+1},M_{t+1}',\ldots, c_n',M_n')
$$
for some $c_{i+1}'$ and $M_i'$ for $i>t$.
\elem

\bproof For $t=2$ the assertion is evidently true. So we assume in
the sequel that $t\ge 3$. Since $X_n$ is a $(-1)$-presentation it
can be obtained by first creating the zigzag $D=C_0\cup\ldots\cup
C_n$ by successive blowups with centers at $c_3,\ldots,c_n$ and
then blowing up the points of $M_2,\ldots,M_n$. Let as before
$C_2\setminus C_1=\{x=0\}$ in coordinates $(x,y)$ in $\A^2$. After
a suitable translation we may suppose that $c_3=(0,0)$. The blowup
with center at $c_3$ can be written in coordinates as
$$
(x_3, y_3)=(x, y/x)\,,\quad\mbox{or, equivalently,}\quad
(x,y)=(x_3, x_3y_3).
$$
In these coordinates, the exceptional curve $C_3$ is given by
$x_3=0$ and the proper transform of $C_2$ by $y_3=\infty$.
The elementary shift $h_{a,t}$ can be written as \be\label{transl}
h_{a,t}: (x_3,y_3)\mapsto (x_3, y_3+ax_3^{t-3})\,. \ee In
particular $h_{a,3}$ yields the identity on the curve
$C_3\backslash C_2\cong\A^1$ if $t>3$ and the translation by $a$
if $t=3$. The formulas (\ref{transl}) remain the same after
replacing the coordinates $(x_3, y_3)$ by the new ones $(x_3,
y_3-\delta_4)$, where $c_4=(0,\delta_4)$. Thus we may assume that
$c_4=(0,0)$ in the coordinate system $(x_3, y_3)$. Now the lemma
follows easily by induction. \eproof

\bsit\label{2.100} Consider now  the action of the 2-torus
$\T=\C^*\times\C^*$ on $\A^2$ given by
$$
(\lambda_1,\lambda_2).(x,y)=(\lambda_1 x,\lambda_2 y),\quad
(\lambda_1,\lambda_2)=\lambda\in\T\,.
$$
It leaves both axes invariant and extends to a biregular action on
the quadric $Q=\PP^1\times\PP^1$. Hence in this case the map
$h$ in (\ref{map}) is biregular, while in general it can
have points of indeterminacy. By Lemma \ref{2.12} we can use the
induced $\T$-action on presentations in the following way. \esit

\bcor\label{2.14}
Let as before $X_n=X(M_2, c_3,\ldots , c_n,M_n)$
be a $(-1)$-presentation.
Then it can be transformed by a suitable sequence
of elementary shifts into a new one such that the points
$c_3,\ldots, c_n$ and one of the points $p_{n,j}\in C_n$ to which
a feather is attached, are fixed by the induced $\T$-action.
\ecor

\bproof
After a suitable translation we may assume that
$c_3=(0,0)$, so this point is fixed by the torus action. Using
induction on $t$, suppose that the presentation is already
transformed by a sequence of elementary shifts into a new one such
that $c_3, \ldots , c_t$ are  invariant under the torus
action. Then for $t<n$, $\T$ acts on $C_t$. Applying Lemma
\ref{2.12}, the elementary shift $h_{a,t}$ with a suitable $a\in
\C$ moves $c_{t+1}$ into the second fixed point of the $\T$-action
on $C_t$, as required. Similarly, if $t=n$ then we can achieve
that $p_{n,j}$ is the second fixed point of $\T$ on $C_n$.
\eproof

\blem\label{2.149} Let
\[
\ba{l}
X_n=X(M_2, c_3,\ldots, M_{t-1}, c_{t}, M_{t},c_{t+1},\ldots ,c_n,M_n)
\quad\mbox{and}\\
X_n'=X(M_2, c_3,\ldots, M_{t-1}, c_{t}, M'_{t}, c'_{t+1},\ldots
,c'_n,M'_n)
\ea
\]
be presentations with reversed presentations\footnote{See
Corollary \ref{2.8}.}
$$\ba{l}
X_n^\vee=X(M_n, c^\vee_3,\ldots, c^\vee_{t^\vee}, M_{t},
c^{\vee}_{t^\vee+1}, M_{t-1},
\ldots ,c^\vee_n,M_2)\qquad\mbox{and}\\
X_n^{\prime\, \vee}=X(M'_n, c^{\prime\, \vee}_3,\ldots,
c^{\prime\, \vee}_{t^\vee}, M_{t}', c^{\prime\, \vee}_{t^\vee+1},
M_{t-1}, \ldots ,c^{\prime\, \vee}_n,M_2)\,, \ea
$$
respectively. Then we have
\be\label{999} c^{ \vee}_{i}=c^{\prime\, \vee}_{i}\qquad
\forall i=t^\vee+1, \ldots, n\,. \ee
\elem

\bproof Starting with the completion $(X_n,  D)$  of $V$ we
consider the correspondence fibration $\psi:(W,E)\to \PP^1$ for
the pair of curves $(C_t, C^\vee_{t^\vee})$ as in Definition \ref{corr}.
Similarly, we let $\psi':(W',E')\to \PP^1$ be the correspondence
fibration associated to  $(X'_n, D')$ for the pair of curves
$(C'_t, C^{\prime\,\vee}_{t^\vee})$. To obtain the part
$D^{\vee\ge t^\vee}=C^\vee_n\cup \ldots \cup C^\vee_{t^\vee}$ of
the reversed zigzag $D^\vee$ only inner elementary transformations
with centers at the components $C_0=C_0', \ldots,
C_{t-1}=C_{t-1}'$ are required. It follows that
$C^\vee_{i^\vee}=C^{\prime\, \vee}_{i^\vee}$ for all $i\le t-1$
i.e., for all $i^\vee\ge t^\vee+1$. In particular (\ref{999})
hold. \eproof

Our next aim is to study the behavior of $(-1)$-presentations
under reversion, see Proposition \ref{2.15}. Let us first give an
example showing that reversion does not necessarily preserve
$(-1)$-type.

\bexa\label{concrete} Consider a $\C^*$-surface $V$ and its
equivariant standard completion $(\bV,D)$ as in Proposition
\ref{eboundary}. By virtue of that Proposition,
if $t=2$ in the extended divisor
(\ref{ezigzag}) then
$V$ is special of $(-1)$-type. Passing
to the inverse $\C^*$-action
$\lambda\mapsto\lambda^{-1}$ amounts to interchanging the
divisors $D_+$ and $D_-$, see  \cite{FlZa1}.
According to Proposition \ref{eboundary}(a)
the reverse equivariant completion $(\bV^\vee,D^\vee)$ has an
extended divisor with $t=n$ and
$(F_{n,0}^\vee)^2=1-n$. Thus for $n\ge 3$, $(\bV^\vee,D^\vee)$ is
not of $(-1)$-type, while $(\bV,D)$ is.
\eexa

\bprop\label{2.15}
Given a $(-1)$-presentation $X_n=X(M_2,
c_3,\ldots , c_n,M_n)$, by applying
a finite sequence of elementary
shifts as in (\ref{elmap}) we can transform $X_n$
into a $(-1)$-presentation
$$
X_n^o:=X(M^o_2, c^o_3,\ldots , c^o_n,M^o_n)
$$
such that the reversion of $(X_n^o, D^o)$ is again of
$(-1)$-type.
\eprop

\bproof
Let
$$
X_n^\vee=X(M_n, c^\vee_3,\ldots, c^\vee_{t^\vee}, M_{t},
c^{\vee}_{t^\vee+1}, M_{t-1}, \ldots ,c^\vee_n,M_2)
$$
be the reversion of $X_n$, see Corollary \ref{2.8}. With a suitable
coordinate on $\A^1\cong C_t\backslash C_{t-1}$, the elementary
shift $h_{a,t}$ transforms
$X_n$ into a $(-1)$-presentation
\be\label{change}
X_n'=X(M_2, c_3,\ldots ,c_t, a+M_t, a+c_{t+1}, M_{t+1}',c_{t+2}',
\ldots c'_n,M'_n)\, \ee
with reversion
$$
X_n^{\prime\, \vee}=X(M'_n, c^{\prime\, \vee}_3,\ldots,
c^{\prime\, \vee}_{t^\vee}, a+M_{t}, c^{ \vee}_{t^\vee+1},
M_{t-1}, \ldots ,c^{\vee}_n,M_2)\,,
$$
see Lemmas \ref{2.12} and \ref{2.149}. Choosing $a$
general we may suppose that
$$
c^\vee_{t^\vee+1} \not\in M_t^o:=a+M_t\,.\leqno (*)_t
$$
Applying successively shifts $h_{a_t, t}$, $t=2,\ldots,
n-1$, with general $a_3, \ldots ,a_n\in\C$ the resulting
surface $X_n^o$ satisfies $(*)_t$ for all $t=2,\ldots n-1$. Thus
$X_n^{o\vee}$ is of $(-1)$-type, as required.
\eproof

\bdefi\label{bkw}
Applying an elementary shift $h=h_{a,t}$ to
the reversed presentation $X_n^\vee$ we obtain a presentation
$$
h_{a,t^\vee}^\vee(X_n):= (h_{a,t}(X_n^\vee))^\vee\,,
$$
which we call a
{\it backward elementary shift}. If $X_n^\vee$ is of $(-1)$-type,
then according to Lemmas \ref{2.12} and \ref{2.149}, $h^\vee_{a,t^\vee}$
transforms a presentation
$$
X_n=
X(M_2, c_3, \ldots, c_{t-1},\,M_{t-1},\, c_{t},\, M_{t},\,
c_{t+1},\,M_{t+1},\ldots,c_n,\, M_n)\,
$$
into
\be\label{fbkw} X_n^{\prime}= X(M'_2, c_3^{\prime}, \ldots,
c_{t-1}^{\prime},\,M_{t-1}',\, c'_{t},\, a+ M_{t},\,
c_{t+1},\,M_{t+1},\ldots,c_n,\, M_n)\,.
\ee
Clearly then $X_n^{\prime\vee}$ is as well of $(-1)$-type.
However note that a backward shift can transform a $(-1)$-presentation $X_n$
into one not of $(-1)$-type.
\edefi

In the sequel we fix, for every $t$ in the range $2\le t<n$, an isomorphism
\be\label{isop}
\bdi\alpha_t: C_t\backslash C_{t-1}&\rTo^\cong& \A^1
\quad \mbox{with}\quad
\alpha_t(C_t\cap C_{t+1})=\{0\}
\edi\ee
so that $[M_t]\in \fM^*_{s_t}$ (see \ref{conf*}).

\blem\label{2.20}
Let
$$
X_n=X(M_2, c_3,M_3\ldots , c_n,M_n)\quad\mbox{ and }\quad
X'_n=X(M'_2, c'_3,M'_3\ldots , c'_n,M'_n)
$$
be $(-1)$-presentations of (possibly different)
Gizatullin surfaces $V$ and $V'$,
respectively.  Assume that the reversed presentations $X_n^\vee$ and
$X_n^{\prime\,\vee}$  are also of $(-1)$-type. If the
corresponding invariants
$$
(s_i)_{2\le i\le n},\quad Q(X_n,D)\in\fM \quad\mbox{and}\quad
(s'_i)_{2\le i\le n},\quad Q(X'_n,D')\in\fM
$$
are equal (cf.\ \ref{coin}), then there exists a
presentation
$$
X_n''=X(M''_2, c''_3,M_3'',\ldots ,c_n'',M''_n)\quad
\mbox{of}\quad V
$$
such that
for every $t \geq 2$,\bnum[(a)]\item
$\alpha'_t(M_t')=\alpha_t''( M_t'')$ when fixing  suitable
isomorphisms $\alpha_t':C'_t\backslash C'_{t-1}\to \A^1$ and
$\alpha_t'':C''_t\backslash C''_{t-1}\to\A^1$ as in (\ref{isop});
\item $X_n''$ and its reversion are both of $(-1)$-type.
\enum
\elem

\bproof
By assumption we have
\be\label{haho}
\lambda_i \alpha'_i(M_i')=\alpha_i(M_i) +a_i\quad\mbox{ for some}
\quad \lambda_i\in
\C^*\quad\mbox{and}\quad  a_i\in\C,\quad i =2,\ldots,n-1\,.
\ee
After changing one of the isomorphisms $\alpha_i$, $\alpha'_i$
appropriately
we may assume that $\lambda_i=1$.
Applying a suitable backward shift
$h^\vee_{a,t^\vee}$ we can
translate $M_t$ to $M_t''=M_t+a_t\cong M_t'$.
Under this transformation $c_i$
and $M_i$ remain unchanged for $i>t$, see (\ref{fbkw}).
Moreover the relations (\ref{haho}) remain valid for $i<t$
with possibly new coefficients $\lambda_i$, $a_i$.
Applying decreasing induction starting with $t=n-1$
we can thus achieve that $a_t=0$ for $t=2,\ldots, n-1$,
as required. The transformed presentation
$X_n''$ is then necessarily of $(-1)$-type since
$M''_t\subseteq C''_t\backslash (C''_{t-1}\cup C''_{t+1})$
by construction. Moreover the reversed presentation is as well of
$(-1)$-type since we only applied backwards shifts, proving also (b).
\eproof

\subsection{Isomorphisms of special surfaces of $(-1)$-type}\label{ss3.7}

The following Proposition is the main result of Section 4.

\bprop\label{main}
Two special smooth Gizatullin surfaces $V$
and $V'$ with standard ($-1$)-completions $(\bV,D)$ and
$(\bV',D')$ \footnote{See Definition \ref{basic}.}
are isomorphic if and only if
$D'\cong D$ or $D'\cong D^\vee$ and the configuration invariants
\be\label{config}
\tQ(\bV, D)\quad\mbox{and} \quad \tQ(\bV', D')
\ee
are equal.
\eprop

\bproof
The `only if'  statement follows from Theorem
\ref{matchmain}. To show the converse we note first that according
to Corollary \ref{DG+special} and Lemma \ref{2.7} $\bar V$ and $\bar V'$ admit
$(-1)$-presentations. Using
Proposition \ref{2.15}, applying appropriate elementary shifts we
can achieve that also the reversed completions are of $(-1)$-type.

Since $V$ and $V'$ are special, the standard zigzags of our
completions have the form
$$
D=[[0,0,(-2)_{t-2}, -2-r, (-2)_{n-t}]]\quad\mbox{and} \quad
D'=[[0,0,(-2)_{t'-2}, -2-r', (-2)_{n'-t'}]]\,,
$$
respectively, where by our assumption either $t'=t$ or $t'=t^\vee$.
Replacing $(\bV,
D)$ by the reversion $(\bV^\vee, D^\vee)$, if necessary, we may
restrict to the case where $t'=t$ and the configurations invariants
$Q(\bar V,D)$ and $Q(\bar V',D')$ are equal, see Theorem \ref{matchmain}.

Assume first that $2<t=t'<n$.
As we already remarked,  the completions $(\bV, D)$ and
$(\bV', D')$ arise from $(-1)$-presentations
$$
\bV=X_n=X(M_2, c_3,\ldots, M_t\ldots , c_n, M_n) \,\,
\mbox{and}\,\, \bV'=X_n'=X(M'_2, c'_3,\ldots, M'_t\ldots ,
c'_n, M_n'),
$$
respectively, where
$$
s_2=s_2'=s_n=s_n'=1\,, \quad s_t=s_t'=r,  \quad \mbox{and}\quad
s_i=s_i'=0\quad\forall i\not\in \{2,t,n\}\,.
$$
According to Lemma \ref{2.20} we may
assume  that for every $j=2,\ldots, n-1$ the
configurations $M_j$, $M_j'$ coincide under appropriate
isomorphisms $ C_j\backslash (C_{j-1}\cup
C_{j+1})\cong \C^*\cong C'_j\backslash (C'_{j-1}\cup
C'_{j+1})$. In particular, they are proportional
whatever are these isomorphisms.

Applying now a sequence of elementary shifts,
by Corollary \ref{2.14}
we may suppose that the points $c_3,\ldots, c_n$ and
the unique point of $M_n$ are fixed under the torus
action,\footnote{Possibly after such transformations $X_n^\vee$ is not
any longer of $(-1)$-type; however this does
not matter in the rest of the proof.}
and similarly for $c'_3,\ldots, c'_n$ and $M'_n$.
In particular $c_i=c_i'$ for $ i=3,\ldots,n$ and
$M_n=M_n'$.

After these shifts the configurations $M_i$ and $M'_i$ are contained
in the same curve $C_i=C_i'$.
By our assumption they are proportional for $2\le i<n$
and equal for $i=n$.

Using further the torus action on one of the surfaces we can move
the point of $M_2$ into that of $M_2'$ so that
$M_2=M_2'$. There is a one-parameter subgroup of the torus
acting trivially on $C_2$. It is easily seen that, since
$t>2$, this subgroup acts nontrivially on $C_t$. With this
$\C^*$-action we can move $M_t$ on $C_t$ into the position of
$M_t'$ (that is proportional to $M_t$),
keeping $M_2$ and $M_n$ fixed. Now the presentations
became equal, and so they define isomorphic Gizatullin surfaces
$V$ and $V'$.

The same argument works also in the case, where $t=t'=2<n$ or $2<t=t'=n$.
We leave the details to the reader.
\eproof

To give right away an application let us deduce
Theorem \ref{0.5}(c) in the Introduction
in a particular case.

\bcor\label{2.3}
Let  $V,\,V'$ be smooth Gizatullin
$\C^*$-surfaces with DPD-presentations
$$
V=\Spec \C[u][D_+, D_-]\qquad\mbox{and}\qquad V'
=\Spec \C[u][D'_+, D'_-]\,.
$$
Suppose that  $D_+=D'_+=0$. Then $V$ and $V'$
are isomorphic if and only if
\be\label{cois}
\deg \{D_-\}=\deg \{D_-'\} \quad \mbox{and}\quad
\lfloor D_-\rfloor=\lfloor \beta^*(D'_-)\rfloor
\ee
for some automorphism $\beta$ of $\A^1$.
\ecor

\bproof Since $D_+=0$, the pair $(D_+,D_-)$ is equivalent to
a pair in (\ref{smoothDPD}) with $t=2$.
Hence by Proposition \ref{crit} either $n=2$ in (\ref{smoothDPD})
i.e., $\{D_-\}=0$,
or
$V$ is special. The same holds for $V'$.

Suppose first that (\ref{cois}) is fulfilled.
In the case where $\{D_-\}=0$ the pairs $(D_+,D_-)$ and
$(D'_+,D'_-)$ are equivalent  in the sense of \ref{eph} and so $V$
and $V'$ are equivariantly isomorphic.  For the remaining
part of the proof we assume that $\{D_-\}\ne 0$. Letting $\lfloor
D_-\rfloor=\{p_1,\ldots, p_r\}$, by Corollary
\ref{specialconfig}(1) the configuration invariants of both $V$
and $V'$ are given by $Q(\bV,D)=(p_i)_{0\le i\le r}\in\fM_{r+1}$.
Since $\deg \{D_-\}=\deg \{D_-'\}$ the zigzags of the standard
completions of $V$ and $V'$ are equal, see
Corollary \ref{specialconfig}(2).
Since $t=2$
in (\ref{ezigzag}), both
$V$ and $V'$ arise from $(-1)$-presentations.
Applying Proposition \ref{main}, $V$ and $V'$ are
isomorphic, as required.

Conversely, if $V\cong V'$
then their standard boundary zigzags coincide up to reversion.
Inspecting Proposition \ref{eboundary}
it follows that $\deg \{D_-\}=\deg \{D_-'\}$.
Moreover, by Proposition \ref{main}
the configuration invariants of $V$ and $V'$
are the same. Since they are given by
$\lfloor D_-\rfloor$ and $\lfloor D'_-\rfloor$, respectively,
(\ref{cois}) follows.

\eproof


%


\brem\label{345}
The $\C^*$-surfaces as in Corollary
\ref{2.3}  are just normalizations of the hypersurfaces in $\A^3$ as in
\ref{0099} in the Introduction.
\erem

\section{Shifting presentations and moving coordinates}

In this section we provide further technical
tools that will enable us in the next section
to deduce Theorem \ref{0.5} from the
Introduction (see Theorem \ref{main2} below).
For certain special Gizatullin surfaces of type I
this theorem was already proved in Section 4.
The strategy of proof for the remaining surfaces is similar,
namely to apply Proposition \ref{main}.
However, given a special Gizatullin surface $V$ of type II
with standard completion $(\bV, D)$ usually neither
$(\bV, D)$ nor its reverse completion is of $(-1)$-type,
as it was the case for the special type I surfaces already treated.
For instance, if in the presentation the blowup centers
$c_3,\ldots, c_{k+1}$ are points of a feather generated
by the blowup of $M_2$ then the elementary shifts
are the identity on many of the subsequent components of the zigzag.
It is not enough to shift just
the first moving blowup center,
either in the presentation or in its reversion, as in the
proof of Proposition \ref{main} above.
We have to take into account also second order motions,
which makes the proof considerably
more involved.

\subsection{Coordinate description of a presentation}\label{ss4.1}

As a principal technical tool we  use a sequence
of coordinate charts on our presentation $X_n$.
They appear in the
recursive construction of $X_n$
when describing the blowups in coordinates.
This procedure is similar to the Hamburger-Noether algorithm
for resolving a plane curve singularity, and our
coordinates are analogous to those in the Hamburger-Noether
tables, see e.g., \cite{Ru}. One of the main insights is that they allow
an explicit description
of the correspondence fibration $\psi:W\to \PP^1$ as in Definition
\ref{corr}, see Proposition \ref{111} below.

\bsit\label{4-1}
On the quadric $X_1=Q=\PP^1\times\PP^1$ we consider the affine chart
$$
U_1=Q\backslash (C_0\cup C_1)\cong\A^2\,
$$
with affine coordinates $(x_1,y_1)=(x,y)$, where as before
\be\label{3c} C_0=\{x=\infty\},\quad
C_1=\{y=\infty\},\quad\mbox{and}\quad C_2=\{x=0\}\,. \ee
We let $V=V_n$ be a special smooth Gizatullin surface with a
presentation $X=X(M_2, c_3,\ldots, c_n, M_n)$
obtained from the quadric $Q$ by a sequence of
blowups as in Definition
\ref{basic0}.
We decompose this presentation into a sequence of single blowups
$$
X=X_N\to X_{N-1}\to \ldots \to X_1=Q\,,
$$
so that first we blow up the points of $M_2$ on $C_2$
to create the corresponding feathers (in any order),
then $c_3$ to create $C_3$, then the points on $M_3\subseteq C_3$
to create feathers etc., as prescribed by Definition \ref{basic}.

We say that the blowup $X_i\to X_{i-1}$ is of type (F) if it
creates a feather, say, $F_i$. Otherwise $X_i\to X_{i-1}$
is called of type (C), in which case we let $F_i=\emptyset$.

Starting with the coordinate system $(x_1,y_1)=(x,y)$ on $X_1=Q$
we construct recursively coordinate charts $U_i\cong \A^2$ on the
intermediate surfaces $X_i$ with coordinates $(x_i,y_i)$.
They satisfy the following properties.
\bnum \item[(1)$_i$]
If $C_s$ is the last curve of the zigzag constructed on $X_i$ then
\be\label{U}
U_i=X_i\backslash \left(C_0\cup C_1\cup\ldots \cup C_{s-1}\cup
\bigcup_{j\le i} F^\vee_j\right)\cong \A^2\,. \ee
where $F_j^\vee=\emptyset$  if the blowup $X_j\to X_{j-1}$ is of type (C),
and $F_j^\vee$ is the curve described explicitly below,
when this blowup is of type (F). As we shall show in \ref{mafe1} $F_j^\vee$
is the matching feather
of  $F_j$.

\item[(2)$_i$] $C_s\cap U_i=\{x_i=0\}$ and $y_i|C_s$ yields
an affine coordinate on $C_s\backslash C_{s-1}\cong \A^1$.

\item[(3)$_i$] If $X_i\to X_{i-1}$ is of type (F) then $F_i\cap
U_i=\{y_i=0\}$. \enum

Clearly these properties are satisfied for $U_1$ and $(x_1,y_1)$
when we let $F_1^\vee=\emptyset$.
Assume that the coordinate chart $U_i\subseteq X_i$ with
coordinates $(x_i,y_i)$ was already constructed so that (1)$_i$,
(2)$_i$, (3)$_i$ are satisfied. We introduce the coordinate chart
$U_{i+1}$ on $X_{i+1}$ as follows.

{\bf Type F:} In the next blowup $X_{i+1}\to X_i$ with center at a
point  $p=(0, d)\in M_s\subseteq C_s$ a feather $F_{i+1}$ is
created. Then we let \be\label{star}
(x_{i+1},y_{i+1})=\left(\frac{x_i}{y_i-d}, \
y_i\right)\,\qquad\mbox{so that}\quad
(x_i,y_i)=\left((y_{i+1}-d)x_{i+1},y_{i+1}\right)\,. \ee We let
$F_{i+1}^\vee$ denote the proper transform on $X_{i+1}$ (and on
all further surfaces $X_{i+j+1}$) of the closure in $X_i$ of the
affine line $\{y_i=d\}$. Clearly, the rational function $x_{i+1}$
has a first order pole along the curve $F_{i+1}^\vee$. So
$(x_{i+1},y_{i+1})$ are coordinates in the affine chart $U_{i+1}$
as in (1)$_{i+1}$ with axes as in (2)$_{i+1}$ and (3)$_{i+1}$. By
construction, the level set $\{y_i=d\}$ of the rational function
$y_i$ on $X_{i+1}$ contains both $F_{i+1}$ and $F_{i+1}^\vee$. We
show in Lemma \ref{mafe1} below  that $(F_i,F_i^\vee)$ is actually
a matching pair as in Definition \ref{mafe}.

{\bf Type C:} In the next blowup $X_{i+1}\to X_i$ with center at a
point  $c_{s+1}=(0,c)\in C_s$  the component $C_{s+1}$ is
created. In this case $F_{i+1}=F_{i+1}^\vee=\emptyset$ and
\be\label{starstar}
(x_{i+1},y_{i+1})=\left(x_i,\frac{y_i-c}{x_i}\right)\,\qquad\mbox{so
that}\quad (x_i,y_i)=\left(x_{i+1},x_{i+1}y_{i+1}+c\right)\,. \ee
These are coordinates in the affine chart $U_{i+1}$ as in
(1)$_{i+1}$. The exceptional curve $C_{s+1}$ is given on
$U_{i+1}$ as   $\{x_{i+1}=0\}$. The rational function $x_{i+1}$
has a first order pole along the curve $C_s$, hence $C_s\cap
U_{i+1}=\emptyset$.
\esit

\subsection{Correspondence fibration revisited}
Let us reverse the zigzag $D=C_0\cup\ldots\cup C_s$ on the surface
$X_i$ by a sequence of inner elementary transformations so
that $C_s$ and its previous (new) component
$C^\vee_{s^\vee-1}$ become $0$-curves. On the resulting surface
$W_i$ the component $C_{s-1}$ is blown down, while the coordinate
chart $U_i$ remains unchanged; indeed,
$$
U_i=W_i\setminus \left(E_i\cup\bigcup_{j\le i} F^\vee_j\right)\,,
\qquad\mbox{where}\quad E_i=C^\vee_{n}\cup\ldots\cup
C^\vee_{s^\vee-1}\cup C_s\,.
$$
The linear systems $|C_{s}|$ and
$|C^\vee_{s^\vee-1}|$ define morphisms $p_1:W_i\to \PP^1$ and
$p_2:W_i\to \PP^1$, respectively.
We may suppose that $C_s$ and
$C^\vee_{s^\vee-1}$ are the fibers of $p_1$ and $p_2$,
respectively, over $\infty\in \PP^1$. Here $p_2$ is
the correspondence fibration for the pair of curves $(C^\vee_{s^\vee},C_s)$
(see Definition \ref{corr}). We note that this fibration differs from
the correspondence fibration for $(C_s, C^\vee_{s^\vee})$
just by a sequence of elementary transformations in the fiber
$C^\vee_{s^\vee}$ of $p_2$.
 The following
proposition gives an insight into the geometric meaning of the
above coordinates $(x_i,y_i)$.

\bprop\label{111} In appropriate coordinates on $\A^1=\PP^1
\backslash\{\infty\}$, the maps $p_1$ and $p_2$ restricted to the
chart $U_i$ are given by $x_i$ and $y_i$, respectively.
Furthermore, in the coordinates $(x_i,y_i)$ the curve
$C^\vee_{s^\vee+1}\cap U_i$ is given by equation
$y_i=0$ and so $q^\vee=c^\vee_{s^\vee+1}=(0,\infty)\in
C^\vee_{s^\vee}\backslash C^\vee_{s^\vee-1}$.\footnote{Cf.\
Lemma \ref{mem}(a) and Corollary \ref{claim}. }
\eprop

\bproof We assume by induction that the assertion holds for $W_i$,
and we show that it holds also for $W_{i+1}$. Suppose first that
the blowup $X_{i+1}\to X_i$ is of type (F) so that the feather
$F_{i+1}$ attached to the point $(0,d)\in C_s$ is created.
Reversing the zigzags, we replace $X_i$ by $W_i$ and $X_{i+1}$ by
$W_{i+1}$.

To pass directly from the surface $W_i$ to $W_{i+1}$, we perform
first an elementary transformation on $W_i$ by blowing up with
center at the point $P'\in C^\vee_{s^\vee-1}\cap C^\vee_{s^\vee}$
(this results in a new curve $\tilde C^\vee_{s^\vee-1}$) and
contracting the proper transform of $C^\vee_{s^\vee-1}$. Then we
blow up with center at the point $(0,d)\in C_s$ to
create the feather $F_{i+1}$. On the resulting surface
$W_{i+1}$ we have $C_s^2=0$ and $(\tilde C^\vee_{s^\vee-1})^2=0$.

Let us consider a coordinate chart on the
surface $W_i$ centered at the intersection point $C_{s}\cap
C^\vee_{s^\vee-1}$ with coordinate functions
$$
(u,v)=\left(x_i, \frac{1}{y_i-d}\right)\,
\quad\mbox{and axes}\quad
C_{s}=\{u=0\}\,,\quad C^\vee_{s^\vee-1} =\{v=0\}\,.
$$
Blowing up $P'$ and contracting $C^\vee_{s^\vee-1}$  leads
to  coordinates
$$
(\tilde u, \tilde v)=(uv,v)=
\left(\frac{x_i}{y_i-d},\, \frac{1}{y_i-d}\right) =\left(x_{i+1},
\frac{1}{y_{i+1}-d}\right)
$$
on the surface $W_{i+1}$  (see (\ref{star})) with axes
$$
C_s=\{\tilde u=0\}=\{x_{i+1}=0\}\qquad\mbox{and}\qquad
\tilde C^\vee_{s^\vee-1}=\{\tilde v=0\}=\{y_{i+1}=\infty\}\,.
$$
Now both assertions follow in this case.

Suppose further that the blowup $X_{i+1}\to X_i$ is of type
(C) so that it creates the new component $C_{s+1}$ of the zigzag.
Similarly as before, we introduce first a coordinates chart on
$W_i$ centered at the point $C_s\cap C^\vee_{s^\vee-1}$, at which
$(x_i,y_i)=(0,\infty)$, with coordinates
$$
(u,v)=\left(x_i, \frac{1}{y_i-c}\right)\,.
$$
Next we create the component $C_{s+1}$ by a blowup
at the point $(0,c)\in C_s$ with $(u,v)$-coordinates
$(0,\infty)$. After performing an elementary transformation at
$P'$ as before we get a new surface $\tilde W_i$ and a coordinate
chart on $\tilde W_i$ with coordinate functions
$$
(\tilde u, \tilde v)=(uv,v)=
\left(\frac{x_i}{y_i-c},\frac{1}{y_i-c}\right)\,
\quad\mbox{and axes}\quad
C_s=\{\tilde u=0\}\,,\quad\tilde C_{s^\vee-1}^\vee
=\{\tilde
v=0\}\,.
$$
On  $\tilde W_i$ we have $(\tilde
C_{s^\vee-1}^\vee)^2=C_{s}^2=0$ and $C_{s+1}^2=-1$.
Around the intersection point $C_s\cap C_{s+1}$ on $\tilde
W_i$ this leads to coordinates
$$
(\tilde x, \tilde
y)=\left(\tilde u, \frac{1}{\tilde v}\right)=
\left(\frac{x_i}{y_i-c},y_i-c\right)
\quad\mbox{and axes}\quad
C_{s}=\{\tilde x=0\}\,,\quad C_{s+1} =\{\tilde
y=0\}\,.
$$
To achieve the equality $C_{s+1}^2=0$ \footnote{That is, to shift
two zero weights in the zigzag one position to the right.} we have
to perform a further elementary transformation by blowing up with
center at the point $C_s\cap \tilde C^\vee_{s^\vee-1}$ and
contracting the proper transform of $C_{s}$. On the resulting
surface $W_{i+1}$ this creates a new component
$C_{s^\vee-2}^\vee$ with $(C_{s^\vee-2}^\vee)^2=C_{s+1}^2=0$.
Proceeding in the same way as before we obtain
on $W_{i+1}$ coordinates (see
(\ref{starstar}))
$$
(\hat u,\hat v)=\left(\tilde x\tilde
y, \tilde x\right) =\left(x_i,\frac{x_i}{y_i-c}\right)=\left(x_{i+1},
\frac{1}{y_{i+1}}\right)\,,
$$ centered at the point $C_{s^\vee-2}^\vee\cap C_{s+1}$, with axes
$$C_{s+1}=\{\hat u=0\}\qquad\mbox{ and}\qquad C_{s^\vee-2}^\vee=\{\hat
v=0\}\,,$$ and again the first assertion follows.

Let us finally check that also in this case the intersection point
$C_{s^\vee}^\vee\cap \tilde C_{s^\vee+1}^\vee$ satisfies the
condition $y_{i+1}=0$. The coordinate chart on $W_{i+1}$ around
this intersection point has coordinates
$$
(u',v')=\left(\frac{1}{\tilde u},\tilde v\right)\,
\quad\mbox{with axes}\quad
\tilde C_{s^\vee-1}^\vee=\{v'=0\}\quad\mbox{and}\quad
C_{s^\vee}^\vee =\{u'=0\}\,.
$$
Since $u'=1/\tilde
u=1/\tilde x=y_{i+1}\,,$ the curve $C_{s^\vee}^\vee\cap
U_{i+1}$ is contained in $\{y_{i+1}=0\}$. Now the proof is
completed. \eproof

The next lemma clarifies the meaning of the curves $F_i^\vee$.

\blem\label{mafe1}
If $F_{i+1}$ and $F_{i+1}^\vee$ are nonempty
then they form a  matching pair on the surface $X=X_N$
as in (\ref{mafe}) above. \elem

\bproof
Consider  the correspondence fibration $\psi:W\to\PP^1$
for the pair $(C_s, C^\vee_{s^\vee})$. With the notation as
in the proof of Proposition \ref{111}, $W$ is a blowup of
$W_{i+1}$ with centers  in points of $C_{s}\cap U_{i+1}$
different from $(0,d)$ and in infinitesimally near points. By
construction, the curve $F_{i+1}^\vee$ on $W_{i+1}$ does not meet
$C_s$. Hence the proper transform of $F^\vee_{i+1}$ on $W$
does not meet $D^{\ge t}$. By Lemma \ref{mem}(b)
$F_{i+1}^\vee$ is a feather of $D_\ext^{\vee\ge
t^\vee}$. Since it meets $F_{i+1}$, the pair
$(F_{i+1},F^\vee_{i+1})$ is a  matching pair, as
desired.
\eproof

\subsection{Coordinates on special Gizatullin surfaces}\label{ss4.32}
\bsit\label{4.5}
Let $V$ be a special smooth
Gizatullin surface with data $(n,r,t)$ as in
Definition \ref{0.2}. This means
that  $V$ admits a standard completion $(\bV,D)$ satisfying: \bnum \item[(a)]
every component $C_i$, $i\ge 3$, of the zigzag $D$ is
created by a blowup on $C_{i-1}$; \item[(b)] the extended divisor of
$(\bV,D)$ is of the form
$$
D_\ext=C_0+\ldots +C_n+F_2+F_{t1}+\ldots F_{tr}+F_n,\quad n\ge
3\,,
$$
where $F_2$, $F_{t\rho}$ and $F_{n}$ are feathers with mother
components $C_2$, $C_t$, and $C_n$, respectively.\enum Due to
Corollary \ref{DG+special} there exists a presentation
\be\label{f4.51} \bV=X=X(M_2, c_3, \ldots, c_n, M_n)\, \ee  as in
\ref{basic0} with $V=X\backslash D$.
Reversing the completion $(\bV,D)$, if necessary, we may
assume that either
\bnum[(i)] \item $2<t<n$, $|M_2|=|M_n|=1$,
$|M_t|=r\ge 0$ and $|M_i|=\emptyset$ $\forall i\ne2,t,n$, or

\item $2<t=n$, $|M_2|=1$, $|M_n|=r+1\ge 1$ and $|M_i|=\emptyset$
$\forall i\ne2,n$. \enum
We suppose in the sequel that the feather $F_2$ is
attached to component $C_{k+1}$ with $k\ge 1$. Thus $c_{i+1}
\in C_i\cap F_2$ $\forall i=2,\ldots,k$, while $c_{k+2}\notin
F_2$.\esit

We note that in general neither $(\bV, D)$ nor its
reversion is of $(-1)$-type. Our main
goal  is to transform such a presentation into one of
$(-1)$-type.

Returning to the procedure as in Definition \ref{basic}, we
adopt the coordinates in \ref{4-1} to our case.

\bsit \label{coord}
In what follows we suppose that $k+1\le t$.
To describe coordinate charts of the surface $X$ as in
(\ref{f4.51}) we proceed as follows.
\begin{enumerate}
\item Let $U_1=X_1\setminus (C_0\cup C_1)$ be the affine
coordinate chart on the quadric $X_1=Q$ as in \ref{ss4.1} with
coordinates $(x_1,y_1)$. The feather $F_2$ is created via the
blowup $X_2\to X_1$ of type (F) with center $c_2=(0,0)\in
C_2$ and $F_2^\vee\subseteq X_2$ is the proper transform of
$\{y_1=0\}$. In the affine chart $U_2$ as in (\ref{U}) we have
coordinates
$$
(x_2,y_2)=(x_1/y_1, y_1) \mbox{ with axes }
C_2=\{x_2=0\}  \mbox{ and } F_2=\{y_2=0\}\,.
$$
\item We perform inner blowups at the subsequent intersection
points $c_3=F_2\cap C_2$, $c_4=F_2\cap C_3, \ldots,
c_{k+1}=F_2\cap C_k$ creating the components $C_3,\ldots,C_{k+1}$
of the zigzag.\footnote{In the case $k=1$ this step is
absent.} Since $c_3=
(0,0)$ in  $(x_2,y_2)$-coordinates, its
blowup results in new coordinates $(x_3, y_3)=\left(x_2,
y_2/x_2\right)$. Continuing in this way, in each step $i=2,\ldots,
k$ we obtain the affine chart
\be\label{UU}
U_{i+1}=X_{i+1}\backslash (C_0\cup\ldots \cup C_i\cup F_2^\vee)
\ee on the corresponding surface $X_{i+1}$ with coordinates
\be\label{cordi} (x_{i+1},y_{i+1})=\left(x_i,
y_i/x_i\right)=\left(x/y, y^i/x^{i-1}\right)\,, \ee
the origin $c_{i+2}$ for $i\le k-1$ and axes
$$
C_{i+1}=\{x_{i+1}=0\}\quad \mbox{and}\quad F_2=\{y_{i+1}=0\}\,.
$$
The converse formulas are: \be \label{eq3.2}
(x,y)=(x_{i+1}^iy_{i+1},x_{i+1}^{i-1}y_{i+1}), \quad i=1,\ldots,k\,. \ee

\item For $i=k+1,\ldots, t-1$ \footnote{This step is absent if
$k+1=t$.} we perform outer blowups with center at $c_{i+1}\in
C_i\setminus C_{i-1}$ creating the components
$C_{k+2},\ldots,C_{t}$ of the zigzag. In particular, in
$(x_{k+1},y_{k+1})$-coordinates we have $c_{k+2}=(0,c_{k+2}')$
with $c_{k+2}'\ne 0$. To reduce notation we will identify
$c_{k+2}$ with its coordinate $c_{k+2}'$. Using this convention
also in the following steps, we get  the affine chart
$U_{i+1}\subseteq X_{i+1}$  as in (\ref{UU}) with coordinates
\be\label{cordi2} (x_{i+1},y_{i+1})=\left(x_i,
\frac{y_i-c_{i+1}}{x_i}\right), \quad \mbox{where}\quad
C_{i+1}=\{x_{i+1}=0\}\, \ee
The converse formulas are
\be\label{cordi3} (x_i,y_i)=(x_{i+1}, x_{i+1}y_{i+1}+c_{i+1})\,.
\ee

\item
In particular, for $i=t$ we get the affine chart
$U_t\subseteq X_t$ with coordinates $(x_t,y_t)$. We perform
outer blowups at the distinct points $d_i\hat = (0,d_i)$, $1\le
i\le r$, of $M_t\subseteq C_t\setminus C_{t-1}$ to create the
feathers $F_{ti}$, and a further blowup at $c_{t+1}\hat =
(0, c_{t+1})$ creating the component $C_{t+1}$. On the
resulting surface $X_{t+1}$ this leads to the new coordinate
system\footnote{The enumeration differs from those in \ref{4-1}
since we are performing several blowups in one step.}
\be\label{polt} (x_{t+1},y_{t+1})=\left(\frac{x_t}{P},
\frac{y_t-c_{t+1}}{x_t}P\right)\,, \quad \mbox{where}\quad P
=\prod_{i=1}^r (y_t-d_i)\in\C[y_t]\,, \ee in the affine chart
\[ U_{t+1}=X_{t+1}\backslash
(C_0\cup\ldots\cup C_t\cup F_2^\vee\cup F_{t1}^\vee\cup\ldots\cup
F_{tr}^\vee)\cong\A^2\,,
\]
where $C_{t+1}=\{x_{t+1}=0\}$ and $F_{ti}^\vee$ denotes the
proper transform of the closure in $X_{t+1}$ of the affine line
$\{y_i=d_i\}\subseteq U_t$, $i=1,\ldots,r$. If $c_{t+1}\notin
M_t$ then $F_{ti}=\{x_{t+1}y_{t+1}=d_i-c_{t+1}\}$, while for
$d_i=c_{i+1}$ the feather $F_{ti}$ is given by $y_{t+1}=0$.
The converse formulas are
\be\label{polt1}
(x_t,y_t)=\left(x_{t+1}P_{t+1} ,\, x_{t+1}y_{t+1}+c'_{t+1}\right),
\,\,\,\,\mbox{where}\,\,\,\, P=\prod_{j=1}^r
(x_{t+1}y_{t+1}+c_{t+1}-d_i)\,.
\ee

\item For $i=t+1,\ldots, n-1$ we blow up at the point $c_{i+1}\in
C_i\setminus C_{i-1}$ creating the component $C_{i+1}$, while for
$i=n$ we blowup at the point $c_{n+1}\in M_n$ creating the
feather $F_n$. In step $i$ with $t+1\le i\le n-1$ we obtain
coordinates $(x_{i+1}, y_{i+1})$ by formula (\ref{cordi2}),
where $c_{i+1}\hat=(0,c_{i+1})$ in $(x_i, y_i)$-coordinates.
The converse formulas (\ref{cordi3}) are still available.
Obviously $(x_{i+1},y_{i+1})$ forms a coordinate system in the
affine chart
$$
U_{i+1}=X_{i+1}\backslash (C_0\cup\ldots\cup
C_i\cup F_2^\vee\cup F_{t1}^\vee\cup\ldots\cup F_{tr}^\vee)\cong\A^2
$$
on the corresponding surface $X_{i+1}$, where $t+1\le i\le n$.
Similarly for $i=n$ the same formulas define coordinates
$(x_{n+1},y_{n+1})$ in the affine chart $U_{n+1}$ on the
terminal surface $X=X_{n+1}$, where
$F_n=\{x_{n+1}=0\}$.
\end{enumerate}
\esit

\subsection{Moving coordinates}\label{ss4.3} In this subsection
we study the effect of an elementary shift\footnote{This
corresponds to the elementary shift $h_{a,m+3}$ from Section 3.}
\be\label{mshift}
h=h_{a,m}: (x,y)\longmapsto (\xi(x,y),\
\eta(x,y))= (x,\ y+ax^{1+m})\,,
\ee where $m\ge 0$ and $a\in\C_+$,
on the sequence of coordinate systems as in \ref{coord}.
According to \ref{2.9}, $h$ transforms the presentation $X$ as
in (\ref{f4.51}) into $X'=h_*(X)$, and yields an isomorphism
of the affine surfaces $V$ onto a new one
$V'=h_*(V)$. The standard completion can change due to the presence of
indeterminacy points of $h$ on $C_0$.
We introduce the following coordinates on $h_*(X)$.

Letting $(\xi_1,\eta_1)=(\xi,\eta)$ and performing the sequence of
blowups as described in \ref{4.5}
in the images of the centers under $h$,
we obtain a new sequence of coordinates
$$
(\xi_i(x_i,y_i),\eta_i(x_i,y_i)), \quad i=1,\ldots, n+1\,.
$$
Our next aim is to give, for every $i\ge 1$, explicit
expressions of the maps $(x_i,y_i)\longmapsto
(\xi_i(x_i,y_i),\,\eta_i(x_i,y_i))$ and of the vector field
$$
\left(\bar\xi_i,\bar\eta_i\right)=\left(\frac{\p\xi_i}{\p a}(0),\,
\frac{\p\eta_i}{\p a}(0)\right)\,.
$$

\bsit\label{recursion}
According to the cases in
\ref{coord}, the following recursive formulas hold.
\bnum[(1)]
\item In the blowup $X_{2}\to X_1$ we have
$(\xi_2,\eta_2)=(\xi_1/\eta_1, \,\eta_1)$. \item  In the inner
blowups $X_{i+1}\to X_i$, $i=1,\ldots, k$, of \ref{4.5}(2) we have
in view of (\ref{cordi})
\[
\left(\xi_{i+1},\eta_{i+1}\right)
=\left(\xi_i,\frac{\eta_i}{\xi_i}\right)
=\left(\frac{\xi}{\eta},\, \frac{\eta^i}{\xi^{i-1}}\right)
=\left(\frac{x}{y+ax^{1+m}},\frac{(y+ax^{1+m})^i}{x^{i-1}}\right)\,.
\]
Using (\ref{eq3.2}) we get
\be\label{rec4} \left(\xi_{i+1}(x_{i+1},
y_{i+1}),\eta_{i+1}(x_{i+1}, y_{i+1})\right)=
\left(\frac{x_{i+1}}{1+ax_{i+1}^{im+1}y_{i+1}^m}
,(1+ax_{i+1}^{im+1}y_{i+1}^m)^i y_{i+1}\right)\,. \ee
\item  In
the blowups $X_{i+1}\to X_i$, $i=k+1,\ldots,t-1$  with centers at
$c_{i+1}\hat=(0,c_{i+1})\in C_i$ we get
\be\label{rec6}
(\xi_{i+1},\eta_{i+1})=\left(\xi_i,\frac{\eta_i-\gamma_{i+1}}{\xi_i}\right),
\quad\mbox{where}\quad \gamma_{i+1}=\eta_i(0,c_{i+1})\,.
\ee
In the same range $i=k+1,\ldots,t-1$ we obtain by
(\ref{cordi3})
\bea \label{rec7}
\xi_{i+1}(x_{i+1},
y_{i+1})&=& \xi_i(x_{i+1},x_{i+1}y_{i+1}+c_{i+1})
\quad\mbox{and}\\
\label{rec8}
\eta_{i+1}(x_{i+1}, y_{i+1})&=&
\frac{\eta_i(x_{i+1},\, x_{i+1}y_{i+1}+c_{i+1})-\gamma_{i+1}}
{\xi_i(x_{i+1},\, x_{i+1}y_{i+1}+c_{i+1})}\,\,.
\eea
Moreover
\be\label{rec9}
\bar\xi_{i+1}=
\bar\xi_i(x_{i+1},x_{i+1}y_{i+1}+c_{i+1})\quad \mbox{and}
\ee
\be\label{rec10}
\bar\eta_{i+1}= \frac{\bar\eta_i(x_{i+1},\,
x_{i+1}y_{i+1}+c)-\bar\eta_{i}(0,c_{i+1})} {x_{i+1}}
-\frac{y_{i+1}}{x_{i+1}}\bar\xi_i(x_{i+1},\,
x_{i+1}y_{i+1}+c_{i+1})\,\,.
\ee

\item Similarly as in (\ref{polt}), with
$\gamma_{t+1}=\eta_t(0,c_{t+1})$ and
$\delta_\rho=\eta_t(0,d_\rho)$ we have
\be\label{poltg}
(\xi_{t+1},\eta_{t+1}) =\left(\frac{\xi_t}{\Pi},
\frac{\eta_t-\gamma_{t+1}}{\xi_t}\Pi
\right),\quad\mbox{where}\quad
\Pi=\prod_{\rho=1}^r(\eta_t-\delta_\rho)\,. \ee The converse
formulas are \be\label{polt1g} (\xi_t,\eta_t)=
\left(\xi_{t+1}\Pi,\xi_{t+1}\eta_{t+1}+\gamma_{t+1}\right),
\quad\mbox{where}\quad
\Pi=\prod_{\rho=1}^r(\xi_{t+1}\eta_{t+1}
+\gamma_{t+1}-\delta_\rho) \,.\ee

\item If $i>t$, the recursion formulas for $(\xi_{i+1}, \eta_{i+1})$
are the same as in step (3).
\enum
\esit

\brems\label{4.6} 1. In all cases $\xi_i=0$ is a local equation for $C_i$ and so
\be\label{barin} \xi_i(0,y_i)=\bar\xi_i(0,y_i)=0\,. \ee
Furthermore,
$$
\p_1\xi_i(0,y_i)=1 \quad\mbox{and}\quad \p_1\bar\xi_i(0,y_i)=0\,,
$$ where $\p_1=\frac{\p}{\p x_i}$ and $\p_2=\frac{\p}{\p y_i}$.
Indeed,  $\p_1\xi_i(0,y_i)\ne 0$
$\forall y_i\in C_i\backslash C_{i-1}\cong\A^1$,
$\forall a\in \A^1$. Hence $\p_1\xi_i(0,y_i)$ is  a
constant equal to its value at $a=0$, which equals 1. The
second relation follows by differentiating the first one.

2. The map $h_i:y_i\mapsto \eta_i(0, y_i)$ yields a translation
$$
\eta_i(0, y_i)=y_i+e_i(a)\quad
\mbox{with}\quad e_i(a)\in \C[a]\,, e_i(0)=0,
$$
of the affine line $C_i\backslash C_{i-1}\cong \A^1$. Indeed,
$\eta_i(0,y_i)$ is an automorphism of $\A^1$ hence it has the form
$c(a)y_i+e_i(a)$, where $c(a),\,e(a)\in \C[a].$ Here $c(a)=c(0)=1$
is constant since it does not vanish for all $a$. Moreover $e_i(0)=0$
since $h_{0,m}=\id$. In particular
\be\label{barin2} \p_2\eta_i(0, y_i)=1 \quad\mbox{and}\quad
\p_2\bar \eta_i(0,y_i)=0\,. \ee \erems

We can now deduce
the following formulas.

\bprop\label{4.7} (a) At step $k+1$ we have
$$
\ba{l}\ds
\xi_{k+1}(x_{k+1},y_{k+1})=
\frac{x_{k+1}}{1+ax_{k+1}^{km+1}y_{k+1}^m}\,,\\
\eta_{k+1}(x_{k+1},y_{k+1})=(1+ax_{k+1}^{km+1}y_{k+1}^m)^k y_{k+1}\,,\\
\bar\xi_{k+1}(x_{k+1},y_{k+1})=-x_{k+1}^{km+2} y_{k+1}^m\,,\\
\bar\eta_{k+1}(x_{k+1},y_{k+1})=kx_{k+1}^{km+1} y_{k+1}^{m+1}\,.\\
\ea
$$
(b) If $k+2\le t$ then we have at step $k+2$
$$
\ba{l}
\bar\xi_{k+2}(x_{k+2},y_{k+2})=-c^{m}x_{k+2}^{s-k} +\hot\,,\\
\bar\eta_{k+2}(x_{k+2},y_{k+2})=kc^{m+1}x_{k+2}^{s-k-2}
+(s-1)c^mx_{k+2}^{s-k-1}y_{k+2}
+\hot\,,\\
\ea
$$
where $\hot$ stands for higher order terms in the first variable
$x_{k+2}$ and
\be\label{s} s=k(m+1)+2\ge k+2\,,\qquad
c=c_{k+2}\ne 0 \,.\ee

\smallskip

(c) With $c$ and $s$ as in (\ref{s}), in the range
$k+3\le i\le\min\{s,t\}$ we have
$$
\ba{l}\ds
\bar\xi_{i}(x_{i},y_{i})=-c^mx_{i}^{s-k}+\hot\qquad\mbox{and}\\
\bar\eta_{i}(x_{i},y_{i})=kc^{m+1}x_{i}^{s-i} +(s-1)
c_{k+3}c^mx_{i}^{s-i+1} +\hot\, \ea
$$
\eprop

\bproof
The first two formulas in (a) are contained in (\ref{rec7})
while the remaining two are obtained by
differentiation.

Inserting $\bar\xi_{k+1}$ and $\bar\eta_{k+1}$ as in (a) into the
recursive formulas (\ref{rec9})-(\ref{rec10}) we get \bea
\label{eq3.5a} \bar\xi_{k+2}(x_{k+2}, y_{k+2}) &=&
-x_{k+2}^{km+2}(x_{k+2}y_{k+2}+c)^m
\qquad\mbox{and}\\
\label{eq3.5b} \bar \eta_{k+2}(x_{k+2}, y_{k+2}) &=&
x_{k+2}^{km}(kc+(k+1)x_{k+2}y_{k+2})(c+x_{k+2}y_{k+2})^m \,.\eea
Now the Taylor expansion
implies (b).

Finally, starting from (b) and using again the recursion formulas
(\ref{rec9})-(\ref{rec10}), (c) follows by an easy computation.
We leave the details to the reader.
\eproof

Proposition \ref{4.7} implies the following
corollary.

\bcor\label{m=0} Suppose that $m=0$ so that $s=k+2$.

(a) On component $C_{k+1}$
we have
\[
\ba{l}\ds
\xi_{k+1}(x_{k+1},y_{k+1})=\frac{x_{k+1}}{1+ax_{k+1}}\,,\\
\eta_{k+1}(x_{k+1},y_{k+1})=(1+ax_{k+1})^{k} y_{k+1}\,.\\
\ea
\]

(b) For every $i=2,\ldots,t$ we have
$$
\ba{l}\ds
\bar\xi_{i}(x_{i},y_{i})=-x_{i}^{2}\,,\\
\bar\eta_{i}(x_{i},y_{i})=(i-1)x_{i}y_i +(i-2)c_{i}\,.\\
\ea
$$
In particular, if $s\le t$ then $\bar\eta_{i}(0,y_{i})=0$ for
$i=2,\ldots,s-1$, while $\bar\eta_{s}(0,y_{s})\neq 0$.
\ecor

\bproof
(a) is just a specialization of Proposition  \ref{4.7} to
the case where $m=0$. The first two formulas in (b) follow in the
case $2\le i\le k+2$ from (\ref{rec4}) by differentiation.
In the general case one can proceed by recursion using
(\ref{rec7}) and (\ref{rec8}). The last assertion holds
since by our assumptions $c_2=\ldots=c_{k+1}=0$ while $c_{k+2}\neq
0$, see \ref{coord}(2), (3).
\eproof

\subsection{Induced motions}
Let us study the map induced by the shift
$$
h=h_{a,m}: (x,y)\mapsto (x, y+ax^{1+m})
$$
on the component $C_i=\{x_i=0\}$. By Remark \ref{4.6}(2) for every
$i=1,\ldots,n$,
$$
h=h|C_i:(0,y_i)\longmapsto(\xi_i(0,y_i),\eta_i(0,y_i))=(0,y_i+e_i(a))
\quad\mbox{ for some}\quad e_i(a)\in\C[a]
$$
with $e_i(0)=0$. We say that $h$
{\em generates a motion on component} $C_i$ if $\deg e_i> 0$.
For instance, using (\ref{rec4}) we obtain
$\eta_i(0,y_i)=0$ for $i=2, \ldots , k+1$, so there is no motion
on components $C_2, \ldots, C_{k+1}$.

In the next lemma we study the motions in the case $m=0$.

\blem\label{4.9} Let $X=X_{n+1}$ be as in \ref{coord}.
Assume as
before that $k+1\le t$. Then $\deg e_q=q-k-1>0$ and $e_q(0)=0$
$\forall q\in [k+2,t]$. In particular, the shift $h=h_{a,0}$
induces a motion on component $C_q$ for every $q$ in this
range.  \elem

\bproof We let $\tilde x=x_2$ and $\xi=\xi_2$. At each step
$i=2,\ldots,t$ in \ref{coord}, only components of the zigzag are
created. Hence we have \be\label{dg0} x_i=\tilde x,\qquad
\xi_i=\xi=\frac{\tilde x}{1+a \tilde x},\quad\mbox{and,
conversely,}\qquad \tilde x=\frac{\xi}{1-a\xi}\,.\ee Furthermore,
by \ref{coord}(3) for $k+2\le q\le t$ \be\label{dg3} \ba{l}\ds
y_q=\frac{y_{q-1}-c_{q}}{\tilde x}=\frac{y_{q-2}-c_{q-1} -c_{q}
\tilde x}{\tilde x^2}=\ldots =\frac{y_{k+1}-\sum_{j=k+2}^q c_{j}
\tilde x^{j-k-2}}{\tilde
x^{q-k-1}}\,\,\\
\quad=\tilde x^{k+1-q}y_{k+1}-\sum_{j=k+2}^q c_{j} \tilde
x^{j-q-1}\,.\ea\ee Conversely, \be\label{dg5} y_{k+1}=\tilde
x^{q-k-1}y_{q} + \sum_{j=k+2}^q c_{j} \tilde x^{j-k-2}\,. \ee
Similarly as in (\ref{dg3}), letting
$\gamma_{j+1}=\eta_j(0,c_{j+1})$ we obtain \be\label{dg6}
\eta_q(\tilde x,y_q)=\xi^{k+1-q}\eta_{k+1}(\tilde
x,y_{k+1})-\sum_{j=k+2}^q \gamma_{j}\xi^{j-q-1}\,,
\ee
where, by virtue of Corollary \ref{m=0}(a) and (\ref{dg5}),
\be\label{dg2} \eta_{k+1}=(1+a \tilde x)^ky_{k+1}=(1+a \tilde
x)^k\left(\tilde x^{q-k-1}y_{q} + \sum_{j=k+2}^q c_{j} \tilde
x^{j-k-2}\right)\,. \ee
According to (\ref{dg0}) we have $\tilde x=\xi(1-a\xi)^{-1}$
and so $1+a\tilde x=(1-a\xi)^{-1}$. Inserting (\ref{dg2}) into (\ref{dg6}) and
replacing $\tilde x$ by $\xi$ leads to
\be\ba{l}\ds\label{dg71}
\eta_q(\tilde x,y_q)=\xi^{k+1-q}(1+a \tilde
x)^ky_{k+1}-\sum_{j=k+2}^q
\gamma_{j}\xi^{j-q-1}\\
=(1-a\xi)^{1-q}y_q+\sum_{j=k+2}^q
c_{j}(1-a\xi)^{2-j}\xi^{j-q-1}-\sum_{j=k+2}^q
\gamma_{j}\xi^{j-q-1}\,\,.\ea\ee Regarding $\eta_q(\tilde x,y_q)$
as a function of $(\xi,\, y_q)$ we need to compute the constant
term $e_q(a)=\eta_q(0,0)$ in its Laurent expansion. We note that
this function is regular in a neighborhood of the axis $\{\tilde
x=0\}=\{\xi=0\}$. Therefore we can omit the first term and the
last sum on the right in (\ref{dg71}). Using the binomial
expansion the remaining sum can be written as the Laurent series
\be\label{dg81} \sum_{j=k+2}^q c_j\sum_{\mu=0}^\infty {2-j\choose
\mu}(-1)^\mu a^\mu \xi^{j-q-1+\mu}\,\,. \ee Its constant term is
\be\label{dg72} e_q(a)=\eta_q(0,0)=\sum_{j=k+2}^q {2-j\choose
q+1-j}(-1)^{q+1-j}c_j a^{q+1-j}\,\,. \ee Since $c_{k+2}\ne 0$, see
\ref{coord}(3), this is a polynomial of degree $q-k-1$ in $a$,
while $e_q(0)=0$ by \ref{4.6}(2),  as stated. \eproof

\bcor\label{DGcase}
If $V$ is a Danilov-Gizatullin surface i.e.\ $r=0$ in
\ref{coord}, then in the range $k+2\le q\le n$ the shift
$h=h_{a,0}$ induces a translation $y_q\mapsto y_q+e_q(a)$ on
component $C_q$, where $\deg e_q=q-k-1$ and $e_q(0)=0$.
\ecor

In the general case, where $m\ge 0$ is arbitrary and $X$ is
a standard completion of a special smooth Gizatullin surface as in
\ref{coord}, we have the following lemma.

\blem\label{4.11}
For every  $i\ge 2$ with $i\ne t$ the
degrees of $e_{i+1}$ and $\p_1\eta_i(0,c_{i+1})$ as polynomials in
$a$ are equal. Moreover the same is true for $i=t$ provided
that  $c_{t+1}\notin M_t$ i.e., $c_{t+1}$ is not a root of $P$.
\elem

\bproof We  give the proof only in the case $i=t$;
in all other cases the same proof applies by replacing $t$ by $i$
and $P$ by $1$. According to (\ref{polt1g}) and (\ref{polt1})
\be\label{4.11.2} \eta_t(x_{t+1}P_{t+1},x_{t+1}y_{t+1}+c_{t+1})=
\xi_{t+1}(x_{t+1},y_{t+1})\cdot
\eta_{t+1}(x_{t+1},y_{t+1})+\gamma_{t+1}\,. \ee By virtue of
Remark \ref{4.6} we have the relations:
$$
\xi_{t+1}(0,y_{t+1})=0,\qquad
\p_1\xi_{t+1}(0,y_{t+1})=1,\quad\mbox{and}\quad
\p_2\eta_t(0,c_{t+1})=1\,.$$ Thus differentiating (\ref{4.11.2}) with
respect to $x_{t+1}$ and evaluating at $x_{t+1}=0$ yields
$$
\p_1\eta_t(0,c_{t+1})\cdot P(c_{t+1})+
y_{t+1} =\eta_{t+1}(0,y_{t+1})\,.
$$
By definition the term on the right is $y_{t+1}+e_{t+1}(a)$. Since
$c_{t+1}\notin M_t$ is not a root of $P$, $P(c_{t+1})\neq 0$,
hence the result follows. \eproof

Returning to the special case $m=0$ (that is, $h=h_{a,0}$)
treated in Lemma \ref{4.9}, let us  show that its conclusion can
be equally applied to the polynomial $e_{t+1}$.

\bprop\label{4.12} Let as before $k+1\le t$. If
$c_{t+1}\notin M_t$  and $m=0$, then $\deg_a e_{t+1}(a) =t-k>0$
and $e_{t+1}(0)=0$. \eprop

\bproof
Let again $r=0$ so that no feather is attached to
component $C_t$ and $V=V^{DG}$ is a Danilov-Gizatullin surface.
The construction of \ref{coord} yields coordinates $(x_i^{DG},
y^{DG}_i)$ and $(\xi_i^{DG}, \eta_i^{DG})$ on
$X^{DG}=\bV^{DG}$. Let us compare these coordinates with
those $(x_i,y_i)$ and $(\xi_i, \eta_i)$
on a general surface $X=\bV$, where $r$ is arbitrary. Since up
to step $t$ in \ref{coord} the constructions in the
Danilov-Gizatullin case and in the general case are identical,
we have
$$
(x^{DG}_t,y^{DG}_t)=(x_t,y_t)
\quad\mbox{and}\quad
(\xi_t^{DG}, \eta_t^{DG})=(\xi_t,\eta_t).
$$
By Lemma \ref{4.11} the degrees of the polynomials $e_{t+1}(a)$ and
$\p_1 \eta_t(0, c_{t+1})$ are equal. Moreover
$$
\p_1 \eta_t(0, c_{t+1})= \p_1 \eta^{DG}_t(0, c_{t+1})\,.
$$
Applying Lemma \ref{4.11} and Corollary \ref{DGcase}, the term on
the right is a polynomial of degree $t-k$ with zero constant
term, which gives the required result.
 \eproof

\subsection{Component of first motion}\label{ss4.5}
We continue to study
an elementary shift
$$
h=h_{a,m}:(x,y)\mapsto (x, y+a x^{1+m})\,
$$
expressed in terms of moving coordinates $(\xi_i,\eta_i)$.

\bprop\label{motion} (Motion Lemma) Suppose that $k+1\le t$ and
$c_{t+1}\notin M_t$. Let as before $s=k(m+1)+2$. Then the
following hold. \bnum[\rm(a)] \item The first motion under
$h_{a,m}$ occurs on component $C_s$ i.e., $e_i=0$ for all $i<s$,
while $\deg e_s>0$.

\item If $k+2\le t$ then after a
general coordinate change
 $(x,y)\to (x,y+bx)$ there is a motion on component
$C_{s+1}$ i.e., $\deg e_{s+1}> 0$. \enum \eprop

The proof of Proposition \ref{motion} is based on the following
lemma.

\blem\label{mainlem} If  $k+1\le t$ and $c_{t+1}\notin M_t$ then
the following hold.
\bnum[\rm(a)] \item $\bar\eta_i(0,y_i)= 0$ for
all $i=1, \ldots, s-1$ while
$\bar\eta_s(0,y_s)=\bar\eta_s(0,0)\ne 0$.

\item If $k+2\le
t$ then after a general coordinate change
 $(x,y)\to (x,y+bx)$ we have
$\p_1\bar\eta_{s}(0,c_{s+1})\ne 0$.
\enum \elem

\bproof[Proof of Proposition \ref{motion}] Lemma \ref{mainlem}(a)
implies that $e_s'(0)=\bar\eta_s(0,0)=\bar\eta_s(0,y_s)\ne
0$ and so there is a motion on the component $C_s$. For every
$i<s$ by Lemma \ref{mainlem}(a), $\bar\eta_i(0,y_i)=0$. This
remains true after a coordinate change $(x,y)\to (x,
y+a'x^{1+m})$, which replaces $a$ by $a+a'$. Consequently,
$e_i'(a)=\frac{\p\eta_i}{\p a}(0,0)=0$ for every $a$ and so
$e_i(a)=e_i(0)=0$. This proves (a).

To deduce (b) we note that by Lemma \ref{mainlem}(b),
$\p_1\bar\eta_s(0,c_{s+1})\ne 0$ and so $\p_1\eta_s(0,c_{s+1})$
is not constant in $a$.
In view of Lemma \ref{4.11} also $e_{s+1}$ is not constant. Now
the assertion follows. \eproof

The rest of this subsection is devoted to the proof of Lemma
\ref{mainlem}.

\bproof[Proof of Lemma \ref{mainlem}]
By (\ref{rec4}) and
Proposition \ref{4.7}, (a)
is true if $s\le t$. Let us deduce (b) for
$s\le t$.
The latter inequality implies that $k+2\le t$. Thus by Lemma
\ref{4.9} a coordinate change $(x,y)\to (x,y+bx)$ results in a
non-trivial translation on component $C_{k+2}$.
Applying such a translation we may assume that $c_{k+3}\ne 0$.
Now $\partial_1\eta_s(0,c_{s+1})\ne 0$ by the second formula in
\ref{4.7}(c), as required.

We assume in the sequel that $s>t$. We have to distinguish several
cases. First we treat in \ref{4.15} the case $k+3\le t$. The
proof in the remaining cases where $k+2=t$ or $k+1=t$ is given
in \ref{4.16}. \eproof

\bsit{\bf The case $k+3\le t$}.\label{4.15} According to
Proposition \ref{4.7}(c), in step $t$

\bea\label{compare}
\bar\xi_{t}(x_{t},y_{t})&=&-c^mx_{t}^{s-k}+\hot \quad{\rm
and}
\\
\label{compare1}
\bar\eta_{t}(x_{t},y_{t})&=&
\alpha x_{t}^{s-t} +\beta c_{k+3} x_{t}^{s-t+1}+\hot\,,\quad{\rm where} \\
\label{compare01}
c=c_{k+2},&&\alpha=kc^{m+1} \quad\mbox{and}\quad
\beta=(s-1)c^m
\eea
are nonzero constants. Let us compute the vector field
$(\bar\xi_{t+1},\bar\eta_{t+1})$.
By (\ref{poltg}),
\be\label{compare4}
(\xi_{t+1},\eta_{t+1})  =\left(
\frac{\xi_t}{\Pi}, \frac{(\eta_t-\eta_t(0,c_{t+1}))\Pi}{\xi_t}
\right), \quad\mbox{where}\quad
\Pi=\prod_{\rho=1}^r(\eta_t-\eta_t(0,d_\rho))\,.
\ee
Since by our
assumption $s-t>0$, by virtue of (\ref{compare1})
$\bar\eta_t(0,d_\rho)=0$ $\forall \rho=1,\ldots,r$. Consequently
$\frac{d}{da} (\Pi)|_{a=0}=P'\cdot\bar\eta_t$,
where $P=P(y_t)$ is as in
(\ref{polt}). Applying in (\ref{compare4}) the derivation
$\frac{d}{da}|_{a=0}$ and expressing $x_t,y_t$  by
 $x_{t+1},y_{t+1}$ as in  (\ref{polt1}) we get
\be\label{compare5}
\ba{rcl}
\bar\xi_{t+1} &=&
\ds \frac{\bar\xi_t}{P(y_t)}-\frac{x_t}{P^2(y_t)}
P'(y_t)\cdot \bar\eta_t \\[8pt]
&=&
\ds\frac{\bar\xi_t}{P_{t+1}}-\frac{x_{t+1}}{P_{t+1}}
P'_{t+1}\cdot (\alpha x_{t+1}^{s-t}P_{t+1}^{s-t}+\hot)\,,
\ea \ee
where
\[
P_{t+1}:=P(c_{t+1}+x_{t+1}y_{t+1})
\quad\mbox{and}\quad P'_{t+1}:=P'(c_{t+1}+x_{t+1}y_{t+1})\,.
\]
The first order Taylor expansion of $P_{t+1}^{s-t}$ is
\be\label{compare6} P_{t+1}^{s-t}=P(c_{t+1})^{s-t}+
(s-t)P'(c_{t+1})P(c_{t+1})^{s-t-1}x_{t+1}y_{t+1}+\hot\, \ee To
compute the lowest order term of $\bar \xi_{t+1}$, the first term
on the right of (\ref{compare5}) is irrelevant by (\ref{compare}).
Therefore from (\ref{compare5}) and (\ref{compare6}) we obtain
\be\label{eqxi} \bar\xi_{t+1}= \gamma x_{t+1}^{s-t+1}+\hot\,,
\qquad\mbox{ where} \quad\gamma= -\alpha
P'(c_{t+1})P^{s-t-1}(c_{t+1})\,. \ee
Likewise we can
differentiate the expression for $\eta_{t+1}$ in (\ref{compare4})
with respect to $a$,  then replace $x_t,y_t$ by
$x_{t+1},y_{t+1}$ according to (\ref{polt1}), and finally use
(\ref{compare}), (\ref{compare1}) to obtain
\be\label{eqcompare}
\ba{l}\ds \bar\eta_{t+1} =
\frac{\bar\eta_t}{x_t}P(y_t)+\frac{y_t-c_{t+1}}{x_t}P'(y_t)\bar\eta_t
-\frac{y_t-c_{t+1}}{x_t^2}P(y_t)\bar\xi_t\\[8pt]
\ds \qquad=  \alpha x_{t+1}^{s-t-1} P_{t+1}^{s-t} + \beta c_{k+3}
x_{t+1}^{s-t}P_{t+1}^{s-t+1} +
y_{t+1}\frac{P'_{t+1}}{P_{t+1}}\bar\eta_t+\hot\, \ea
\ee
Inserting (\ref{compare1}) and the Taylor expansion (\ref{compare6}) into
this formula yields \be\label{compare7} \ba{rcl} \bar\eta_{t+1}&=&
\ds\alpha P(c_{t+1})^{s-t}x_{t+1}^{s-t-1} +(s-t)\alpha
P'(c_{t+1})P(c_{t+1})^{s-t-1}y_{t+1}x_{t+1}^{s-t}
\\
&& + \beta c_{k+3} P(c_{t+1})^{s-t+1}x_{t+1}^{s-t}
+ \alpha P'(c_{t+1})P(c_{t+1})^{s-t-1}y_{t+1}
x_{t+1}^{s-t}+\hot\,
\\&=&
\tilde\alpha x_{t+1}^{s-t-1}+ (\tilde\beta c_{k+3}  +
\tilde\gamma y_{t+1}) x_{t+1}^{s-t}+\hot\,, \ea \ee
where
\be\label{compare8}
\tilde\alpha= \alpha P(c_{t+1})^{s-t},\quad\tilde \beta= \beta
P(c_{t+1})^{s-t+1}\,,\quad\mbox{and} \quad
\tilde\gamma=(s-t+1)\frac{P'(c_{t+1})}{P(c_{t+1})}\tilde\alpha\,
\ee are constants with $\tilde \alpha$, $\tilde\beta\ne 0$.

In the range $t+1\le j<s$ we have $x_{j+1}=x_j$ and
$\xi_{j+1}=\xi_j$. In view of (\ref{eqxi}) this yields
\be\label{compare9} \bar\xi_{j+1}=-\gamma x_{j+1}^{s-t+1}+\hot
\quad \mbox{for}\quad k+1\le j<s\,. \ee
To compute
$\bar\eta_j$ for $j\ge t+2$ we first consider the step from
$t+1$ to $t+2$. Differentiating the recursion formula (\ref{rec6})
for $\bar\eta_{t+2}$ we obtain
\[
\bar\eta_{t+2}=
\frac{\bar\eta_{t+1}}{x_{t+1}}-\frac{y_{t+1}-c_{t+2}}
{x_{t+1}^2}\bar\xi_{t+1}\,.
\]
Using (\ref{eqxi}) and (\ref{compare7}) and replacing $x_{t+1}, y_{t+1}$
by $ x_{t+2}, y_{t+2}$ as in (\ref{cordi3}) we get
\[
\bar\eta_{t+2}=
\tilde\alpha x_{t+2}^{s-t-2}+(\tilde\beta c_{k+3}  + \tilde\gamma
c_{t+2}) x_{t+2}^{s-t-1}+\hot\,
\]
Recursively the same arguments yield
\be\label{coeff00}
\bar\eta_{j}= \tilde\alpha x_{j}^{s-j}+(\tilde\beta c_{k+3}+
\tilde\gamma c_{t+2}) x_{j}^{s-j+1}+\hot \quad\mbox{for}\quad
t+2\le j\le s\,.
\ee
Using (\ref{coeff00}) in the case $s\ge t+2$
and (\ref{compare7}) in the case $s=t+1$, assertion (a) of Lemma
\ref{mainlem} follows.

Let us show part (b) of the lemma.  For $j=s\ge t+2$ (\ref{coeff00})
yields
\be\label{coeff}
\p_1\bar\eta_{s}(0,c_{s+1}) =\tilde\beta
c_{k+3} +\tilde\gamma c_{t+2}\,.
\ee
Because of (\ref{compare7})
this formula remains valid in the case $s=k+1$. To deduce (b) we
have to check that this quantity is nonzero after an appropriate
coordinate change
$$
(x,y)\lto (x(b),y(b))=(x,y+bx)\,.
$$
For this we perform the sequence of blowups as above with
$(x(b),y(b))$ instead of $(x,y)$ so that the centers $c_i(b)\in
C_{i-1}$ and $d_\rho(b)\in C_t$ of the blowups now depend on $b$.
According to Remark \ref{4.6}.2 these centers can be written
as
$$
c_i(b)=c_i+e_{i-1}(b) \quad\mbox{and }\quad d_\rho(b)=d_\rho+e_t(b)
$$
with polynomials $e_{i}(b)$ satisfying $e_{i}(0)=0$. Now
(\ref{coeff}) can be written as
\be\label{coeff1} \p_1\bar\eta_{s}(0,c_{s+1})(b) =\tilde\beta
c_{k+3} +\tilde\gamma c_{t+2}+\tilde\beta e_{k+2}(b) +\tilde\gamma
e_{t+1}(b) \,. \ee By Lemma \ref{4.9} and Proposition \ref{4.12},
$$
e_{k+1}(b)=0, \quad \deg e_{k+2}(b)=1,\quad \mbox{and}\quad \deg
e_{t+1}(b)=t-k\ge 3\,.
$$
In particular, $c_{k+2}(b)=c_{k+2}+e_{k+1}(b)=c_{k+2}$
does not depend on $b$,
and hence also the constants $\alpha$ and $\beta$ in (\ref{compare01})
do not depend on $b$.

We claim that the polynomial $P$ as in (\ref{polt}) does not
depend on $b$ either. Indeed, with $P(b,T)=\prod_{i=1}^r
(T-d_\rho(b))$ we have
$$
P(b,y_t(b))=\prod_{i=1}^r (y_t(b)-d_\rho(b)) =
\prod_{i=1}^r
(y_{t}-e_t(b)-d_\rho+e_t(b))
=\prod_{i=1}^r (y_{t}-d_\rho)=P(y_t)\,.
$$
In particular, $P(b,c_{t+1}(b))=P(c_{t+1})$ and
$P'(b,c_{t+1}(b))=P'(c_{t+1})$. Hence $\tilde\alpha$,
$\tilde\beta$, and $\tilde\gamma$ in (\ref{compare8}) do not
depend on $b$.

It follows that (\ref{coeff1}) is a nonzero polynomial  in
$b$ of degree $1$ if $\tilde \gamma=0$ and of degree $t-k$
otherwise.  Anyway, $\p_1\bar\eta_{s}(0,c_{s+1})(b)\neq 0$
for a suitable
choice of $b$, proving (b). \esit

\bsit\label{4.16} {\bf The cases $k+2=t$ and $k+1=t$.} If
$k+2=t$ then according to Proposition \ref{4.7}(b)
$$
\ba{l}
\bar\xi_{t}(x_{t},y_{t})=-c^{m}x_{t}^{s-t+2} +\hot\,,\\
\bar\eta_{t}(x_{t},y_{t})=\alpha x_{t}^{s-t}
+\beta x_{t}^{s-t+1}y_{t}
+\hot\,\\
\ea
$$
with $\alpha$ and $\beta$ as before. The formulas (\ref{eqxi}) for
$\xi_{t+1}$, (\ref{compare7}) for $\eta_{t+1}$, (\ref{compare9})
for $\xi_{j+1}$ and (\ref{coeff00}) for $\eta_{j+1}$ are
applicable again; note that $c_{k+3}=c_{t+1}$ since $t=k+2$.
So the proof of (a) proceeds as before. Also the proof
of (b) applies if we take into account that
$\deg\,e_{t+1}(b)=t-k=2$ in our case.

In the case $k+1=t$ we only need to establish (a). Therefore it
suffices to control the terms of lowest order of $\eta_j$ for
$j=t+1, \ldots, s$. By Proposition \ref{4.7}(a), instead of
(\ref{compare}) and (\ref{compare1}) we have to use the formulas
$$
\ba{l}
\bar\xi_{t}(x_{t},y_{t})=-x_{t}^{s-t+1} y_{t}^m \quad\mbox{and}\\
\bar\eta_{t}(x_{t},y_{t})=kx_{t}^{s-t} y_{t}^{m+1}\,. \ea
$$
Proceeding as before we obtain that $\bar\xi_{t+1}$ is a multiple
of  $x_{t+1}^{s-t+1}$, while
$$
\bar\eta_{t+1}(x_{t+1},y_{t+1})=\tilde\alpha x_{t+1}^{s-t-1} +
\hot\,,
$$
where $\tilde\alpha $ is as in (\ref{compare8}). By recursion in
the range $j=t+2,\ldots,s$, $$\bar\eta_{j}(x_{j},y_{j})
=\tilde\alpha x_{t+1}^{s-j}+\hot\, $$ (cf.  (\ref{coeff00})) and
so (a) follows. Now the proof of Lemma \ref{mainlem} is
completed.\qed \esit

\brem In case $k+1=t$ the distinguished feather $F_2$
is attached to component $C_t$. It may happen that
$\p_1\eta_{s}(0,y_s)=0$ so that there is no motion on component
$C_{s+1}$. Indeed, the formula (\ref{compare7}) holds with
$\tilde{\alpha}$ as before and $ \tilde{\beta}=0$ while $
\tilde{\gamma}$ is a bit more complicated than in (\ref{compare7})
since the term $\frac{y_t-c_{t+1}}{x_t^2}P(y_t)\bar\xi_t$ in
 (\ref{eqcompare}), which contributed before only to
higher order terms, cannot be ignored any more.
More precisely,
$$
\tilde{\gamma}=c_{t+1}^mP(c_{t+1})^{km}\left(k(km+2)
c_{t+1}P'(c_{t+1})+(km+k+1)P(c_{t+1})\right)
$$
and so $\tilde{\gamma}$ vanishes for an appropriate choice of
$c_{t+1}$. Note that $\tilde\gamma$ also vanishes after a
linear change of coordinates as in Lemma \ref{mainlem}(b), since
such a coordinate change induces no motion on component $C_{t}$.
Therefore the second highest coefficient of the expansion for
$\eta_j$, $j\le s$, vanishes  as well in all following steps. Thus
$\p_1\eta_{s}(0,c_{s+1})= \tilde{\gamma}c_{s+1}=0$ i.e., the
derivative vanishes for this choice of $c_{t+1}$. \erem

\section{Applications}

\subsection{Moving feathers}

Given a presentation
$$
X=X(M_2,c_3,\ldots,c_{t}, M_{t}, c_{t+1},\ldots, c_n,M_n)
$$
of a special smooth Gizatullin surface
$V=X\backslash D$ with data $(n,r,t)$ (see \ref{4.5}),
we consider the sequence of coordinate
systems $(x_i,y_i)$ as
in \ref{coord}. For a  fixed $i$
the curve $C_i\backslash
C_{i-1}\cong\A^1$ is the axis $x_i=0$. So it is equipped with
the coordinate $y_i$ such that $C_i\cap C_{i-1}=\{y_i=\infty\}$.
With respect to this coordinate, the data
$(M_i,c_{i+1})$ correspond to a collection of complex
numbers. We also deal with the reversed presentation
$$
X^\vee=X(M_n ,c^\vee_3,\ldots,c^\vee_{t^\vee}, M_{t^\vee},
c^\vee_{t^\vee+1},\ldots, M_2)\,,
$$
where as before $t^\vee=n-t+2$ and the curves $C_i$ and
$C^\vee_{i^\vee}$ are identified via the correspondence
fibration as in \ref{corr}. Under this identification $y_i$
yields a coordinate on $C^\vee_{i^\vee}$ so that the
data $(M^\vee_{i^\vee}=M_i, c^\vee_{i^\vee+1})$
are as well expressed by complex numbers.
The reader should
keep in mind that according to Proposition
\ref{111}  $c^\vee_{i^\vee+1}=0$.

Let $F_2$, $F_{t\rho}$ ($1\le\rho\le r$) and $F_n$
be the feathers of $X$ corresponding to the
points of $M_2$, $M_t$ and $M_n$, respectively.
The dual feathers $F^\vee_2$, $F^\vee_{t\rho}$
and $F^\vee_n$ have then mother components $C^\vee_n$,
$C^\vee_{t^\vee}$ and $C^\vee_2$, respectively.

To study the effect of an elementary shift
$$
h_{a,m}: (x,y)\to (x, y+ax^{1+m})
$$
on these presentations, we exploit as in \ref{recursion} the induced
coordinate systems $(\xi_i,\eta_i)$ depending on $a$. The data
$c_{i+1}(a)$, $M_i(a) $ and $c^\vee_{i^\vee+1}(a)$ expressed in
the coordinate $y_i$
also depend on $a$. The following result is immediate
from Proposition \ref{111} and Remark \ref{4.6}(2).

\blem\label{4-23}
If $\eta_i(0,y_i)=y_i+e_i(a)$ then
$$
c_{i+1}(a)=c_{i+1}+e_i(a), \quad M_i(a)=M_i+e_i(a), \quad
\mbox{while}\quad
c^\vee_{i^\vee+1}(a)=c^\vee_{i^\vee+1}=0\,\,\forall a\in\C\,.
$$
\elem

Consequently, if $\deg e_i(a)>0$ then
the shift
$h=h_{a,m}$ translates
$M_i=M^\vee_{i^\vee}$ and $c_{i+1}$ while keeping the point
$c^\vee_{i^\vee+1}=0$ fixed. Thus for general $a$ we have
$c_{i+1}\ne c^\vee_{i^\vee+1}$ and
$c^\vee_{i^\vee+1}\not\in M^\vee_{i^\vee}$.
Inspecting Proposition \ref{4-6} this means
that for every feather $F^\vee$ attached to a component
$C^\vee_{j^\vee}$ with $j^\vee>i^\vee$ one has $\mu^\vee> i^\vee$,
where $C^\vee_{\mu^\vee}$ is the mother
component of $C^\vee_{j^\vee}$. In other words,
feathers cannot ``cross" $C^\vee_{i^\vee}$.
Using this idea to remove jumping feathers
we can prove the following lemma.

\blem\label{4-24}
Assume that $F_2$ sits on component $C_{k+1}$
and $F_n^\vee$ on component $C^\vee_{l^\vee+1}$.
Applying suitable elementary shifts and backward shifts we can achieve
that
\bnum[(a)]
\item $k+1<\max(3,t)$ and, dually, $l^\vee +1<\max(3,t^\vee)$;
\item all feathers $F_{t\rho}$ of $X$ and $F^\vee_{t^\vee\rho}$ of
$X^\vee$ are attached to their mother components
$C_t$ and $C^\vee_{t^\vee}$, respectively.
\enum
\elem

\bproof
By Corollary \ref{4-7} we have $k+1\le l-1$ or, equivalently,
$k+1+l^\vee+1\le n+2=t+t^\vee$.
Thus $k+1\le t$ or $l^\vee +1\le t^\vee$. By symmetry
we may suppose that $k+1\le t$. We proceed in several steps.

(1) {\em After a suitable shift we can achieve that
$l^\vee +1\le t^\vee$.} Indeed, if $k+1=t$ then
by the above inequality we have $l^\vee +1\le t^\vee$.
If $k+1<t$ then using Lemma \ref{4.9}, after a shift
$c^\vee_{t^\vee +1}=0\ne c_{t+1}$.
Inspecting Lemma \ref{4-6}(b)
this shows that $l^\vee+1\le t^\vee$.

(2) {\em After suitable shifts and backward shifts (b) holds.}
If $k+1= t$ then applying  Corollary \ref{4-7} to the pairs
$(F_2, F^\vee_{t\rho})$, all feathers $F^\vee_{t\rho}$ are
attached to their mother component $C^\vee_{t^\vee}$.
If $k+1< t$ then by Lemma \ref{4.9}, after a suitable shift $h_{a,0}$
we have $c^\vee_{t^\vee+1}=0\not\in M_t$. Thus by Lemma \ref{4-6}(b) all
feathers $F^\vee_{t\rho}$ are attached to their mother component
$C^\vee_{t^\vee}$.  Applying the same arguments to the reversion,
after suitable backward shifts all feathers $F_{t\rho}$
are attached to their mother component $C_t$.

(3) {\em If $t^\vee\ne 2$ then after a suitable shift
$l^\vee +1<t^\vee$. } Indeed, using (b)
we have $c_{t+1}\not\in M_t$ so that
Proposition \ref{4.12} can be applied.
The shift $h_{a,0}$ creates a motion on component
$C_{t+1}$.
After such a motion we may assume that $c_{t+2}$ is nonzero, i.e.\
$c_{t+2}\ne c^\vee_{t^\vee}$.
By Lemma \ref{4-6}(b) this forces $l^\vee+1<t^\vee$.

Applying now (3) and its dual statement, (a) follows.
\eproof

\bcor\label{4-25}
Given a presentation $X$ of a special smooth Gizatullin surface,
by performing shifts and backward shifts we can
transform  $X$  into a $(-1)$-presentation.
\ecor

\bproof
After applying suitable shifts and backward shifts we may
assume that (a) and (b) in
Lemma \ref{4-24} are fulfilled.
Interchanging $X$ and $X^\vee$, if necessary,
we may assume that $k\le l^\vee$.
We may also suppose that $k\ge 2$,
since otherwise $X$ has $(-1)$-type and we are done.
We choose now $m$ in such a way that
$$
n\ge s=k(m+1)+2\ge n-l^\vee +1 ;
$$
this is always possible because of our assumption $k\le l^\vee $.
By Proposition \ref{motion}, the elementary shift
$h_{a,m}$ induces a motion on the curve $C_{s}$.
After this motion we will have $c_{s+1}\ne c^\vee_{s^\vee+1}=0$.
Inspecting Lemma \ref{4-6}(b), on the new surface $X^\vee$
the feather $F_n^\vee$ will sit on a component
$C^\vee_{\tilde l^\vee}$ with $\tilde l^\vee< l^\vee$.
Thus after several such shifts we can achieve that $l^\vee <k$.
Interchanging now $X$ and $X^\vee$ and continuing
as before we obtain after a finite number of steps that $k=l^\vee=1$,
as required.
\eproof

Combining this with Proposition \ref{main}
leads to the following result.

\bcor\la{cormain}
Assume that $V$ and $V'$
are smooth special Gizatullin surfaces such
that the zigzags of standard completions are equal up to reversion.
Then $V$ and $V'$ are isomorphic if and only
if the configuration invariants of $V$ and $V'$ coincide.
\ecor

\bproof
The `only if" part follows from Theorem \ref{matchmain}.
To prove the converse, let $(\bV, D)$ and $(\bV', D')$
be standard completions of $V$ and $V'$, respectively.
Reversing one of them, if necessary,
we may suppose that the dual graphs of $D$ and $D'$ are equal.
By Proposition \ref{DG+special}
both surfaces admit presentations.
Moreover by Corollary \ref{4-25} we may assume
that both presentations are of $(-1)$-type.
Applying Proposition \ref{main} the result follows.
\eproof

\subsection{Main theorem and its corollaries}

Let us now deduce the Isomorphism Theorem
\ref{0.5}(c) in the Introduction.

\bthm\label{main2}  Given a special smooth Gizatullin
$\C^*$-surface $V=\Spec\C[u][D_+,D_-]$, the isomorphism type
of $V$ is uniquely determined by the unordered pair of numbers $(\deg \{D_+\},
\deg \{D_-\})$
and the configuration of
points
\be\label{pco} \supp(\lfloor
-D_+-D_-\rfloor)=\{p_1,\ldots,p_r\}\,,\ee
up to the natural
action of the automorphism group $\Aut (\A^1)$ on such
configurations.
\ethm

\bproof
This follows immediately from Corollary \ref{cormain}.
Indeed, according to
Corollary \ref{specialconfig} the configuration invariant
is given by $\{p_1,\ldots,p_r\}\in \fM_r^+$, whereas the boundary zigzag is
up to reversion uniquely  determined by the numbers
$\deg \{D_+\}$ and $\deg \{D_-\}$.
\eproof

The next result and Corollary \ref{5.6.6} below yield the
first assertion of Theorem \ref{0.3} in the Introduction.

\bthm\label{0.33}
Every special smooth Gizatullin
surface carries a $\C^*$-action. Moreover, if
$V$ is of type I then the conjugacy classes
of $\C^*$-actions on $V$ form in a natural way a one-parameter
family, while in case of type II they form a two-parameter family.
\ethm

\bproof
By Corollaries \ref{DG+special} and \ref{4-25}, $V$ admits a presentation
of $(-1)$-type $X=X(M_2, c_3, \ldots, M_t, \ldots M_n)$  as in \ref{4.5},
where $M_t=\{p_1, \ldots, p_r\}$. After reversing the presentation,
if necessary, we may suppose that $t<n$.

First assume that also $t>2$. Let us consider the  $\C^*$-surface
$V'=\Spec \C[t][D_+,D_-]$ with
$$
D_+=-\frac{1}{t-1} [p_+] \quad\mbox{and}\quad
D_-=-\frac{1}{n-t+1} [p_-]-\sum_{\rho=1}^r [p_\rho]\,,
$$
where $p_+, p_-$ are different and not contained in $M_t$.
Applying Corollary \ref{specialconfig}(a)
the configuration invariant of this surface
is given by $M_t\in \fM_r^+$.
 Moreover by (\ref{ezigzag})
the boundary zigzags of $V$ and $V'$ coincide. Now by
virtue of Corollary
\ref{cormain}, $V$ and $V'$ are isomorphic.

Furthermore, the equivariant isomorphism type of $V'$ depends on
the position of $p_+$, $p_-$ and $M_t$ while the abstract
isomorphism type only depends on the configuration $M_t$.
Thus in the case $r\ge 2$ the
family of conjugacy classes of $\C^*$-actions on
the surface $V$ depends on
two parameters,
while in the case $r=1$ it depends on just one parameter.
Since for $r\ge 2$ the surface is of type II while
for $r=1$ it is of type I, the result follows in this case.

Finally if $t=2$ and  $r\ge 2$
then the above reasoning together with Corollary
\ref{2.3} shows again that the  $\C^*$-actions on $V$ form a
one-parameter family.
\eproof

\bsit\label{dgs}
Any Danilov-Gizatullin surface $V$ admits a presentation
\be\label{DGcond}
X=X(M_2, c_3, \ldots, c_n,M_n)\mbox{ with }
n\ge 2,\, |M_2|=|M_n|=1\mbox{ and } M_i=\emptyset \mbox{
otherwise} .\ee
Indeed, let $V=\Sigma_k\backslash S$,
where $S$ is an ample section in the Hirzebruch surface
$\Sigma_k\to \PP^1$ with $S^2=n>k$. Blowing up a point of $S$
successively yields a semistandard boundary zigzag
$[[0, -1, (-2)_{n-1}]]$ and so a standard zigzag $[[0, 0,
(-2)_{n-1}]]$. Hence by Lemma \ref{2.7}, $V$ admits a
presentation as in (\ref{DGcond}). Applying
Corollary \ref{cormain} we recover the
theorem of Danilov and Gizatullin \cite{DaGi}
cited in the Introduction\footnote{See also \cite{CNR,
FKZ4}.}.
\esit

\bcor\label{DG-uniproof} The isomorphism type of a
Danilov-Gizatullin surface $V=\Sigma_k\backslash S$ depends
only on $n=S^2$. \ecor

The following corollary completes the proof of Theorem
\ref{0.3} in the Introduction.

\bcor\label{5.6.6} A special smooth Gizatullin surface $V$ of type
I admits a one-parameter family of pairwise non-equivalent
$\A^1$-fibrations $V\to\A^1$, while for  type II it admits
such a family depending on two parameters.
\ecor

\bproof
Every $\C^*$-action on $V$ extends to two
mutually reversed equivariant standard completions $(\bV, D)$
and $(\bV^\vee, D^\vee)$. On each of them there is an
associated $\A^1$-fibration of $V$ induced by the linear system
$|C_0|$. By Lemma 5.12 in \cite{FKZ3},  if these
$\A^1$-fibrations are conjugated then  the associated
extended divisors are isomorphic as reduced curves.
Moreover, by Proposition 5.12 in \cite{FKZ2}
 and its proof, the latter holds if and only
if the associated $\C^*$-actions are conjugated. Consequently,
there are at least as many conjugacy classes of
$\A^1$-fibrations as of $\C^*$-actions. \eproof

Every $\A^1$-fibration on $V$ arising as in Corollary
\ref{5.6.6} is compatible with a certain $\C^*$-action. In the
next subsection we exhibit further $\A^1$-fibrations that do not
appear in this way.

\subsection{Applications to $\C_+$-actions and
$\A^1$-fibrations }

 Every $\A^1$-fibration $V\to\A^1$ on
a Gizatullin surface $V$
arises as the orbit map of a $\C_+$-action $\Lambda$ on $V$,
see e.g.\ Lemma 1.6 in \cite{FlZa2}.
Of course, $\Lambda$ is not unique. If,
say, $\partial$ is the locally nilpotent derivation associated to
$\Lambda$ and $a \in \Ker \partial$ is non-zero then
$\partial'=a\partial$ is also locally nilpotent and generates a
$\C_+$-action with the same general orbits.
It is well known that this  is
the only ambiguity in associating
a $\C_+$-action
to a given $\A^1$-fibration, see e.g., \cite{KML}. Thus
to classify $\C_+$-actions up to conjugation
is essentially equivalent
to determining all $\A^1$-fibrations over $\A^1$ up to conjugation.

According to Theorem 5.2 in \cite{FKZ3}
for a wide class of Gizatullin surfaces
there are at most 2 conjugacy class of $\A^1$-fibrations;
see also Section 6.4 below.
However, for surfaces arising from presentations
the assumptions of this theorem
are almost never satisfied.
Thus we concentrate below on surfaces
of the latter class and even of a subclass called
quasi-special surfaces.
In Theorem \ref{fibration.140} we give a
complete classification of
$\A^1$-fibrations on such surfaces
up to conjugation.
Our methods can be applied
as well to other classes of Gizatullin surfaces.

To start with,
let us recall some notation.
Let $V$ be a smooth Gizatullin surface and
$(\bV,D)$ be a standard completion of $V$.
The linear system $|C_0|$ defines
a $\PP^1$-fibration $\Phi_0 :
\bV \to \PP^1$ called the {\em standard fibration}. Its
 restriction $\varphi_0=\Phi_0|V:V \to \A^1$ will be called
the {\em associated $\A^1$-fibration} on $V$. A key observation is
the following result, which summarizes
5.11, 5.12 in \cite{FKZ3} and their proofs.

\bprop\la{fibration.10} Let $V$ be a normal Gizatullin surface
equipped with an $\A^1$-fibration $\pi:V\to \A^1$.
Then $V$ admits
a standard completion $(\bV,D)$ such that
$\pi$ is induced by
the standard fibration $\Phi_0:\bV\to\PP^1$.
Furthermore, if\/ $V\not\cong \A^1\times\A^1_*$
and $\Phi_0':\bV'\to\PP^1$ is a second extension of $\pi$ to another
standard completion $(\bV',D')$ of $V$ then  there
exists an isomorphism $\psi : \bV \setminus
C_0 \to \bV' \setminus C_0'$ with $\psi|V=\id_V$ and $\Phi_0'\circ \psi=\Phi_0$.
In particular, the extended graphs of both completions coincide.
\eprop

Thus every $\A^1$-fibration $V\to\A^1$ arises as a standard fibration
from a standard completion $(\bV,D)$,
which is unique up to a modification at $C_0$.

The pair of linear systems
$|C_0|$, $|C_1|$ defines a
birational morphism $\Phi=(\Phi_0,\Phi_1) :
\bV \to Q=\PP^1 \times \PP^1$ to
the quadric $Q$ called the standard morphism in \ref{1.2}.
In suitable coordinates $(x,y)$ on $Q\setminus(C_0\cup C_1)\cong\A^2$ we have
$\Phi_0=x\circ\Phi$ on $Q\setminus(C_0\cup C_1)$. Moreover
$$
C_0\cong \{\infty\}\times \PP^1,\quad
C_1 \cong \PP^1 \times\{\infty\},\quad\mbox{and}\quad
C_2 \cong \{0\}\times\PP^1\,,
$$
while the curves $C_3,\ldots , C_n $ are contained in the preimage of the origin.

\bsit\la{fibration.30} Let as in \ref{2.9}
$\Aut_y (\A^2)$ stand for
the group of all automorphisms
of $\A^2$ stabilizing the $y$-axes, and let
$\Aut_{y,0} (\A^2)\subseteq \Aut_y (\A^2)$
denote the stabilizer of the origin.
Every automorphism $\alpha \in \Aut_{y,0}$ can be written as
$$\alpha . (x, y)=( \lambda_1 x, \lambda_2 y+q(x))\quad \mbox{with} \quad
\lambda_1,\lambda_2\in \C^* \,\,\,\mbox{and}\,\,\, q\in\C[x],\,\,q(0)=0\,,$$
cf.\ (\ref{auto}).
Hence $\Aut_{y,0}(\A^2)=H\rtimes \T$ is a semidirect product
of the torus $\T=\C^{*2}$
acting on $\A^2$ by
\be\label{fibration.40}\lambda.(x,y)=(\lambda_1 x, \lambda_2 y), \quad
\lambda=(\lambda_1, \lambda_2)  \in \T\,,
\ee
and the abelian group $H$ of all triangular automorphisms
\be\label{fibration.50}
h_q: (x,y)
\mapsto (x,y+q(x)),\qquad q\in\C[x],\,\,q(0)=0\,.\ee
Here $\T$ acts on $H$ by conjugation
$(\lambda, \mu).h_{q(x)} =
h_{\mu^{-1}q(\lambda x)}\,.$
\esit

\bsit\la{1350} Given a standard
completion $(\bV,D)$ of a Gizatullin surface $V$,
we consider the restriction $\Psi=\Phi|_{\bV
\setminus (C_0 \cup C_1)}: \bV \setminus (C_0 \cup C_1) \to
\A^2$ of the standard morphism.
Clearly $\Psi$ is the contraction of the feathers
and the components $C_3, \ldots ,C_n$ of the zigzag.
Every automorphism $\alpha \in \Aut_{y,0}(\A^2)$ extends to
$\bV \setminus (C_0 \cup C_1)$
inducing a commutative diagram
\be
\bdi\label{fibration.15}
\bV \setminus (C_0 \cup C_1) & \rTo^{\tilde
\alpha}& \bV' \setminus (C_0' \cup C_1') \\
\dTo>\Psi && \dTo{\Psi'} \\
\A^2 & \rTo^{\alpha } & \A^2,
\edi
\ee
where $\bV'$ is a completion of another Gizatullin
surface $V'$ isomorphic to $V$, and $\Psi'$
is defined similarly as $\Psi$. In fact $\tilde
\alpha$ can be extended to an automorphism $\psi : \bV \setminus
C_0 \to \bV' \setminus C_0'$. Note that $\alpha \in
\Aut_{y,0}(\A^2) $ in (\ref{fibration.15})
is compatible with the standard fibration $\Phi_0$
on $V$ as it
preserves the $x$-coordinate up to a multiple. Proposition
\ref{fibration.10} implies the following result.
\esit

\bcor\label{fibration.20}
Every isomorphism $\psi : \bV \setminus C_0
\to \bV' \setminus C_0'$ as in Proposition \ref{fibration.10} is
induced by an automorphism $\alpha\in\Aut_{y,0}(\A^2)$.
\ecor

In terms of presentations this
leads to the following proposition.

\bprop\label{fibration.90} Let $V$ be a smooth Gizatullin
surface, and let
$$
X_n=X(M_2, c_3,\ldots, M_i, \ldots ,c_n, M_n)
\,\,\,\,\mbox{and} \,\,\,\,
X'_n=X(M'_2, c'_3,\ldots, M'_i, \ldots ,c'_n, M'_n)
$$
be two presentations
of $V$. Then the associated $\A^1$-fibrations
$\Phi_0|V,\Phi'_0|V:V\to\A^1$ are conjugated
if and only if there is an automorphism
$\alpha \in \Aut_{y,0}(\A^2)$ with
$$
X'_n=\alpha_*(X_n).
$$
\eprop

The Motion Lemma \ref{motion}
was stated for presentations of
special Gizatullin surfaces. However, part (a)
remains true more generally for quasi-special
surfaces which we introduce below.

\bdefi\la{2000.000} A presentation
$X=X(M_2,c_3,\ldots, M_n)$ as in \ref{basic0} will be called
{\em quasi-special of type} $(n,k)$ (or simply quasi-special) if
$0\in M_2$ and $c_3=\ldots =c_{k+1}=0$ but $c_{k+2}\ne 0$ and
$c_{i+1}\not\in M_i$ for $i=k+1, \ldots, n-1$.
As usual, here the data $(M_i,c_{i+1})$ are considered
as collections of complex numbers expressed in
the coordinates introduced in Section 5.1.
It will be convenient to complete these
data by introducing also the point $c_{n+1}$
as the center of mass of $M_n$.

Thus all components $C_{i+1}, \, i \geq 2$ of the zigzag $D$
are of $+$-type i.e., $C_{i+1}$ is created
by blowing up a point $c_{i+1}\in
C_i\setminus C_{i-1}$.
Furthermore, for $k\ge 2$ the
extended divisor $D_\ext$ of $(\bV,D)$
has exactly one feather, denoted $F_2$,
with $F_2^2\le-2$; the mother
component of this feather is $C_2$
and the neighbor in the zigzag is
$C_{k+1}$. Otherwise the presentation can be arbitrary.

In particular,
every $(-1)$-presentation is quasi-special,
cf.\ Definition \ref{basic} and Corollary \ref{DG+special}.
By Corollary
\ref{4-25} every special surface admits
quasi-special presentations. However
there exist also presentations of special surfaces
which are not quasi-special.
\edefi

\blem\label{motion2}
Given a quasi-special presentation $X=X(M_2,c_3,\ldots, M_n)$,
the first motion under the elementary shift
$h_{a,m}$ as in (\ref{mshift})
occurs on component $C_s$, where as in Section 5 $s=k(m+1)+2$,
i.e., $e_i=0$ for all $i<s$,
while $\deg e_s>0$.
\elem

\bproof
The map $a\mapsto h_{a,m}$
yields a $\C_+$-action on $\A^2$.
We can lift its infinitesimal generator $\partial$
to the surface
$X$ as a meromorphic vector field. The curves on which
$h_{a,t}$ is constant are characterized by the fact that
$\partial$ is regular and identically zero on them.
Moreover the curve of first motion $C_l$ is characterized
by the property that $\partial$ (tangent to $C_l$)
is regular and nonzero in
the general points of $C_l$. However this property is
independent of blowdowns of $(-1)$-feathers.
Thus it is enough to find the curve of first motion
in the case where $X$ has no $(-1)$-feathers,
which is just the
Danilov-Gizatullin case. Applying Proposition
\ref{motion}(a), the result follows.
\eproof

We let below $m_0(n,k)=\lfloor {\frac{n-2}{k}}\rfloor\,.$

\bcor\label{fibration.100}
With $X=X_n$ and  $h=h_{a,m}$ as in Lemma \ref{motion2},
if $m\geq m_0(n,k)$ then $X_n=h_*(X_n)$ i.e.,
$h$ generates an automorphism of $X_n \setminus C_0$.
\ecor

\bproof By the
Motion Lemma \ref{motion2} the
data $M_i,c_{i+1}$ remain unchanged for $i \leq k (m+1) +2$,
in particular for $i\leq n$ so that $X_n=h_*(X_n)$.
\eproof

Given $m\in\N$, we let
$H_{m}$ denote the subgroup of $H$
generated by all elements $h_q$ as in
(\ref{fibration.50}) with $\deg q(x)\le m$.

\bdefi\la{3000.000} A presentation
$X=X_n$ of type $(n,k)$
will be called {\em semi-canonical} if it is quasi-special
with $c_{k+2}=1$ and $c_{j+1}=0$ for all
$j=k+2,\, 2k+2,\,
\ldots , m_0k +2$, where $m_0=m_0(n,k)$. \footnote{In particular,
for $m_0={\frac{n-2}{k}}$ the center of mass $c_{n+1}$ of $M_n$
should be also $0$.}
\edefi

\blem\label{fibration.110} Every quasi-special presentation
$X=X_n$ of type $(n,k)$ can be transformed into a semi-canonical
one by applying a suitable automorphism
$h_0\circ \lambda\in\Aut_{y,0}(\A^2)$, where $\lambda\in\T$ and
$h_0 \in H_{m_0}$.
Furthermore, such an element $h_0\in H_{m_0}$ is unique.
\elem

\bproof The torus $\T$ acts non-trivially on
$C_{k+1}\setminus C_{k} \cong\A^1$ with the fixed point $0$.
Hence a suitable $\lambda\in\T$ sends $c_{k+2}\neq 0$ to the
point $1$.

Consider further an elementary shift $h_{q}$ as in
(\ref{fibration.50}) with $q(x)=ax^{m+1}$, where $m\le m_0$.
By the generalized Motion Lemma \ref{motion2} it
does not change the data $(M_i,c_{i+1})$ for $i \leq s-1$, where
$s=k(m+1)+2$, and induces a nontrivial translation on
$C_s \setminus C_{s-1}\cong \A^1$.  So $h_{q}$ sends $c_{s+1}$
to $0$ for a suitable value of
$a\in\C_+$.
Applying such actions repeatedly for $m=0,1, \ldots , m_0$ we
obtain the desired semi-canonical presentation.
The uniqueness part is easy and left to the reader.
\eproof

\bsit\la{2000.020} Given an arbitrary presentation $X_n=X(M_2,c_3,\ldots,M_n)$
we can adopt the construction
of canonical coordinate systems on $X_n$ in \ref{coord}
as follows.
Starting with the affine coordinates $(x_2,y_2)=(x,y)$ on the quadric $Q$,
at the first step we blow up $M_2\cup\{c_3\}$ getting a coordinate chart
$(x_3,y_3)$, at the second one we blow up $M_3\cup\{c_4\}$
getting $(x_4,y_4)$, and so forth.
\footnote{In the case of special surfaces the enumeration is
slightly different from \ref{coord} since we perform here
the first two blowups in one step.}

Under this procedure the coordinates
are given recursively by the formulas in \ref{coord}(4).
So  letting $M_i=\{d_{i1},\ldots d_{is_i}\}$ for $i=2,\ldots, n$,
by virtue of (\ref{polt})  we have
\be\label{polt111}
(x_{i+1},y_{i+1})=\left(\frac{x_i}{P},
\frac{y_i-c_{i+1}}{x_i}P\right)\,, \quad \mbox{where}\quad P
= \prod_{j=1}^{s_i} (y_i-d_{ij})\in\C[y_i]\,. \ee
\esit

Let us investigate the torus action in terms of these coordinates.
In the following, for $m=(m_1,m_2)\in \Z^2$ and
$\lambda =(\lambda_1, \lambda_2)\in\T$ the power  $\lambda ^m$ stands for
$\lambda_1^{m_1} \lambda_2^{m_2}$.

\blem\label{fibration.125} \bnum[(a)]
\item  $x_i$ and $y_i$ ($i=2,\ldots, n+1$) are quasi-invariant functions under
the induced $\T$-action, i.e.\ there are vectors $a_i$, $b_i\in \Z^2$ with
\be\la{form6}
\lambda.x_i=\lambda^{a_i}x_i\quad\mbox{and}\quad
\lambda.y_i=\lambda^{b_i}y_i
\quad\mbox{for}\quad\lambda\in \T
\ee
In particular, this action leaves the point on $C_i$ with
$y_i$-coordinate zero unchanged.
\item If we form the matrix $A_i=(a_i,b_i)$ with
column vectors $a_i,b_i$ then $\det A_i=1$ for $2\le i\le n+1$.

\item $\T$ acts on the set of all quasi-special presentations.
\enum
\elem

\bproof
(a) is certainly true for $i=2$ with $a_2$, $b_2$ being
the standard basis of $\Z^2$.
Assume that for $i\ge 2$ this action can be expressed in
$(x_i,y_i)$-coordinates as in (\ref{form6}).
In particular, with the notations as in \ref{2000.020},  $\lambda.d_{ij}= \lambda^{b_i}d_{ij}$ and
$\lambda.c_{i+1}= \lambda^{b_i}c_{i+1}$. Hence $\lambda.P= \lambda^{s_ib_i}P$ and so
$$
\lambda.x_{i+1}=\frac{\lambda.x_i}{\lambda.P}=
\lambda^{a_i-s_ib_i}\frac{x_i}{P}= \lambda^{a_i-s_ib_i}x_{i+1}\,
$$
and similarly
$$
\lambda.y_{i+1}=\frac{\lambda.(y_i-c_{i+1})}{\lambda.x_i}\lambda.P=
\lambda^{b_i-a_i+s_ib_i}\frac{y_i-c_{i+1}}P= \lambda^{b_i-a_i+s_ib_i}y_{i+1}\,.
$$
Thus
\be\la{form66}
\lambda.(x_{i+1},y_{i+1})=(\lambda ^{a_{i+1}}x_{i+1},
\lambda ^{b_{i+1}}y_{i+1})
\mbox{ with } a_{i+1}=a_i-s_ib_i\,,
b_{i+1}=b_i(1+s_i)-a_i\,,
\ee
proving (a). To deduce (b) we note that
the recursion (\ref{form66}) can be expressed as
$$
A_{i+1}=B_iA_i=B_i\cdot B_{i-1}\cdots \cdot B_2,
\quad\mbox{where}\quad
B_i:=\left(\begin{array}{cc}1 & -s_i \\-1 & 1+s_i\end{array}\right)
\quad 2\le i\le n\,.
$$
Since $A_2$ is the identity matrix and $\det B_i=1$ we obtain by induction that also $\det A_i=1$ for $i= 2,\ldots, n+1$.
Finally (c) was already observed in \ref{2.9}.
\eproof

\bsit\la{4000.000}
Given a semi-canonical presentation
of type $(n,k)$ we
let $\T'\subseteq\T$ denote the stabilizer subgroup of the point
$c_{k+2}=1$. According to (\ref{form6})
$$
\T'=\{\lambda\in \T\, | \,
\lambda^{b_{k+1}}=1\}\,.
$$
Hence $\T'\cong\C^*$; indeed,
by Lemma \ref{fibration.125}(c)  the
columns of the matrix $A_i$ with $\det A_i=1$
are unimodular i.e., primitive lattice vectors.
Due to Lemma \ref{fibration.125}(c)
the 1-torus $\T'$ is independent
of the choice of a semi-canonical presentation
as soon as the sequence $(s_i=|M_i|)_{i=2, \ldots, n}$
is fixed. Furthermore, it acts
on the set of all such presentations.

It also acts non-trivially
on every component $C_i$ of the zigzag except for
the component  $C_{k+1}$.
To show this, let us first grow in the blowup process
the zigzag $D$ together with the feather $F_2$.
At this point
our $\T'=\C^*$-action lifts to the resulting surface, say $X'$
with $C_t$ as the unique parabolic component.
Since our surface $X=X(M_2,c_3,\ldots, M_n)$ is
obtained from $X'$ by blowing up the remaining $(-1)$-feathers,
$\T'$ will act  also non-trivially on $C_i$,
$i\ne k+2$, when considered as curves in  $X$.

Combining now Lemmas \ref{fibration.110},
\ref{fibration.125} and Proposition \ref{fibration.90}
leads to the following corollary.
\esit

\bcor\la{fibration.cor}
Given two semi-canonical presentations $X_n$ and $X_n'$
of the same Gizatullin surface
$V$, the associated
$\A^1$-fibrations
$\varphi_0,\,\varphi'_0:V\to\A^1$
are conjugated
if and only if $X'_n=\lambda_*(X_n)$
for some $\lambda\in \T'$.
\ecor

We can normalize the presentation further
using the action of $\T'$.

\bdefi\label{fibration.130}
A semi-canonical presentation $X_n$
(see Definition \ref{3000.000})
will be called $c_{l+1}$-{\em canonical},  if
$c_{i+1} =0$ for $k+2\le i \leq l-1$ and
$c_{l+1}=1$. Such a $c_{l+1}$-canonical presentation can only exist
for $k+2\le l\le n$
and  $l \ne k+2,\, 2k+2, \ldots ,\, m_0k+2$.

It can happen that $c_{i+1}=0$ for every $i\neq k+1$ so that
the presentation
is not $c_{l+1}$-canonical whatever $l$ is.
In this case we let $l$ ($2\le l\le n$)
be the minimal index with $l\ne k+1$ such that
$M_l^*:=M_l\backslash \{0\} \ne \emptyset$.
\footnote{We note
that $M_l^*=M_l$ if $3\le l\le n-1$ and $l\ne k+1$.}
The presentation
$X_n$ will be called $M_l$-{\em canonical}\/ if
$$
m_l:=\prod_{m\in M_l^*}m=1\,.
$$
In the remaining case, where
$$
c_{i+1}=0\mbox{ for all } i=k+2, \ldots ,
n\quad\mbox{and}\quad
M_l^*=\emptyset \quad\mbox{for all } l
\ne k+1\mbox{ with }2\le l\le n,
$$
$X_n$ will be called $*$-{\em canonical}.

A presentation will be called for short {\em canonical}\/ if it is $a$-canonical
for some $a\in \{c_{l+1},M_l,*\}$.
\edefi

The next lemma is immediate from the fact that
the 1-torus $\T'$ acts in a nontrivial way
on each component $C_i$, $i\ne k+1$,
with the only fixed points $0$ and $\infty$.

\blem\la{5000.000} Every semi-canonical presentation
can be transformed into
a canonical one by an element $\lambda\in \T'$.
\elem

We also need below the following simple lemma.

\blem\la{5000.010} Let $X_n$ be a
canonical presentation  of type $(n,k)$. If it
is $c_{l+1}$-canonical or $M_l$-canonical for
some  $l\le n$, then the subgroup
$$
G_{kl}=
\{\lambda\in\T\mid \lambda^{b_{k+1}}=
\lambda^{b_l}=1\}\subseteq\T'\,
$$
is a finite cyclic group of order
$|\det\, (b_{k+1}, b_{l})|\neq 0$.
\elem

\bproof
Since $G_{kl}$ is contained in $\T'\cong\C^*$ it is either
cyclic or equal to $\T'$. As observed in \ref{4000.000}, $\T'$ acts
non-trivially on $C_l$ for $l\ne k+1$. As by (\ref{form6}) $G_{kl}$ is the
subgroup of all elements of $\T'$ acting trivially on $C_l$, it
is finite.
\eproof

We now come to our main classification results.
By Proposition \ref{fibration.10}, if
the standard $\A^1$-fibrations associated to
two different presentations of the
same Gizatullin surface $V$ are conjugated,
then the corresponding extended divisors
are isomorphic. Therefore, if these presentations are
quasi-special then they  are of  the same type
$(n,k)$.

\bthm \la{fibration.140} Let $X_n$ be an $a$-canonical and $X_n'$
an $a'$-canonical presentations of the same type $(n,k)$ of
a smooth Gizatullin surface $V$, where
$a\in \{c_{l+1},M_l,*\}$ and  $a'\in \{c_{l'+1},M_{l'},*\}$
are as in Definition \ref{fibration.130}.
Let $\varphi_0,\,\varphi_0':V\to
\A^1$ be the associated $\A^1$-fibrations.

\bnum[(a)]
\item If $a\ne a'$ then  $\varphi_0$ and $\varphi_0'$
are not conjugated.

\item If $a=a'=c_{l+1}$ or $a=a'=M_l$ then $\varphi_0$ and $\varphi_0'$ are conjugated
if and only if  $X'_n=\lambda_*(X_n)$ for some
$\lambda\in G_{kl}$.

\item If $X_n$ and $X_n'$
are both $*$-canonical,
then  $\varphi_0$ and $\varphi_0'$ are conjugated
if and only if $X_n=X_n'$.
Furthermore $\T'$
leaves $X_n$ invariant and yields a $\C^*$-action on $V$.
Conversely, if a canonical presentation $X_n$
admits an effective $\C^*$-action
then it is $*$-canonical.
\enum
\ethm

\bproof Assume that $\varphi_0$ and $\varphi_0'$ are conjugate.
Then by Corollary \ref{fibration.cor}
there exists $\lambda \in
\T'$ transforming $X_n$ into $X_n'$.
To show (a) let us first suppose that $X_n$
is $c_{l+1}$-canonical and $X_n'$ is $c_{l'+1}$-canonical with $l< l'$.
Then $c_{l+1}=1$, which implies by (\ref{form6}) that
$c'_{l+1}=\lambda^{b_{l+1}}\neq 0$. The latter
contradicts the assumption
that $X_n'$ is $c_{l'+1}$-canonical,
see Definition \ref{fibration.130}.
A similar argument yields the other cases in (a).

In (b), assuming again that $\varphi_0$ and $\varphi_0'$
are conjugated,
in the case of $c_{l+1}$-canonical presentations we obtain
as before that $X'_n=\lambda_*(X_n)$ for some $\lambda \in
\T'$. However, since $c_{l+1}=c'_{l+1}=1$,  by (\ref{form6})
we have
$\lambda^{b_{l+1}}=1$. Hence $\lambda\in G_{kl}$.
The proof in the case of $M_l$-canonical presentations
is similar. The converse can be easily deduced
along the same lines.

Finally, if both presentations are $*$-canonical
then $\T'$
leaves all points of $M_i$ and $c_{i+1}$ fixed. Hence it
induces a $\C^*$-action on $X_n$ and on $V$. Conversely, if
$X_n$ is canonical and
admits an effective $\C^*$-action then this action
(which has just one parabolic
component $C_{k+1}$) must stabilize
$M_i$ and $c_{i+1}$ for all $i\neq k+1$.
This can occur only if $X_n$ is $*$-canonical.
\eproof


\bcor\la{7000.000} Let $\cX(V)$ denote the set
of all canonical presentations
of a smooth Gizatullin surface $V$
and $\cY(V)$ denote the set of all conjugacy
classes of $\A^1$-fibrations
$V\to\A^1$. Then
the natural correspondence $\cX(V)\to\cY(V)$,
which sends
a canonical presentation into
the conjugacy class of its standard
fibration, is surjective and finite-to-one.
\ecor

In the case of Danilov-Gizatullin surfaces this
leads to the following complete classification of $\A^1$-fibrations;
cf.\ \ref{1000.000} in the Introduction.

\bcor\la{7000.010}
The $\A^1$-fibrations on the Danilov-Gizatullin
surface $V(n)$ are completely classified by the following
canonical presentations.
\bnum[(a)]\item
For every type $(n,k)$, $k=1,\ldots,n-1$, there is exactly one
$*$-canonical presentations
of $V(n)$.

\item There are no $M_l$-canonical presentations on $V(n)$.

\item
For every type $(n,k)$ with $1\le k\le n-1$ the surface $V(n)$ has a
$c_{l+1}$-canonical presentation
if and only if
\be\la{form77}
k+2\le l \le n\quad\mbox{and}\quad
l\ne k+2,2k+2,\ldots, m_0k+2\,,
\ee
where as before $m_0=\lfloor \frac{n-2}{k}\rfloor$.
Furthermore, if $l\le n$ and  $ak+2< l< (a+1)k+2$ with $1\le a\le  m_0$,
then these $c_{l+1}$-canonical presentations
form a family of dimension $r(n)=(n-l)-(m_0-a)$.
\enum
\ecor

\bproof
Every surface admitting a presentation
with $s_2=s_n=1$ and $s_i=0$  for $i\neq 2,n$, is a
Danilov-Gizatullin surface (see \ref{dgs}).
Clearly one can choose such a presentation of any given type $(n,k)$.
In view of Theorem \ref{fibration.140} this proves (a).

(b) and also the first part of (c) are immediate from
Definition \ref{fibration.130}. The second part of
(c) is a consequence of the fact,
that for a $c_{l+1}$-canonical presentation
the positions of the remaining points
$c_{i+1}$, $l<i\le n$ and $i\ne 2k+2,\ldots,m_0k+2$,
can be freely chosen and so give rise to
a family of the claimed dimension.
\eproof

\bexa\label{fibration.150}
In particular, choosing $k=2$ and $l=5$ in Corollary \ref{7000.010}(c)
any Danilov-Gizatullin surface $V(n)$ of index $n\ge 7$ carries
continuous families of pairwise non-conjugated
$\A^1$-fibrations with the number
of parameters being an increasing function of $n$.

If $n<6$ then according to Corollary \ref{7000.010}
all canonical presentations of
$V(n)$ are $*$-canonical
and so related to a $\C^*$-action. If $n=6$ then besides
the five $*$-canonical presentations there are only
two further canonical presentations, one for $k=2$, $l=5$ and
another one for $k=3$, $l=6$. Thus for any $n\ge 6$ there exist
$\A^1$-fibrations on $V(n)$ not related to $\C^*$-actions.
\eexa

With a similar reasoning we obtain the following result
for special surfaces.

\bcor
Let $V$ be a special surface of type I or II.
Then there are families of pairwise non-conjugated
$\A^1$-fibrations $V\to\A^1$ depending on
$r(n)\ge 1$ parameters with $\lim_{n\to \infty}r(n)=\infty$.
\ecor

\bproof
By Corollary \ref{5.6.6} $V$ admits at least
a one-parameter family of
pairwise non conjugated $\A^1$-fibrations.
To construct other such families
we proceed similarly as in Corollary \ref{7000.010}
for Danilov-Gizatullin surfaces.
We assume first that $V$ is special of type II so that
there is a presentation
$X=X(M_2, c_3,\ldots, M_t, \ldots  c_n, M_n)$ with
$|M_2|=|M_n|=1$,  $s_t=|M_t|\ge 2$ and
$M_i=\emptyset$ otherwise, where $2<t<n$.
According to Corollary \ref{cormain} the isomorphism class
of the affine part $V$ is uniquely determined
by its configuration invariant $[M_t]\in \fM^+_{s_t}$.

This implies in particular that for any $k=1,\ldots, n-1$
there are
quasi-special presentations of type $(n,k)$ of $V$.
Furthermore, different choices of the points $c_{i+1}$
lead to the same surface as long as we keep
the configuration $[M_t]\in \fM^+_{s_t}$ fixed.
Thus one can construct
families of $\A^1$-fibrations depending on the
same number $r(n)$ of parameters
as in Corollary \ref{7000.010}(c).

In the case of special surfaces of type I
the reasoning is similar; we leave the details to the reader.
\eproof

Let us give another application to
Danielewski-Fieseler surfaces; see \cite{Du}.
Recall that such a surface $V$ is
an affine surface
equipped with an $\A^1$-fibration $\pi:V\to \A^1$ such
that all scheme-theoretic fibers
$\pi^{-1}(a)$, $a\ne 0$, are smooth affine lines,
while $\pi^{-1}(0)$ is a disjoint union
of smooth affine lines. It is necessarily smooth.
We give below an alternative proof of the following result of A.\ Dubouloz,
see Corollary 4.13 in \cite{Du2}.

\bprop\la{dub}
For a Danielewski-Fieseler surface $V$ the following are equivalent.
\bnum[(a)]
\item Its $\A^1$-fibration is uniquely determined up to conjugation;
\item either $V$ is not Gizatullin or $V$ is isomorphic to a hypersurface
in $\A^3$ given by an equation $\{xy=P(z)\}$, where $P\in \C[z]$
has only simple roots.
\enum
\eprop

\bproof
By the result of Gizatullin \cite{Gi2}
a non-Gizatullin surface $V$ carries at most  one
$\A^1$-fibration $V\to\A^1$ up to an automorphism of the base.
Moreover, by Proposition \ref{stren} below
(see also \cite{Dai,ML})
the surface $V=\{xy=P(z)\}$ in $\A^3$
has only one $\A^1$-fibration $V\to\A^1$
up to conjugation.
This proves (b)$\To$(a).

To show the converse we may assume that $V$ is
Gizatullin and not the surface $\{xy=P(z)\}$ in $\A^3$.
By Proposition \ref{fibration.10} there is a standard completion
$(\bV,D)$ of $V$ such that $\pi$ extends to the standard
fibration on $\bV$.

The surface $\{xy=P(z)\}$ in $\A^3$ has DPD presentation with
$D_+=0$ and $D_-=-\div (P)$, see e.g., \cite{FlZa1}.
Hence inspecting Proposition \ref{eboundary} the boundary zigzag $D$
is of length $n=2$, and every surface with a boundary zigzag
of length $2$ arises in this way. Hence under our assumptions
$n>2$.

The feathers of the extended divisor
$D_\ext$ are just the components of the fiber
$\pi^{-1}(0)\subseteq V$. Since they are all reduced in the fiber, to
create the surface $\bV$ from the quadric only outer
blowups can occur. Thus $\bV$ arises form a
presentation of $(-1)$-type say, $X_n=X(M_2, c_3, \ldots, M_n)$.
As before we suppose that the data $M_i$ and $c_i$
are represented by complex numbers in the standard
coordinates of Section 5.1. Applying an elementary
shift we may as well assume that $0\in M_n$.
Inspecting Proposition \ref{4-6}, the reversion
$X_n^\vee$ is then not any longer of $(-1)$-type.
Thus the extended divisors
of $X_n$ and $X_n^\vee$ cannot be isomorphic.
In particular, by Proposition \ref{fibration.10}
its associated $\A^1$-fibration cannot be conjugated to $\pi$.
Now the proposition follows.
\eproof

\brem\label{fibration.160}
1. It is interesting
to compare our results with those in Section 5 of the recent
paper \cite{GMMR}, devoted to the study
of affine lines on smooth Gizatullin surfaces.
It turns out that, if
the Picard group $\Pic (V)$ of such a surface $V$
is not a torsion group (i.e., if $V\not\cong\A^2$),
then
there exists an affine line $\A^1\hookrightarrow V$
which is not a component
of a fiber
of any $\A^1$-fibration $V\to\A^1$.
Consequently, there is no analogue of the
Abhyankar-Moh-Suzuki Theorem for such surfaces,
and the classification of affine lines on them cannot be
deduced from that of $\A^1$-fibrations $V\to\A^1$.

2. Let $V$ be a Gizatullin surface admitting
continuous families of $\A^1$-fibrations $V\to\A^1$.
Let us deduce the existence
on
$V$ of a continuous family of
affine lines $\A^1\hookrightarrow V$ such that
for any two of them, there is
no automorphism of $V$ sending one into another.

Indeed, let $\varphi:V\to \A^1$ and $\varphi':V\to \A^1$
be two non-conjugated $\A^1$-fibrations and let $\ell$
and $\ell'$ be smooth fibers of $\varphi$ and $\varphi'$,
respectively. Let us show that these two affine
lines on $V$ are not conjugated
in the automorphism group.
Assuming the contrary, after an automorphism
we may suppose that $\ell=\ell'$
is a smooth fiber of two different $\A^1$-fibrations.
Let $(\bV, D)$ be a completion of $V$ such that $\varphi$
is the restriction of the associated standard fibration
$\Phi_0$. If $\ell'' = \varphi^{\prime -1}(a)$
is a general fiber of $\varphi'$ then  $\ell''$
cannot be a fiber component of $\varphi $.
Hence its closure $\bar\ell''$ in $\bV$ is horizontal with respect to
$\Phi_0$ and so intersects every fiber of $\Phi_0$. In particular,
$\bar\ell''$ meets both
$C_0\subseteq D$ and $\bar\ell'$.
Since $\ell'$ and $\ell''$ are disjoint,
this shows that $\bar\ell''$ meets $D$ in two different
points, which is impossible.
\erem

\subsection{ Uniqueness of
$\A^1$-fibrations on singular surfaces}
Here we consider more generally normal
Gizatullin $\C^*$-surfaces, which
are not necessary smooth. Using the technique developed in the
previous sections we are able to strengthen our previous
uniqueness result for $\A^1$-fibrations on such surfaces, see
Corollary 5.13 in \cite{FKZ3}.

 We recall the following notation from \cite{FKZ3}. Let
us consider a DPD presentation $V=\Spec\, A_0[D_+,D_-]$ of a
Gizatullin $\C^*$-surface $V$, where $A_0=\C[u]$. We let
$(\bV,D)$ be a $\C^*$-equivariant standard completion of $V$. Such
a completion is unique up to reversion of the boundary zigzag
$D=C_0\cup\ldots\cup C_n$, cf.\ (\ref{reverse}). The linear system
$|C_0|$ defines an $\A^1$-fibrations $\Phi_0:V\to\A^1$, and
similarly the linear system $|C_0^\vee|$ on the reversed
equivariant completion $(\bar V^\vee,D^\vee)$ provides a second
$\A^1$-fibration $\Phi_0^\vee:V\to\A^1$.

In \cite{FKZ3} we introduced the following two conditions:

\smallskip

($\alpha_*$) {\it $\supp \{D_+\}\cup \supp \{D_-\}$ is empty or
consists of one point  $p$, where either  $D_+(p) + D_-(p) \le -1$
or both fractional parts $\{D_+(p)\}$, $\{ D_-(p)\}$ are nonzero.}

\medskip

($\alpha_+$) {\it $\supp \{D_+\}\cup \supp \{D_-\}$ is empty or
consists of one point $p$, where either $D_+(p)+D_-(p)=0\,$  or
$$
D_+(p)+D_-(p)\le -\max\left(\frac{1}{{m^+}^2},\,
\frac{1}{{m^-}^2}\right)\,,$$
where $\pm m^\pm$ denote the
minimal positive integers such that $m^\pm D_\pm(p)\in\Z$.}

\smallskip

Corollary 5.13 in \cite{FKZ3} asserts the uniqueness of the
$\A^1$-fibration $V\to \A^1$, up to conjugation and reversion,
under condition ($\alpha_+$). In the  next proposition we
show that the latter uniqueness holds as well under the weaker
condition ($\alpha_*$).

\bprop\label{stren} Let $V$ be a Gizatullin $\C^*$-surface as
above. If condition ($\alpha_*$) is fulfilled then the
following hold.
\begin{enumerate} \item[(a)] Every $\A^1$-fibration
$V\to\A^1$ is conjugated either to $\Phi_0$ or to $\Phi_0^\vee$.
\item[(b)] Suppose that $V$ is non-toric. Then $\Phi_0$ and
$\Phi_0^\vee$ are conjugated if and only if $\{D_+(p)\}=\{
D_-(p)\}$.
\end{enumerate}\eprop

\bsit\label{prob} The proof of (b) is the same as in
\cite[Corollary 5.13]{FKZ3}. To deduce (a) we start with some
preliminary observations. Comparing ($\alpha_*$) and
($\alpha_+$) it is enough to suppose that $\supp \{D_+\}=\supp
\{D_-\}=\{p\}$ and \be\label{(*)}
0<-(D_++D_-)(p)<\max\left(\frac{1}{{m^+}^2},
\frac{1}{{m^-}^2}\right)<1\,. \ee Indeed, assuming (\ref{(*)}),
condition ($\alpha_*$) is fulfilled, but ($\alpha_+$) fails. We
precede the proof with some necessary preliminaries.

Let $D_\ext$ and $D_\ext^\vee $ denote the extended divisors of
$(\bV,D)$ and $(\bV^\vee,D^\vee)$, respectively. According to
Corollary 3.26 in \cite{FKZ3}, in our case at least one of these
divisors is rigid\footnote{See the terminology in \cite{FKZ3}.}.
We may assume that $D_\ext$ is. By Proposition 3.10 in
\cite{FKZ3}, its dual graph is

\vskip0.3truecm \be\label{ezigzag1} D_{\rm ext}:
\qquad\quad\cou{C_0}{0}\lin \cou{C_1}{0}\vlin{18}
\boxou{}{\{D_+(p)\}^*} \vlin{18} \cou{C_s\quad\;\;}\,
\cou{}{w_s}\nlin\xbshiftup{}{\qquad \{\fF_\rho\}_{\rho \ge
1}}\vlin{20} \boxou{}{\{D_-(p)\}} \nlin\boxshiftup{}{ \fF_0}
\quad,\quad \ee where $\{\fF_\rho\}_{\rho\ge 1}$ is a  collection
of $A_{k_\rho}$-feathers. Since $D_\ext$ is rigid, by Proposition
2.15 in \cite{FKZ3} the bridge curve $\tO_p^-$ of the feather
$\fF_0$ is a $(-1)$-curve. We recall (see Definition 3.20 in {\it
loc.cit.}) that the tail $L=L_{s+1}$ of $D_{\rm ext}$ is the
linear chain \be\label{tail} \boxo{\{D_-(p)\}} \vlin{17}\boxo{
\fF_0} \qquad\mbox{=} \qquad \co{C_{s+1}} \lin \ldots \lin
\co{C_n} \lin\co{F_0}\lin \ldots\lin\co{F_k}\quad, \ee where the
feather $\fF_0$ is formed by the curves $F_0,\ldots, F_k$ with
bridge curve $F_0=\tO_p^-$, cf. Remark \ref{2801}. By our
assumption, $F_0^2=-1$. According to Lemmas 3.21 and 3.22(c) in
{\em loc.cit.} we have \bnum \item $L$ is not contractible; \item
some subtail $L_t$ of $L$, $s+2\le t\le n$ can be contracted to a
smooth point, where $L_t$ is the chain starting with $C_t$ to the
right in (\ref{tail}). \enum
Let $C_r$ be the mother component of
$F_k$. By (2) $r\ge t$. Furthermore, since the subtail $L_t$ is
contractible, after contracting the divisor $L_r\setminus F_k$ and
all $(-1)$-feathers $\fF_\rho$ ($\rho\ge 1$) we obtain the
linear chain
\be\label{gr2} \cou{0}{C_0}\lin\cou{0}{C_1}
\lin\ldots\lin\cou{-1}{C_s}\lin\ldots\llin
\cou{\le-3}{C_{t-1}}\llin\cou{-2}{C_t}
\lin\cou{-2}{C_{t+1}}\lin\ldots\lin\cou{-2}{C_r}\lin\cou{-1}{F_k}\quad.
\ee
Here $C_{t-1}^2\le -3$ since by our minimality assumption
the subtail $L_{t-1}$ is not contractible. Moreover, in this chain
necessarily $C_s^2=-1$, since otherwise it cannot be contracted to
$[[0,0,0]]$.

The mother component of $C_{t-1}$  is $C_2$, since otherwise,
blowing down successively $(-1)$-vertices in (\ref{gr2}), we
arrive at a chain $[[0,0,w_2, \ldots, w_l,-1]]$ with $l\ge 2$,
$w_i\le -2 \,\,\forall i$, which cannot be contracted to
$[[0,0,0]]$.

Therefore (\ref{gr2}) can be obtained starting from the chain
$C_0\cup C_1\cup C_2$ on the quadric $Q$ and performing a sequence
of outer blowups, which create the component $C_{t-1}, \ldots,
C_r$, $F_k$, followed by a sequence of inner blowups to create the
components $C_3,\ldots,C_{t-2}$. Thus we can obtain the
completion $\bV$ from $Q$ via the following 3 steps.
\smallskip

\noindent {\bf Step 1.} Performing a sequence of outer blowups, we
create first the components $C_{t-1},C_t,\ldots,C_{r}, F_k$
getting a surface $\bV_{(1)}$ together with a linear chain
\be\label{gr4}
\cou{C_0}{0}\lin\cou{C_1}{0}\lin\cou{C_2}{-1}\lin\cou{C_{t-1}}{-2}
\lin\ldots\lin\cou{C_r}{-2}\lin\cou{F_k}{-1} \quad.
\ee
Applying on $Q$ a sequence of suitable shifts $h_{a,m}$ as in
Lemma \ref{2.12} we can subsequently move the centers
of outer blowups in $\bV_{(1)}\to Q$ into the
fixed points of the torus action on
$\bV_{(1)}$ induced, step by step, by the standard $\T$-action on
$Q$.
\smallskip

\noindent {\bf Step 2.} Starting from $\bV_{(1)}$ we create the
components $C_3,\ldots,C_{t}$ and also $C_{r+1},\ldots,C_{n}$,
$F_1,\ldots,F_{k-1}$ by a sequence of inner blowups, which results
in a surface  $\bV_{(2)}$ and a $\T$-equivariant morphisms
$\bV_{(2)}\to\bV_{(1)}\to Q$. We note that the $\T$-action on
$\bV_{(2)}$ restricts to a non-trivial $\T$-action on the
component $C_s$.
\smallskip

\noindent {\bf Step 3.} Starting from $\bV_{(2)}$ we create the
$A_{k_\rho}$-feathers $\fF_\rho$, $\rho\ge 1$, attached to the
component $C_s$. This requires outer blowups, one at each point
where a feather is attached. The remaining blowups are inner,
resulting in the completion $\bV$.\esit

\noindent {\it Proof of Proposition \ref{stren}.} Suppose that we
are given an $\A^1$-fibration $\varphi:V\to\A^1$.
By Proposition \ref{fibration.10} there
exists a standard completion $(\bV',D')$ of $V$ such that
$\varphi=\Phi_0'|V$ is defined by the linear system $|C_0'|$ on
$\bV'$. After
replacing $(\bV,D)$, if necessary, by the reversed completion we may
assume that $(\bV',D')$ is obtained from $(\bV,D)$ via a symmetric
reconstruction, see Lemma \ref{equivariant.6}.

Let us now look at the blowup process as in \ref{1.2} which
creates $\bV$ and $\bV'$, respectively, starting from the quadric
$Q=\PP^1\times\PP^1$. It suffices to prove the following claim.

\smallskip

\noindent {\it Claim. There is an automorphism $\alpha\in\Aut
(Q\setminus C_0)\cong \Aut (\A^1\times \PP^1)$ which
stabilizes the curves $C_1=\PP^1\times\{\infty\}$,
$C_2=\{0\}\times\PP^1$, fixes the point $(0,0)$ and maps the
centers of successive blowups in the decomposition of the standard
morphism $\Phi:\bV\to Q$ into the respective centers of blowups in
the decomposition of $\Phi':\bV'\to Q$.}

\smallskip

\noindent
Assuming the claim, such an automorphism $\alpha$
preserves the $\A^1$-fibration associated to the first projection
pr$_1:Q\to\PP^1$. Hence it can be lifted to an automorphism of $V$
which conjugates the standard $\A^1$-fibrations $\Phi_0: V\to
\A^1$ and $\varphi=\Phi_0': V\to \A^1$, as required.

Thus it remains to prove the claim.

\smallskip

\no {\it Proof of the claim.} We decompose the standard morphisms
$\Phi,\Phi'$ following Steps 1-3:
$$\Phi:\bV\to\bV_{(2)}\to\bV_{(1)}\to Q\qquad\mbox{ and }\qquad\Phi':
\bV'\to\bV_{(2)}'\to\bV_{(1)}'\to Q\,.$$ On Step 1 the
$\T$-equivariant morphisms $\bV_{(1)}\to Q$ and $\bV'_{(1)}\to Q$
coincide, since they have the same centers of blowups at the fixed
points of the $\T$-action. Hence  $\bV_{(1)}=\bV_{(1)}'\,.$ Then
also $\bV_{(2)}=\bV_{(2)}'\,,$ since Step 2 involves only inner
blowups.
 We let further
$$P_1,\ldots,P_k\in C_s,\qquad P_1',\ldots,P_k'\in C_s'=C_s$$
denote the base points of the $A_{k_\rho}$-feathers collections
$\{\fF_{s,\rho}\}_{\rho\ge 1}$ and $\{\fF_{s',\rho'}\}_{\rho\ge
1}$ attached to the components $C_s\subseteq \bV_{(2)}$ and
$C_s'\subseteq \bV_{(2)}'= \bV_{(2)}$, respectively. Since
$D_\ext$ is rigid, by  Proposition 5.3 in \cite{FKZ3} the
configurations
$$
\{P_1,\ldots,P_k\}
\qquad\mbox{and}\qquad\{P_1',\ldots,P_k'\}
$$
of points in $C_s\setminus (C_{s-1}\cup C_{s+1})\cong\C^*$ must be
equivalent under the $\C^*$-action on $C_s\setminus (C_{s-1}\cup
C_{s+1})$. Hence, changing suitably the enumeration and using the
induced $\T$-action on $C_s$ in $\bV_{(2)}$ we can achieve that
$P_i=P_i'$ for all $i=1,\ldots,k$ (cf.\ the proof of Proposition  \ref{main}).
Now the
surfaces $\bV$ and $\bV'$ become isomorphic via an isomorphism
which conjugates the induced $\A^1$-fibrations on $V$, as
required. \qed

\end{document}